\documentclass[11pt, a4paper]{article}
\usepackage{amsmath}
\usepackage{amsfonts}
\usepackage{amssymb}
\usepackage{ulem}

\usepackage{amsthm}
\usepackage{url}
\usepackage{dcolumn}

\usepackage{cite}
\usepackage{fancyhdr}

\usepackage{graphicx, color}
\usepackage{hyperref}
\usepackage{xcolor}

\newtheorem*{rem*}{Remark}

\renewcommand{\thefigure}{\arabic{section}.\arabic{figure}}

\makeatletter
\@addtoreset{figure}{section}
\@addtoreset{table}{section}
\@addtoreset{table}{subsection}

\usepackage{here}

\begin{document}

\title{Converting tilings with multiple types of rhombuses to pentagonal tilings}
\author{ Teruhisa SUGIMOTO$^{ 1), 2)}$ }
\date{}
\maketitle

{\footnotesize

\begin{center}
$^{1)}$ The Interdisciplinary Institute of Science, Technology and Art

$^{2)}$ Japan Tessellation Design Association

E-mail: ismsugi@gmail.com
\end{center}

}

{\small
\begin{abstract}
\noindent
The results involving rotationally symmetric tilings with multiple types of 
rhombuses, discovered by Penrose, Ammann, Beenker, or Socolar, are 
converted to tilings with multiple types of pentagons are presented. 
The pentagons can be convex or concave, and can be degenerated 
into a trapezoid or parallelogram. If the pentagons are convex, they 
belong to the Type 2 family.
\end{abstract}
}



\section{Introduction}
\label{section1}

Pentagonal tiles\footnote{ 
A \textit{tiling} (or \textit{tessellation}) of the plane is a collection of sets 
that are called tiles, which covers a plane without gaps and overlaps, except 
for the boundaries of the tiles. The term ``tile'' refers to a topological disk, 
whose boundary is a simple closed curve. If all the tiles in a tiling are of the 
same size and shape, then the tiling is \textit{monohedral}~\cite{G_and_S_1987, 
wiki_pentagon_tiling}. In this study, a polygon that admits a monohedral tiling 
is called a \textit{polygonal tile} ~\cite{Sugimoto_NoteTP, Sugimoto_2017_2}. 
Note that, in monohedral tiling, it admits the use of use reflected tiles.
} 
that can form rotationally symmetric tiling in \cite{Sugimoto_2020_1} have 
been introduced, and regarding the tilings formed by this pentagonal tile, 
two adjacent pentagons in the tiling can be converted to a rhombus; 
hence, they are closely related to rhombic monohedral tiling.

The rhombus is a convex quadrilateral tile. Several types of tilings that 
use multiple types of rhombuses (i.e. tilings, that are not monohedral) are 
also known. As introduced in Section 6 of \cite{Sugimoto_2020_1}, tilings 
with multiple types of rhombuses can be converted into those with multiple 
types of pentagons, with the aid of pentagons that can be converted into 
rhombuses; however, not all the tilings allow this conversion. In this study, we 
present the results that indicate that a rotationally symmetric tiling with 
multiple types of rhombuses is converted to a pentagonal tiling (tiling 
formed by pentagons). Specifically, we present tilings with two types 
of pentagons based on five-fold or eight-fold rotationally symmetric tilings 
with the aid of two types of rhombuses, and tilings with three types 
of pentagons based on 12-fold rotationally symmetric tiling with the 
aid of three types of rhombuses.


\section{Conditions of pentagon related to rhombic tiling}
\label{section2}

Note that the content of this section is almost the same as that in Section 2 
of \cite{Sugimoto_2020_1}.

\begin{figure}[t]
 \centering\includegraphics[width=11.5cm,clip]{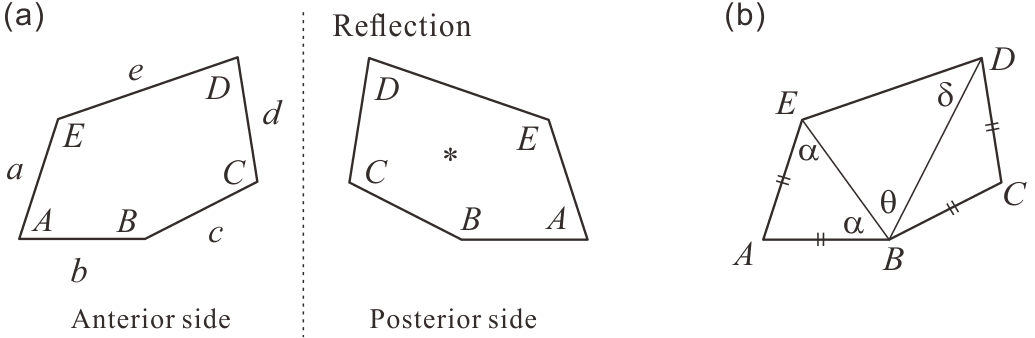} 
  \caption{{\small 
Nomenclature for vertices and edges of convex pentagon, and 
three triangles in the convex pentagonal tile C20-T2
} 
\label{Fig.2-1}
}
\end{figure}

In this study, the vertices (interior angles) and edges of the pentagon will be 
referred to using the nomenclature shown in Figure~\ref{Fig.2-1}(a). As introduced 
in \cite{Sugimoto_2020_1}, if a pentagon related to a rhombic tiling (tiling formed 
by rhombuses) is convex, it is a convex pentagonal tile belonging to the Type 2 
family.\footnote{
To date, fifteen families of convex pentagonal tiles, each of them referred 
to as a ``Type,'' are known~\cite{G_and_S_1987, Sugimoto_NoteTP, 
Sugimoto_2017_2, wiki_pentagon_tiling}. For example, if the sum of three 
consecutive angles in a convex pentagonal tile is $360^ \circ $, the pentagonal 
tile belongs to the Type 1 family. Convex pentagonal tiles belonging to some 
families also exist. Known convex pentagonal tiles can form periodic tiling. 
In May 2017, Micha\"{e}l Rao declared that the complete 
list of Types of convex pentagonal tiles had been obtained (i.e., they have 
only the known 15 families), but it does not seem to be fixed as of January 
2022~\cite{wiki_pentagon_tiling}.
} 
In \cite{S_and_O_2009}, we consider edge-to-edge\footnote{ 
A tiling by convex polygons is \textit{edge-to-edge} if any two convex polygons 
in a tiling are either disjoint or share one vertex or an entire edge in common. 
Then other case is \textit{non-edge-to-edge}~\cite{G_and_S_1987, 
Sugimoto_NoteTP, Sugimoto_2017_2}.
} 
tilings with a convex pentagon and named the convex pentagonal tile 
``C20-T2.'' C20-T2 shown in \cite{S_and_O_2009} is a convex pentagon 
that satisfies the conditions

\renewcommand{\figurename}{{\small Figure.}}

\begin{equation}
\label{eq1}
\left\{ {\begin{array}{l}
 B + D + E = 360^ \circ , \\ 
 a = b = c = d, \\ 
 \end{array}} \right.
\end{equation}

\noindent
and can form the representative tiling (tiling of edge-to-edge version) of Type 2 
that has the relations ``$B + D + E = 360^ \circ ,\;2A + 2C = 360^ \circ $.'' 
Because this convex pentagon has four equal-length edges, it can be divided 
into an isosceles triangle \textit{BCD}, an isosceles triangle \textit{ABE} with 
a base angle $\alpha $, and a triangle \textit{BDE} with $\angle DBE = \theta$ 
and $\angle BDE = \delta $, as shown in Figure~\ref{Fig.2-1}(b). Accordingly, 
using the relational expression for the interior angle of each vertex of C20-T2, 
the conditional expressions of (\ref{eq1}) can be rewritten as follows:

\begin{equation}
\label{eq2}
\left\{ {\begin{array}{l}
 A = 180^ \circ - 2\alpha , \\ 
 B = 90^ \circ + \theta , \\ 
 C = 2\alpha , \\ 
 D = 90^ \circ - \alpha + \delta , \\ 
 E = 180^ \circ + \alpha - \theta - \delta , \\ 
 a = b = c = d, \\ 
 \end{array}} \right.
\end{equation}

\noindent
where

\[
\delta = \tan ^{ - 1}\left( {\frac{\sin \theta }{\tan \alpha - \cos \theta 
}} \right)
\]

\noindent
and $0^ \circ < \alpha < 90^ \circ $ because $A > 0^ \circ$ and 
$C > 0^ \circ $~\cite{S_and_O_2009}. 
This pentagon has two degrees of freedom ($\alpha $ and $\theta $ 
parameters), besides its size. If the edge $e$ of this pentagon exists and 
the pentagon is convex, then $0^ \circ < \theta < 90^ \circ $. But 
depending on the value of $\alpha $, even if $\theta $ is selected in 
$\left( {0^ \circ ,\;90^ \circ } \right)$, the pentagon may not be convex 
or may be geometrically nonexistent. If $a = b = c = d = 1$, then the 
length of edge $e$ can be expressed as follows:

\[
e = 2\sqrt {1 - \sin (2\alpha )\cos \theta } .
\]

Let the interior angle of vertex $A$ be $\frac{360^ \circ }{n}$ (i.e., 
$\alpha = 90^  \circ-\frac{180^ \circ }{n}$) so that pentagons 
satisfying (\ref{eq2}) can form an $n$-fold rotationally 
symmetric tiling. (Remark: Due to the properties of the pentagons, 
the interior angle of vertex $C$, and not vertex $A$, will be 
$\frac{360^ \circ }{n}$.) Note that $n$ is an integer greater than or 
equal to three, because $C > 0^ \circ $. Therefore, the conditions 
of pentagonal tiles that can form $n$-fold rotationally symmetric 
tilings are expressed in (\ref{eq3}).

\begin{equation}
\label{eq3}
\left\{ {\begin{array}{l}
 A = \dfrac{360^ \circ }{n\strut}, \\ 
 B = 90^ \circ + \theta , \\ 
 C = 180^ \circ - \dfrac{360^ \circ }{n\strut}, \\ 
 D = \delta + \dfrac{180^ \circ }{n\strut}, \\ 
 E = 270^ \circ  - \theta - \delta - \dfrac{180^ \circ }{n\strut}, \\ 
 a = b = c = d. \\ 
 \end{array}} \right.
\end{equation}


\section{Relationships between pentagon and rhombus}
\label{section3}

Note that the content of this section is almost the same as that in Section 
3 of \cite{Sugimoto_2020_1}.

\begin{figure}[H]
 \centering\includegraphics[width=15cm,clip]{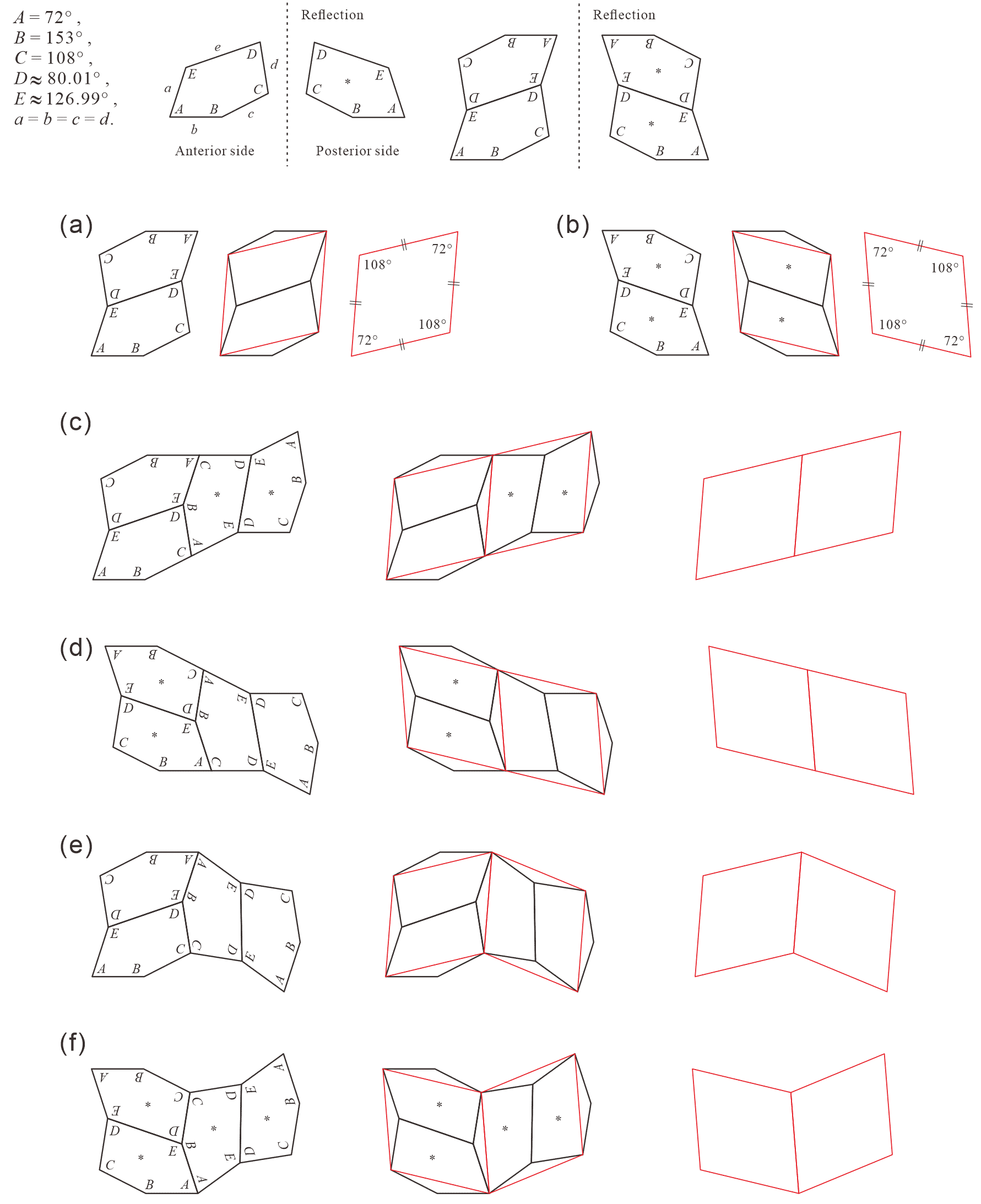} 
  \caption{{\small 
Relationships between pentagonal pair (basic unit) and rhombus} 
\label{Fig.3-1}
}
\end{figure}

The convex pentagon shown in Figure~\ref{Fig.2-1} satisfies (\ref{eq3}), where 
$n = 5$ and $\theta = 63^ \circ$. Note that it is equivalent to the case where 
$\alpha = 54^ \circ$ and $\theta = 63^ \circ$, in (\ref{eq2}). By using this convex 
pentagon of Figure~\ref{Fig.2-1}, the method of forming tilings with pentagons 
satisfying the conditions of (\ref{eq2}) or (\ref{eq3}) is described below. In 
accordance with the relationship between the five interior angles of the pentagon, 
the vertex concentrations that can be always used in tilings are 
``$A + C = 180^ \circ ,\;B + D + E = 360^ \circ ,\;2A + 2C = 360^ \circ $.'' 
According to (\ref{eq2}) and (\ref{eq3}), the edge $e$ of the pentagon is the sole 
edge of different length. Therefore, the edge $e$ of one convex pentagon is always 
connected in an edge-to-edge manner with the edge $e$ of another convex pentagon. 
A pentagonal pair with their respective vertices $D$ and $E$ concentrated forms 
the basic unit of the tiling. This basic unit can be made of two types: a 
(anterior side) pentagonal pair as shown in Figure~\ref{Fig.3-1}(a) and a reflected 
(posterior side) pentagonal pair as shown in Figure~\ref{Fig.3-1}(b). Four different 
types of units, as shown in Figures~\ref{Fig.3-1}(c), \ref{Fig.3-1}(d), \ref{Fig.3-1}(e), 
and \ref{Fig.3-1}(f), are obtained by combining two pentagonal pairs shown in 
Figures~\ref{Fig.3-1}(a) and \ref{Fig.3-1}(b), so that $B + D + E = 360^ \circ $ 
can be assembled.

As shown in Figures~\ref{Fig.3-1}(a) and \ref{Fig.3-1}(b), a rhombus (red line), with 
an acute angle of $72^ \circ $, formed by connecting the vertices $A$ and $C$ of 
the pentagon, is applied to each basic unit of the pentagonal pair. (Remark: In 
this example, because the interior angle of the vertex $A$ is $72^ \circ $, the 
rhombus has an acute angle of $72^ \circ $. That is, the interior angles of the 
rhombus corresponding to the pair of pentagons in Figures~\ref{Fig.3-1}(a) and 
\ref{Fig.3-1}(b) are the same as the interior angles of vertices $A$ and $C$ in 
(\ref{eq2}) and (\ref{eq3}).) Consequently, the parts of pentagons that protrude 
from the rhombus match exactly with the parts that are more dented than the 
rhombus (refer to Figures~\ref{Fig.3-1}(c), \ref{Fig.3-1}(d), \ref{Fig.3-1}(e), and 
\ref{Fig.3-1}(f)). In fact, tilings in which 
``$B + D + E = 360^ \circ ,\;2A + 2C = 360^ \circ $'' using pentagons satisfying 
(\ref{eq2}) and (\ref{eq3}) are equivalent to rhombic tilings. (Though a rhombus is 
a single entity, considering its internal pentagonal pattern, it will be considered 
as two entities.)

\bigskip
\bigskip

\renewcommand{\figurename}{{\small Figure.}}

\renewcommand{\figurename}{{\small Figure.}}
\begin{figure}[H]
 \centering\includegraphics[width=15cm,clip]{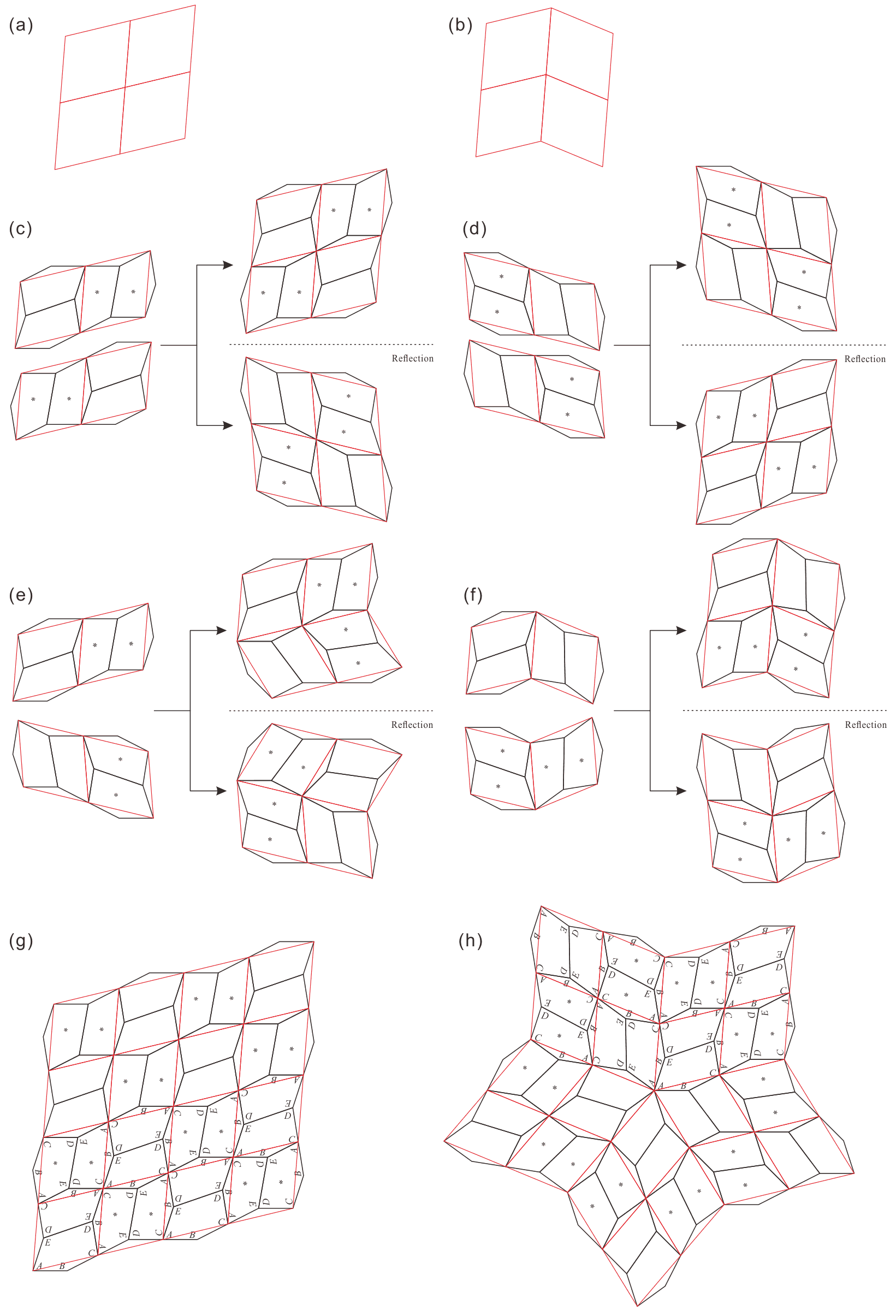} 
  \caption{{\small 
Combinations of rhombuses and pentagons} 
\label{Fig.3-2}
}
\end{figure}

Rhombuses have two-fold rotational symmetry and two axes of reflection 
symmetry passing through the center of the rotational symmetry (hereafter, 
this property is described as $D_{2}$ symmetry\footnote{
``$D_{2}$'' is based on the Schoenflies notation for symmetry in a 
two-dimensional point group~\cite{wiki_point_group, wiki_schoenflies_notation}. 
``$D_{n}$'' represents an $n$-fold rotation axis with $n$ reflection symmetry 
axes. The notation for symmetry is based on that presented in \cite{Klaassen_2016}.
}). 
Therefore, the rhombus and the reflected rhombus have identical outlines. 
Thus, the two methods of concentrating the four rhombic vertices at a point without 
gaps or overlaps are: Case (i) an arrangement by parallel translation as shown 
in Figure~\ref{Fig.3-2}(a); Case (ii) an arrangement by rotation (or reflection) 
as shown in Figure~\ref{Fig.3-2}(b). This concentration corresponds to 
forming a ``$2A + 2C = 360^ \circ $'' at the center by four pentagons. In 
Case (i), because the pentagonal vertices circulate as ``$A \to C \to A \to 
C$'' at the central ``$2A + 2C = 360^ \circ $,'' one combination (refer to 
Figure~\ref{Fig.3-2}(c)) is obtained by using two units of Figure~\ref{Fig.3-1}(c) and 
another combination (refer to Figure~\ref{Fig.3-2}(d)) is obtained by using two units 
of Figure~\ref{Fig.3-1}(d). In Case (ii), because the pentagonal vertices circulate as 
``$A \to A \to C \to C$'' at the central ``$2A + 2C = 360^ \circ $,'' one 
combination (refer to Figure~\ref{Fig.3-2}(e)) is obtained by using units of 
Figures~\ref{Fig.3-1}(c) and \ref{Fig.3-1}(d), and another combination (refer to 
Figure~\ref{Fig.3-2}(f)) is obtained by using units of Figures~\ref{Fig.3-1}(e) 
and \ref{Fig.3-1}(f).

Only when the unit comprising eight pentagons in Figures~\ref{Fig.3-2}(c) or 
\ref{Fig.3-2}(d) are arranged in a parallel manner, a tiling, as shown in 
Figure~\ref{Fig.3-2}(g), is formed that represents a tiling of Type 2, in which 
``$B + D + E = 360^ \circ ,\;2A + 2C = 360^ \circ $.'' Because rhombuses 
can form rhombic belts by translation in the same direction vertically, rhombic 
tilings can also be formed by the belts that are freely connected horizontally 
by the connecting method shown in Figures~\ref{Fig.3-2}(a) and \ref{Fig.3-2}(b). 
Further, pentagonal tilings corresponding to those rhombic tilings can be formed.

When $n$ vertices, with interior angles of $\frac{360^ \circ}{n}$, of $n$  
rhombuses are concentrated at a point, an $n$-fold rotationally symmetric 
arrangement is formed, with adjacent rhombuses connected as shown in 
Figure~\ref{Fig.3-2}(b). Therefore, an $n$-fold rotationally symmetric tiling 
with rhombuses can be formed by dividing each rhombus, in that arrangement, 
into similar shapes. By converting the rhombuses of such rhombic tiling into 
pentagons satisfying (\ref{eq3}), the rotationally symmetric tilings with convex 
pentagons can be obtained (refer to Figure~\ref{Fig.3-2}(h)). Therefore, 
when forming $n$-fold rotationally symmetric tilings from a pentagon 
satisfying (\ref{eq3}), the pentagonal arrangement can be known from the 
corresponding $n$-fold rotationally symmetric tilings with a rhombus.


\section{Conversion method of tilings with multiple types of rhombuses}
\label{section4}


\subsection{Preparation}
\label{subsection4.1}

Countless tilings with multiple types of rhombuses exist, and several rotationally 
symmetric tilings with multiple types of rhombuses are known as well. By using 
pentagons that satisfy (\ref{eq2}) or (\ref{eq3}), the values of $\theta$ in 
(\ref{eq2}) and (\ref{eq3}) are important when a tiling with multiple types of 
rhombuses is converted into one with multiple types of pentagons, which is simply 
referred to as ``pentagonal tiling." The reason is because each pentagon must 
have the same value of $\theta$ (vertex $B$) to be ``$B + D + E = 360^ \circ$" 
to hold within that pentagonal tiling.

As mentioned in \cite{Sugimoto_2020_1}, pentagons that satisfy (\ref{eq2}) or 
(\ref{eq3}) can be concave, convex, trapezoidal, or parallelogram depending 
on the parameters. Therefore, when a tiling with multiple types of rhombuses 
is converted into a pentagonal tiling, the pentagons in the tiling could be solely 
convex or concave, possess both a concave and convex shape simultaneously, 
etc. The pentagons that satisfy (\ref{eq2}) with two parameters, $\alpha$ 
and $\theta$, are difficult to treat as an example because it leaves the 
value of $\alpha $ without a clear guideline after determining the value of 
$\theta$, which is an important factor in generating a pentagonal tiling. 
The pentagons that satisfy (\ref{eq3}) have two parameters, $n$ and 
$\theta $, where $n$ is an integer of three or more. Specifically, a pentagon 
that satisfies (\ref{eq3}) is easy to treat as an example, because the shape 
of the pentagon can be determined only by the value of $\theta$ after 
the value of $n$ is determined. (Of course, the pentagon that satisfies 
(\ref{eq3}) is also the pentagon that satisfies (\ref{eq2}), but here it is 
distinguished for the sake of explanation.)

Therefore, in this study, we utilize rotationally symmetric tilings with pentagons 
that satisfy (\ref{eq3}) a case example where a tiling with multiple types 
of rhombuses is converted into a pentagonal tiling. Specifically, we 
introduce pentagonal tilings based on five-fold or eight-fold rotationally 
symmetric tilings with two types of rhombuses (i.e., the tilings are formed 
by two types of pentagons that satisfy (\ref{eq3})), and pentagonal tilings based on 
a 12-fold rotationally symmetric tiling with three types of rhombuses 
(i.e., the tilings are formed by three types of pentagons that satisfy (\ref{eq3})).

The properties of the shape of pentagons that satisfy (\ref{eq3}) depending on the 
values of $n$ and $\theta$ are arranged below. If the pentagons that satisfy 
(\ref{eq3}) exist geometrically, then $n > 2$, because $C > 0^ \circ$. It is 
$0^ \circ < \theta < 180^ \circ$ because the edge $e$ of the pentagon that 
satisfies (\ref{eq3}) exists and does not intersect the edge $b$. As mentioned in 
Section 4 of \cite{Sugimoto_2020_1}, the pentagons that satisfy (\ref{eq3}) with 
$n = 3$ and $n = 6$ correspond to rhombuses with an acute angle of $60^ \circ$, 
and these pentagons are opposite to each other. The shapes of the pentagons 
that satisfy (\ref{eq3}) change depending on the values of $n$ and $\theta$ 
are as follows:

\vspace{1pt}

\begin{enumerate}
\setlength{\itemindent}{0pt}
\item[]

\textbf{Case where the pentagons that satisfy (\ref{eq3}) admit $n = 4$}

\begin{itemize}
	\item $0^ \circ < \theta < 90^ \circ$ : Convex pentagons
	\item $\theta = 90^ \circ$ : Parallelograms ($B = 180^ \circ$)
	\item $90^ \circ < \theta < 180^ \circ$ : Concave pentagons ($B > 180^ \circ$)
\end{itemize}

\end{enumerate}

\vspace{1pt}

\begin{enumerate}
\setlength{\itemindent}{0pt}
\item[]

\textbf{Cases where the pentagons that satisfy (\ref{eq3}) admit $n \ge 5$}

\begin{itemize}
	\item $0^ \circ < \theta < 90^ \circ - \dfrac{360^ \circ }{n}$ : Concave pentagons 
($E > 180^ \circ$)

	\item $\theta = 90^ \circ - \dfrac{360^ \circ }{n}$ : Trapezoids ($E = 180^ \circ$)

	\item $90^ \circ - \dfrac{360^ \circ }{n} < \theta < 90^ \circ$ : Convex pentagons

	\item $\theta = 90^ \circ$ : Parallelograms ($B = 180^ \circ$)

	\item $90^ \circ < \theta < 180^ \circ$ : Concave pentagons ($B > 180^ \circ$)
\end{itemize}

\end{enumerate}

\vspace{1pt}

\begin{enumerate}
\setlength{\itemindent}{0pt}
\item[]

\textbf{Cases where the pentagons that satisfy (\ref{eq3}) become the polygons with 
$a = b = c = d = e$}

\begin{itemize}
	\item Case where $n = 4$ and $\theta \approx 41.41^ \circ $ (Convex pentagon with 
$E \approx 114.30^ \circ$)

	\item Case where $n = 5$ and $\theta \approx 37.95^ \circ $ (Convex pentagon with 
$E \approx 149.76^ \circ$)

	\item Case where $n = 6$ and $\theta = 30^ \circ $ (Trapezoid of $E = 180^ \circ$)

	\item Case where $n = 7$ and $\theta \approx 16.41^ \circ $ (Concave pentagon with 
$E \approx 213.69^ \circ$)
\end{itemize}

\end{enumerate}

If the pentagons that satisfy (\ref{eq3}) have $\theta = 90^ \circ $ (i.e., 
$B = 180^ \circ$), then they are parallelograms (the case of $n = 4$ is a 
rectangle). Such parallelograms correspond to half of the rhombus 
corresponding to the basic unit of Figures~\ref{Fig.3-1}(a) or \ref{Fig.3-1}(b).

If the pentagons that satisfy (\ref{eq3}) have $E = 180^ \circ $, then 
``$A = D,\;B = C$.'' Thus, they are trapezoids with a line of symmetry and 
$\theta = 90^ \circ - \frac{360^ \circ }{n}$ holds. It is when $n \ge 5$ that the 
pentagons that satisfy (\ref{eq3}) become trapezoids, because $\theta > 0^ \circ $.

If the pentagons that satisfy (\ref{eq3}) have $a = b = c = d = e = 1$, then 
$1 = 2\sqrt {1 - \sin (2\alpha )\cos \theta } $ where 
$\alpha = 90^  \circ-\frac{180^ \circ }{n}$. In this case, 
$n \le 7$ for $\theta$ to exist.


\subsection{Conversion method}
\label{subsection4.2}

\renewcommand{\figurename}{{\small Figure.}}
\begin{figure}[t]
 \centering\includegraphics[width=15cm,clip]{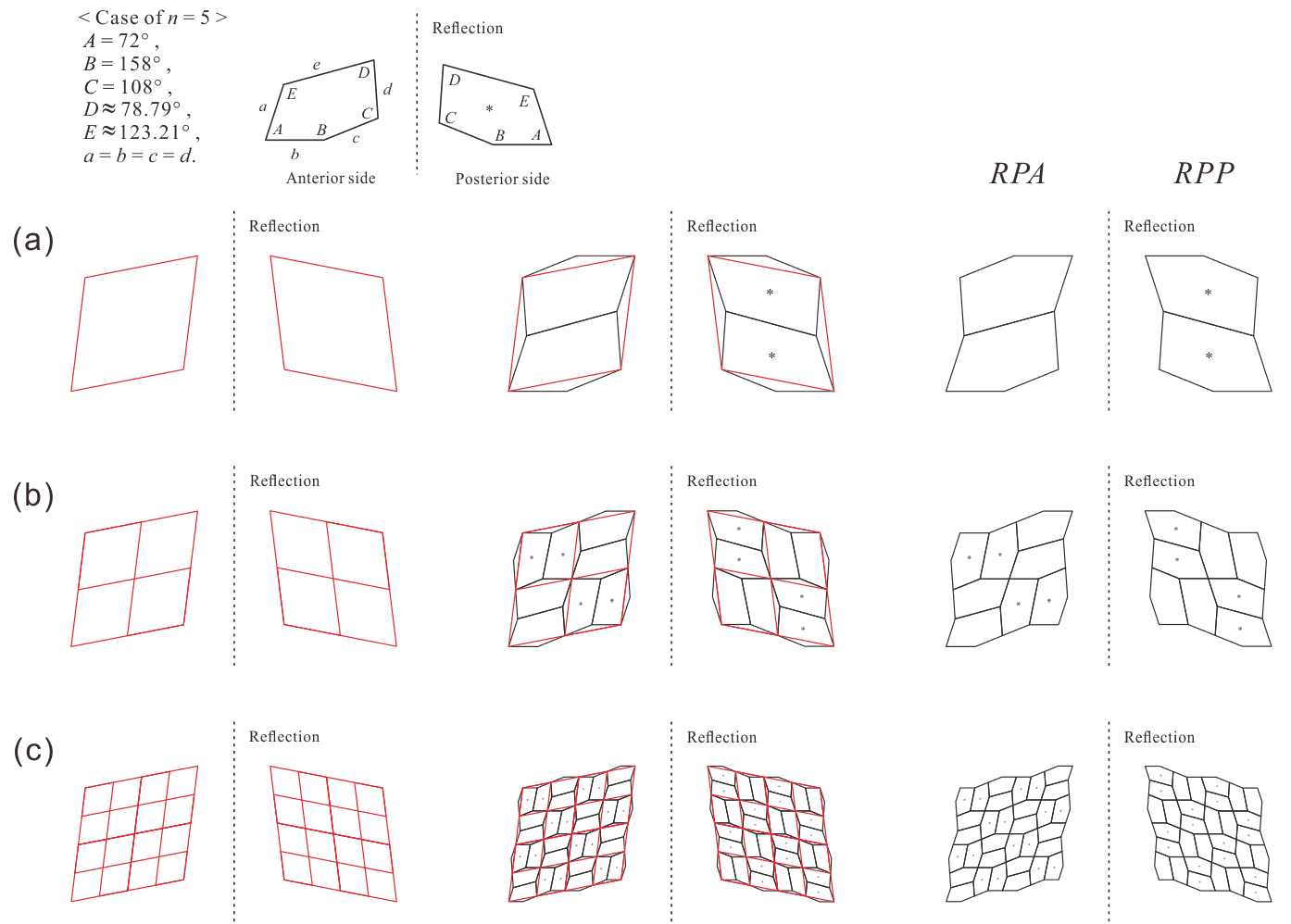} 
  \caption{{\small 
Explanation of how to divide a rhombus which can be converted 
into pentagonal tiling by its own similar figures} 
\label{Fig.4.2-1}
}
\end{figure}

By using pentagons that satisfy (\ref{eq2}) or (\ref{eq3}), when a tiling with 
multiple types of rhombuses is converted into a pentagonal tiling, there 
are two possible methods. The first method is to convert each rhombus in 
the tiling into two pentagons, as illustrated in Figure~\ref{Fig.4.2-1}(a). 
The second method is to divide each rhombus in the tiling by its own similar 
figures and convert each rhombus formed by division into two pentagons. If a 
tiling with rhombuses is monohedral, the pentagonal tiling generated by 
conversion of the above methods is the same, but if it is a tiling with 
multiple types of rhombuses, the pentagonal tiling generated by the 
conversion differs depending on how the rhombuses are divided.

In the tiling with multiple types of rhombuses, the method of dividing a 
rhombus into its own similar figures has two rules that must be followed to 
be able to convert it into a pentagonal tiling. First, it is necessary to divide the 
base rhombus so that it contains $4\times u^2$ for $u = 1, 2, 3,\ldots$ of 
its own similar rhombuses. Accordingly, when converted to a pentagon, the 
base rhombus contains $8\times u^2$ pentagons. To form a pentagonal tiling 
as illustrated in Figures~\ref{Fig.3-2}(g), \ref{Fig.3-2}(h), etc., we combine basic 
units of the pentagonal pair with the anterior and posterior sides as illustrated 
in Figures~\ref{Fig.3-1}(a) and \ref{Fig.3-1}(b). Specifically, when dividing the base 
rhombus into its own similar figures, the four divisions in Figure~\ref{Fig.4.2-1}(b) 
are the standard units of division. Because the standard unit of division 
contains eight pentagons, with four pentagons each on the anterior and 
posterior sides, it is considered as ``divide a rhombus into eight (congruent) 
pentagons'' or ``convert a rhombus into eight (congruent) pentagons.''
Therefore, dividing a rhombus by its own similar figures, which can be converted 
into a pentagonal tiling, is allowed as follows: four divisions with one standard unit 
of division (refer to Figure~\ref{Fig.4.2-1}(b)), 16 divisions with four standard units 
of division arranged in two rows and two columns (refer to Figure~\ref{Fig.4.2-1}(c)), 
36 divisions with nine standard units of division arranged in three rows and 
three columns, etc. The second rule is that each rhombus in the tiling must 
have the same number of divisions. For example, in the case where a tiling 
with two types of rhombuses, $H$ and $R$, is converted into a pentagonal 
tiling, if $H$ is divided into 16 by its own similar figures, $R$ is also divided 
into 16 by its own similar figures (i.e., unifying all rhombuses $H$ and $R$ 
in the tiling into 16 divisions).


\section{Conversion results}
\label{section5}


\subsection{Five-fold rotationally symmetric tilings}
\label{subsection5.1}

\begin{figure}[tb]
 \centering\includegraphics[width=8cm,clip]{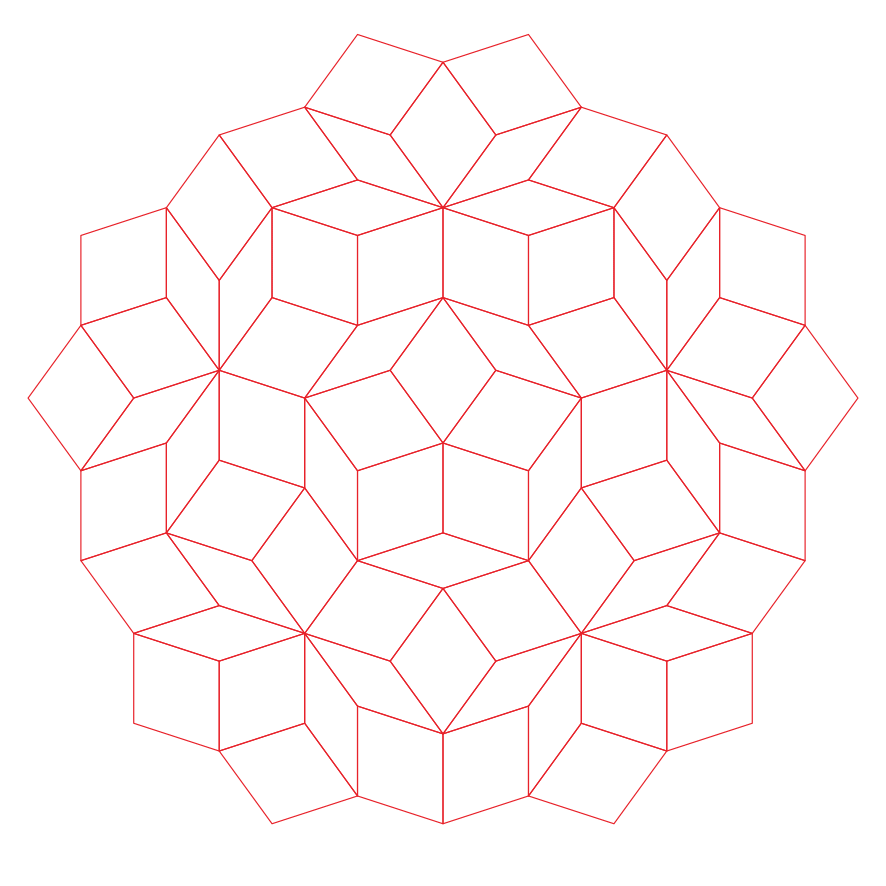} 
  \caption{{\small 
Five-fold rotationally symmetric tiling by rhombuses with 
acute angles $72^ \circ $ and $36^ \circ $ (The figure is solely a depiction 
of the area around the rotationally symmetric structure, and the tiling can 
be spread in all directions) } 
\label{Fig.5.1-1}
}
\end{figure}

In this subsection, we present pentagonal tilings based on the five-fold rotationally 
symmetric tiling with rhombuses of acute angles $72^ \circ $ and $36^ \circ$ 
discovered by Penrose, as illustrated in Figure~\ref{Fig.5.1-1}~\cite{F_H_F_2021}. 
For this conversion, two types of pentagons corresponding to $n = 5, 10$ in 
(\ref{eq3}) are utilized. When the pentagons that satisfy (\ref{eq3}) are degenerated 
into trapezoids, the trapezoids are axisymmetric; hence, they have no distinction 
between their anterior and posterior sides. A shape of two trapezoids based on 
rhombuses, as illustrated in Figures~\ref{Fig.3-1}(a) and \ref{Fig.3-1}(b), has two-fold 
rotational symmetry but no axis of reflection symmetry (refer to Figure~\ref{Fig.5.1.2-1}, 
etc.). Specifically, the shape has a distinction between the anterior and posterior sides. 
Therefore, in the figures of tilings in this study, the trapezoids that 
correspond to the posterior side arrangement are given an asterisk mark ``*,''
indicating that they are for the posterior side.

When combining two types of pentagons corresponding to $n = 5, 10$ in (\ref{eq3}), 
their shapes in a tiling change depending on the value of $\theta$ are as 
follows:

\begin{itemize}

\item $0^ \circ < \theta < 18^ \circ$: 
$n = 5, 10$ are both concave pentagons with $E > 180^ \circ$

\item $\theta = 18^ \circ$: 
$n = 5$ is a trapezoid ($E = 180^ \circ$) and $n = 10$ is a concave pentagon 
with $E > 180^ \circ $

\item $18^ \circ < \theta < 54^ \circ$: 
$n = 5$ is a convex pentagon and $n = 10$ is a concave pentagon with $E > 180^ \circ$

\item $\theta = 54^ \circ$: 
$n = 5$ is a convex pentagon and $n = 10$ is a trapezoid ($E = 180^ \circ$)

\item $54^ \circ < \theta < 90^ \circ$: 
$n = 5, 10$ are both convex pentagons

\item $\theta = 90^ \circ$: 
$n = 5, 10$ are both parallelograms ($B = 180^ \circ$)

\item $90^ \circ < \theta < 180^ \circ$:
$n = 5, 10$ are both concave pentagons with $B > 180^ \circ $

\end{itemize}

\renewcommand{\figurename}{{\small Figure.}}

\vspace{0.5ex}

\subsubsection{$0^ \circ < \theta < 18^ \circ $: $n = 5, 10$ are both concave 
pentagons with $E > 180^ \circ $}
\label{subsubsection5.1.1}

In this study, as examples, the tilings in the case of $\theta = 9^ \circ $ 
are drawn. The interior angles of pentagons satisfying (\ref{eq3}) with 
$n = 5, 10$ and $\theta = 9^ \circ$ are the values in Table~\ref{tab5.1.1}.

\begin{table}[H]
\begin{center}
{\small
\caption{
Values of interior angles of pentagons satisfying (\ref{eq3}) with 
$n = 5, 10$ and $\theta = 9^ \circ $
}
\label{tab5.1.1}
}

\begin{tabular}
{wc{55pt}|wc{55pt}wc{55pt}wc{55pt}wc{55pt}wc{55pt}}
\hline
$n$& 
$A$& 
$B$& 
$C$& 
$D$& 
$E$ \\
\hline

5& 
$72^\circ${ }& 
$99^\circ${ }& 
$108^\circ${ }& 
$57.92^\circ${ }& 
$203.08^\circ${ } \\
\hline

10& 
$36^\circ${ }& 
$99^\circ${ }& 
$144^\circ${ }& 
$22.28^\circ${ }& 
$238.72^\circ${ } \\
\hline

\end{tabular}

\end{center}
\end{table}

\vspace{-1.5ex}

\noindent
Corresponding tiling figures:

\begin{itemize}

\item[$\triangleright$]  Case converting a rhombus to two pentagons: 
Figure~\ref{Fig.5.1.1-1}

\item[$\triangleright$]  Case converting a rhombus divided into four by its own 
similar figures to eight pentagons: Figure~\ref{Fig.5.1.1-2}

\end{itemize}

\vspace{0.5ex}

\subsubsection{
$\theta = 18^ \circ $: $n = 5$ is a trapezoid ($E = 180^ \circ$) and 
$n = 10$ is a concave pentagon with $E > 180^ \circ $
}
\label{subsubsection5.1.2}

The interior angles of pentagons satisfying (\ref{eq3}) with $n = 5, 10$ 
and $\theta = 18^ \circ $ are the values in Table~\ref{tab5.1.2}.

\begin{table}[H]
\begin{center}
{\small
\caption{
Values of interior angles of pentagons satisfying (\ref{eq3}) with 
$n = 5, 10$ and $\theta = 18^ \circ $
}
\label{tab5.1.2}
}

\begin{tabular}
{wc{55pt}|wc{55pt}wc{55pt}wc{55pt}wc{55pt}wc{55pt}}
\hline
$n$& 
$A$& 
$B$& 
$C$& 
$D$& 
$E$ \\
\hline

5& 
$72^\circ${ }& 
$108^\circ${ }& 
$108^\circ${ }& 
$72^\circ${ }& 
$180^\circ${ } \\
\hline

10& 
$36^\circ${ }& 
$108^\circ${ }& 
$144^\circ${ }& 
$26.27^\circ${ }& 
$225.73^\circ${ } \\
\hline

\end{tabular}

\end{center}
\end{table}

\vspace{-1.5ex}

\noindent
Corresponding tiling figures (Note that the line corresponding to edge $e$ of 
pentagons degenerated into trapezoids in the figures is shown as a blue line):

\begin{itemize}

\item[$\triangleright$]  Case converting a rhombus to two pentagons: 
Figure~\ref{Fig.5.1.2-1}

\item[$\triangleright$]  Case converting a rhombus divided into four by its own 
similar figures to eight pentagons: Figure~\ref{Fig.5.1.2-2}

\end{itemize}

\vspace{0.5ex}

\subsubsection{
$18^ \circ < \theta < 54^ \circ $: $n = 5$ is a convex pentagon and 
$n = 10$ is a concave pentagon with $E > 180^ \circ $
}
\label{subsubsection5.1.3}

In this study, as examples, the tilings in the case of $\theta = 35^ \circ $ 
are drawn. The interior angles of pentagons satisfying (\ref{eq3}) with $n = 5, 10$ 
and $\theta = 35^ \circ $ are the values in Table~\ref{tab5.1.3}.

\begin{table}[H]
\begin{center}
{\small
\caption{
Values of interior angles of pentagons satisfying (\ref{eq3}) with 
$n = 5, 10$ and $\theta = 35^ \circ $
}
\label{tab5.1.3}
}

\begin{tabular}
{wc{55pt}|wc{55pt}wc{55pt}wc{55pt}wc{55pt}wc{55pt}}
\hline
$n$& 
$A$& 
$B$& 
$C$& 
$D$& 
$E$ \\
\hline

5& 
$72^\circ${ }& 
$125^\circ${ }& 
$108^\circ${ }& 
$81.83^\circ${ }& 
$153.17^\circ${ } \\
\hline

10& 
$36^\circ${ }& 
$125^\circ${ }& 
$144^\circ${ }& 
$32.25^\circ${ }& 
$202.75^\circ${ } \\
\hline

\end{tabular}

\end{center}
\end{table}

\vspace{-1.5ex}
\noindent
Corresponding tiling figures:

\begin{itemize}

\item[$\triangleright$]  Case converting a rhombus to two pentagons: 
Figure~\ref{Fig.5.1.3-1}

\item[$\triangleright$]  Case converting a rhombus divided into four by its own 
similar figures to eight pentagons: Figure~\ref{Fig.5.1.3-2}

\end{itemize}

\vspace{0.5ex}

\subsubsection{
$\theta = 54^ \circ $: $n = 5$ is a convex pentagon and $n = 10$ is a 
trapezoid ($E = 180^ \circ$)
}
\label{subsubsection5.1.4}

The interior angles of pentagons satisfying (\ref{eq3}) with $n = 5, 10$ and 
$\theta = 54^ \circ$ are the values in Table~\ref{tab5.1.4}.

\begin{table}[H]
\begin{center}
{\small
\caption{
Values of interior angles of pentagons satisfying (\ref{eq3}) with 
$n = 5, 10$ and $\theta = 54^ \circ $
}
\label{tab5.1.4}
}

\begin{tabular}
{wc{55pt}|wc{55pt}wc{55pt}wc{55pt}wc{55pt}wc{55pt}}
\hline
$n$& 
$A$& 
$B$& 
$C$& 
$D$& 
$E$ \\
\hline

5& 
$72^\circ${ }& 
$144^\circ${ }& 
$108^\circ${ }& 
$81.73^\circ${ }& 
$134.27^\circ${ } \\
\hline

10& 
$36^\circ${ }& 
$144^\circ${ }& 
$144^\circ${ }& 
$36^\circ${ }& 
$180^\circ${ } \\
\hline

\end{tabular}

\end{center}
\end{table}

\vspace{-1.5ex}
\noindent
Corresponding tiling figures (Note that the line corresponding to edge $e$ of 
pentagons degenerated into trapezoids in the figures is shown as a blue line):

\begin{itemize}

\item[$\triangleright$]  Case converting a rhombus to two pentagons: 
Figure~\ref{Fig.5.1.4-1}

\item[$\triangleright$]  Case converting a rhombus divided into four by its own 
similar figures to eight pentagons: Figure~\ref{Fig.5.1.4-2}

\end{itemize}

\vspace{0.5ex}

\subsubsection{
$54^ \circ < \theta < 90^ \circ $: $n = 5, 10$ are both convex 
pentagons
}
\label{subsubsection5.1.5}

In this study, as examples, the tilings in the case of $\theta = 68^ \circ$ 
are drawn. The interior angles of pentagons satisfying (\ref{eq3}) with 
$n = 5, 10$ and $\theta = 68^ \circ$ are the values in Table~\ref{tab5.1.5}.

\begin{table}[H]
\begin{center}
{\small
\caption{
Values of interior angles of pentagons satisfying (\ref{eq3}) with 
$n = 5, 10$ and $\theta = 68^ \circ $
}
\label{tab5.1.5}
}

\begin{tabular}
{wc{55pt}|wc{55pt}wc{55pt}wc{55pt}wc{55pt}wc{55pt}}
\hline
$n$& 
$A$& 
$B$& 
$C$& 
$D$& 
$E$ \\
\hline

5& 
$72^\circ${ }& 
$158^\circ${ }& 
$108^\circ${ }& 
$78.79^\circ${ }& 
$123.21^\circ${ } \\
\hline

10& 
$36^\circ${ }& 
$158^\circ${ }& 
$144^\circ${ }& 
$36.93^\circ${ }& 
$165.07^\circ${ } \\
\hline

\end{tabular}

\end{center}
\end{table}

\vspace{-1.5ex}
\noindent
Corresponding tiling figures:

\begin{itemize}

\item[$\triangleright$]  Case converting a rhombus to two pentagons: 
Figure~\ref{Fig.5.1.5-1}

\item[$\triangleright$]  Case converting a rhombus divided into four by its own 
similar figures to eight pentagons: Figure~\ref{Fig.5.1.5-2}

\item[$\triangleright$]  Case converting a rhombus divided into 16 by its own 
similar figures to 32 pentagons: Figure~\ref{Fig.5.1.5-3}

\item[$\triangleright$]  Case converting a rhombus divided into 36 by its own 
similar figures to 72 pentagons: Figure~\ref{Fig.5.1.5-4}

\end{itemize}

\vspace{0.5ex}

\subsubsection{
$\theta = 90^ \circ $: $n = 5, 10$ are both parallelograms 
($B = 180^ \circ$)
}
\label{subsubsection5.1.6}

Both are parallelograms, and the original rhombuses are formed by connecting 
edge $e$ of the parallelograms. The interior angles of pentagons 
(parallelograms) satisfying (\ref{eq3}) with $n = 5, 10$ and $\theta = 90^ \circ$ 
are the values in Table~\ref{tab5.1.6}. The converted tiling is equal to one of the 
tilings in which the rhombus of the original rhombic tiling is bisected into a 
parallelogram. The figure is omitted.

\begin{table}[H]
\begin{center}
{\small
\caption{
Values of interior angles of pentagons satisfying (\ref{eq3}) with 
$n = 5, 10$ and $\theta = 90^ \circ $
}
\label{tab5.1.6}
}

\begin{tabular}
{wc{55pt}|wc{55pt}wc{55pt}wc{55pt}wc{55pt}wc{55pt}}
\hline
$n$& 
$A$& 
$B$& 
$C$& 
$D$& 
$E$ \\
\hline

5& 
$72^\circ${ }& 
$180^\circ${ }& 
$108^\circ${ }& 
$72^\circ${ }& 
$108^\circ${ } \\
\hline

10& 
$36^\circ${ }& 
$180^\circ${ }& 
$144^\circ${ }& 
$36^\circ${ }& 
$144^\circ${ } \\
\hline

\end{tabular}

\end{center}
\end{table}


\subsubsection{
 $90^ \circ < \theta < 180^ \circ $: $n = 5, 10$ are both concave 
pentagons with $B > 180^ \circ $
}
\label{subsubsection5.1.7}

In this study, as examples, the tilings in the case of $\theta = 110^ \circ$ 
are drawn. The interior angles of pentagons satisfying (\ref{eq3}) with 
$n = 5, 10$ and $\theta = 110^ \circ$ are the values in 
Table~\ref{tab5.1.7}.

\begin{table}[H]
\begin{center}
{\small
\caption{
Values of interior angles of pentagons satisfying (\ref{eq3}) with 
$n = 5, 10$ and $\theta = 110^ \circ $
}
\label{tab5.1.7}
}

\begin{tabular}
{wc{55pt}|wc{55pt}wc{55pt}wc{55pt}wc{55pt}wc{55pt}}
\hline
$n$& 
$A$& 
$B$& 
$C$& 
$D$& 
$E$ \\
\hline

5& 
$72^\circ${ }& 
$200^\circ${ }& 
$108^\circ${ }& 
$64.67^\circ${ }& 
$95.33^\circ${ } \\
\hline

10& 
$36^\circ${ }& 
$200^\circ${ }& 
$144^\circ${ }& 
$33.36^\circ${ }& 
$126.64^\circ${ } \\
\hline

\end{tabular}

\end{center}
\end{table}

\vspace{-1.5ex}
\noindent
Corresponding tiling figures:
\begin{itemize}

\item[$\triangleright$] Case converting a rhombus to two pentagons: 
Figure~\ref{Fig.5.1.7-1}

\item[$\triangleright$]  Case converting a rhombus divided into four by its own 
similar figures to eight pentagons: Figure~\ref{Fig.5.1.7-2}

\end{itemize}

\vspace{1ex}

\subsection{Eight-fold rotationally symmetric tilings}
\label{subsection5.2}

\begin{figure}[tb]
 \centering\includegraphics[width=8cm,clip]{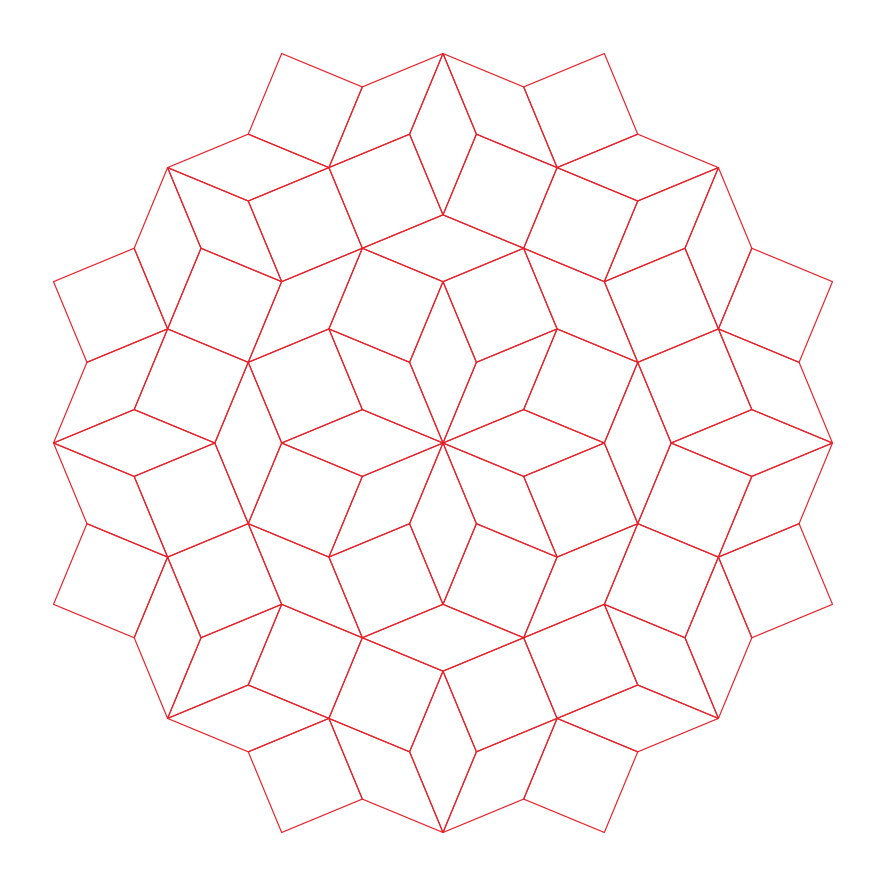} 
  \caption{{\small 
Eight-fold rotationally symmetric tiling by rhombuses with 
acute angles $90^ \circ$ and $45^ \circ$ (The figure is solely a depiction 
of the area around the rotationally symmetric structure, and the tiling can 
be spread in all directions) } 
\label{Fig.5.2-1}
}
\end{figure}

In this subsection, we present pentagonal tilings based on the eight-fold rotationally 
symmetric tiling with rhombuses of acute angles $90^ \circ $ and $45^ \circ$ discovered 
by Ammann and Beenker, as illustrated in Figure~\ref{Fig.5.2-1}~\cite{F_H_F_2021}. 
Of course, a rhombus with $90^ \circ $ angles is a square. For this conversion, 
two types of pentagons corresponding to $n = 4, 8$ in (\ref{eq3}) are utilized. 
As described in Subsection~\ref{subsection5.1}, the trapezoids corresponding to 
the posterior side in the tiling are given an asterisk mark. The pentagons 
corresponding to $n = 4$ in (\ref{eq3}) have a line of symmetry connecting 
the vertex $B$ to the midpoint of the edge $e$, i.e., there is no distinction between 
its anterior and posterior sides. A shape of square converted into two pentagons 
corresponding to $n = 4$ in (\ref{eq3}), as illustrated in Figures~\ref{Fig.3-1}(a) 
and \ref{Fig.3-1}(b), has $D_{2}$ symmetry. In other words, because the shape 
has no distinction between its anterior and posterior sides (refer to 
Figure~\ref{Fig.5.2.1-1}, etc.), a unit comprising two anterior pentagons can be freely 
replaced with one comprising two posterior pentagons. Therefore, in the 
figures of tilings in this study, the pentagons with $n = 4$ in (\ref{eq3}) that 
correspond to the posterior side are not given an asterisk mark. 

When combining two types of pentagons corresponding to $n = 4, 8$ in (\ref{eq3}), 
their shapes in a tiling change depending on the values of $\theta $ are as 
follows:

\begin{itemize}

\item $0^ \circ < \theta < 45^ \circ$: 
$n = 4$ is a convex pentagon and $n = 8$ is a concave pentagon with $E > 180^ \circ $

\item $\theta = 45^ \circ$: 
$n = 4$ is a convex pentagon and $n = 8$ is a trapezoid ($E = 180^ \circ$)

\item $45^ \circ < \theta < 90^ \circ$: 
$n = 4, 8$ are both convex pentagons

\item $\theta = 90^ \circ$: 
$n = 4, 8$ are both parallelograms ($B = 180^ \circ$)

\item $90^ \circ < \theta < 180^ \circ$: 
$n = 4, 8$ are both concave pentagons with $B > 180^ \circ $

\end{itemize}

\renewcommand{\figurename}{{\small Figure.}}

\vspace{0.5ex}

\subsubsection{
$0^ \circ < \theta < 45^ \circ $: $n = 4$ is a convex pentagon and $n = 8$ 
is a concave pentagon with $E > 180^ \circ $
}
\label{subsubsection5.2.1}

In this study, as examples, the tilings in the case of $\theta = 26^ \circ$ 
are drawn. The interior angles of pentagons satisfying (\ref{eq3}) with 
$n = 4, 8$ and $\theta = 26^ \circ$ are the values in 
Table~\ref{tab5.2.1}.

\begin{table}[H]
\begin{center}
{\small
\caption{
Values of interior angles of pentagons satisfying (\ref{eq3}) with 
$n = 4, 8$ and $\theta = 26^ \circ$
}
\label{tab5.2.1}
}

\begin{tabular}
{wc{55pt}|wc{55pt}wc{55pt}wc{55pt}wc{55pt}wc{55pt}}
\hline
$n$& 
$A$& 
$B$& 
$C$& 
$D$& 
$E$ \\
\hline

4& 
$90^\circ${ }& 
$116^\circ${ }& 
$90^\circ${ }& 
$122^\circ${ }& 
$122^\circ${ } \\
\hline

8& 
$45^\circ${ }& 
$116^\circ${ }& 
$135^\circ${ }& 
$38.63^\circ${ }& 
$205.37^\circ${ } \\
\hline

\end{tabular}

\end{center}
\end{table}

\vspace{-1.5ex}
\noindent
Corresponding tiling figures:
\begin{itemize}

\item[$\triangleright$] Case converting a rhombus to two pentagons: 
Figure~\ref{Fig.5.2.1-1}

\item[$\triangleright$]  Case converting a rhombus divided into four by its own 
similar figures to eight pentagons: Figure~\ref{Fig.5.2.1-2}

\end{itemize}

\vspace{0.5ex}

\subsubsection{
$\theta = 45^ \circ $: $n = 4$ is a convex pentagon and $n = 8$ is a 
trapezoid ($E = 180^ \circ$)
}
\label{subsubsection5.2.2}

The interior angles of pentagons satisfying (\ref{eq3}) with $n = 4, 8$ and 
$\theta = 45^ \circ$ are the values in Table~\ref{tab5.2.2}. In this case, 
because the shapes of the pentagon of $n = 4$ and trapezoid of $n = 8$ have 
no distinction between its anterior and posterior sides, it can be considered 
as a tiling that formed of only two complete pieces, the same as the original 
rhombic tiling. (Because the rhombuses have $D_{2}$ symmetry, the original 
rhombic tiling forms two complete pieces. However, if the pentagons of the 
anterior and posterior sides are not congruent when stacked on top of each 
other, we can treat them as separate pieces and consider the tiling to be 
formed of three or more different pieces.)

\begin{table}[H]
\begin{center}
{\small
\caption{
Values of interior angles of pentagons satisfying (\ref{eq3}) with 
$n = 4, 8$ and $\theta = 45^ \circ $
}
\label{tab5.2.2}
}

\begin{tabular}
{wc{55pt}|wc{55pt}wc{55pt}wc{55pt}wc{55pt}wc{55pt}}
\hline
$n$& 
$A$& 
$B$& 
$C$& 
$D$& 
$E$ \\
\hline

4& 
$90^\circ${ }& 
$135^\circ${ }& 
$90^\circ${ }& 
$112.5^\circ${ }& 
$112.5^\circ${ } \\
\hline

8& 
$45^\circ${ }& 
$135^\circ${ }& 
$135^\circ${ }& 
$45^\circ${ }& 
$180^\circ${ } \\
\hline

\end{tabular}

\end{center}
\end{table}

\vspace{-1.5ex}
\noindent
Corresponding tiling figures (Note that the line corresponding to edge $e$ of 
pentagons degenerated into trapezoids in the figures is shown as a blue line):
\begin{itemize}

\item[$\triangleright$] Case converting a rhombus to two pentagons: 
Figure~\ref{Fig.5.2.2-1}

\item[$\triangleright$]  Case converting a rhombus divided into four by its own 
similar figures to eight pentagons: Figure~\ref{Fig.5.2.2-2}

\item[$\triangleright$] Case converting a rhombus divided into 16 by its own 
similar figures to 32 pentagons: Figure~\ref{Fig.5.2.2-3}

\end{itemize}

\vspace{0.5ex}

\subsubsection{
$45^ \circ < \theta < 90^ \circ $: $n = 4, 8$ are both convex 
pentagons
}
\label{subsubsection5.2.3}

In this study, as examples, the tilings in the case of $\theta = 61^ \circ$ 
are drawn. The interior angles of pentagons satisfying (\ref{eq3}) 
with $n = 4, 8$ and $\theta = 61^ \circ$ are the values in 
Table~\ref{tab5.2.3}.

\begin{table}[H]
\begin{center}
{\small
\caption{
Values of interior angles of pentagons satisfying (\ref{eq3}) with 
$n = 4, 8$ and $\theta = 61^ \circ $
}
\label{tab5.2.3}
}

\begin{tabular}
{wc{55pt}|wc{55pt}wc{55pt}wc{55pt}wc{55pt}wc{55pt}}
\hline
$n$& 
$A$& 
$B$& 
$C$& 
$D$& 
$E$ \\
\hline

4& 
$90^\circ${ }& 
$151^\circ${ }& 
$90^\circ${ }& 
$104.5^\circ${ }& 
$104.5^\circ${ } \\
\hline

8& 
$45^\circ${ }& 
$151^\circ${ }& 
$135^\circ${ }& 
$46.89^\circ${ }& 
$162.11^\circ${ } \\
\hline

\end{tabular}

\end{center}
\end{table}

\vspace{-1.5ex}
\noindent
Corresponding tiling figures:
\begin{itemize}

\item[$\triangleright$] Case converting a rhombus to two pentagons: 
Figure~\ref{Fig.5.2.3-1}

\item[$\triangleright$]  Case converting a rhombus divided into four by its own 
similar figures to eight pentagons: Figure~\ref{Fig.5.2.3-2}

\item[$\triangleright$] Case converting a rhombus divided into 16 by its own 
similar figures to 32 pentagons: Figure~\ref{Fig.5.2.3-3}

\end{itemize}

\vspace{0.5ex}

\subsubsection{
 $\theta = 90^ \circ $: $n = 4, 8$ are both parallelograms ($B = 180^ \circ$)
}
\label{subsubsection5.2.4}

Both are parallelograms ($n = 4$ is a rectangle), and the original rhombuses 
are formed by connecting edge $e$ of the parallelograms. The interior angles of 
pentagons (parallelograms) satisfying (\ref{eq3}) with $n = 4, 8$ and 
$\theta = 90^ \circ $ are the values in Table~\ref{tab5.2.4}. The converted 
tiling is equal to one of the tilings in which the rhombus of the original rhombic 
tiling is bisected into a parallelogram. The figure is omitted.

\begin{table}[H]
\begin{center}
{\small
\caption{
Values of interior angles of pentagons satisfying (\ref{eq3}) with 
$n = 4, 8$ and $\theta = 90^ \circ $
}
\label{tab5.2.4}
}

\begin{tabular}
{wc{55pt}|wc{55pt}wc{55pt}wc{55pt}wc{55pt}wc{55pt}}
\hline
$n$& 
$A$& 
$B$& 
$C$& 
$D$& 
$E$ \\
\hline

4& 
$90^\circ${ }& 
$180^\circ${ }& 
$90^\circ${ }& 
$90^\circ${ }& 
$90^\circ${ } \\
\hline

8& 
$45^\circ${ }& 
$180^\circ${ }& 
$135^\circ${ }& 
$45^\circ${ }& 
$135^\circ${ } \\
\hline

\end{tabular}

\end{center}
\end{table}


\subsubsection{
$90^ \circ < \theta < 180^ \circ$: $n = 4, 8$ are both concave 
pentagons with $B > 180^ \circ$
}
\label{subsubsection5.2.5}

In this study, as examples, the tilings in the case of $\theta = 134^ \circ$ 
are drawn. The interior angles of pentagons satisfying (\ref{eq3}) with 
$n = 4, 8$ and $\theta = 134^ \circ$ are the values in Table~\ref{tab5.2.5}.

\begin{table}[H]
\begin{center}
{\small
\caption{
Values of interior angles of pentagons satisfying (\ref{eq3}) with 
$n = 4, 8$ and $\theta = 134^ \circ $
}
\label{tab5.2.5}
}

\begin{tabular}
{wc{55pt}|wc{55pt}wc{55pt}wc{55pt}wc{55pt}wc{55pt}}
\hline
$n$& 
$A$& 
$B$& 
$C$& 
$D$& 
$E$ \\
\hline

4& 
$90^\circ${ }& 
$224^\circ${ }& 
$90^\circ${ }& 
$68^\circ${ }& 
$68^\circ${ } \\
\hline

8& 
$45^\circ${ }& 
$224^\circ${ }& 
$135^\circ${ }& 
$35.53^\circ${ }& 
$100.47^\circ${ } \\
\hline

\end{tabular}

\end{center}
\end{table}

\vspace{-1.5ex}
\noindent
Corresponding tiling figures:
\begin{itemize}

\item[$\triangleright$] Case converting a rhombus to two pentagons: 
Figure~\ref{Fig.5.2.5-1}

\item[$\triangleright$]  Case converting a rhombus divided into four by its own 
similar figures to eight pentagons: Figure~\ref{Fig.5.2.5-2}

\end{itemize}

\vspace{1ex}

\subsection{12-fold rotationally symmetric tilings}
\label{subsection5.3}

\begin{figure}[tb]
 \centering\includegraphics[width=8cm,clip]{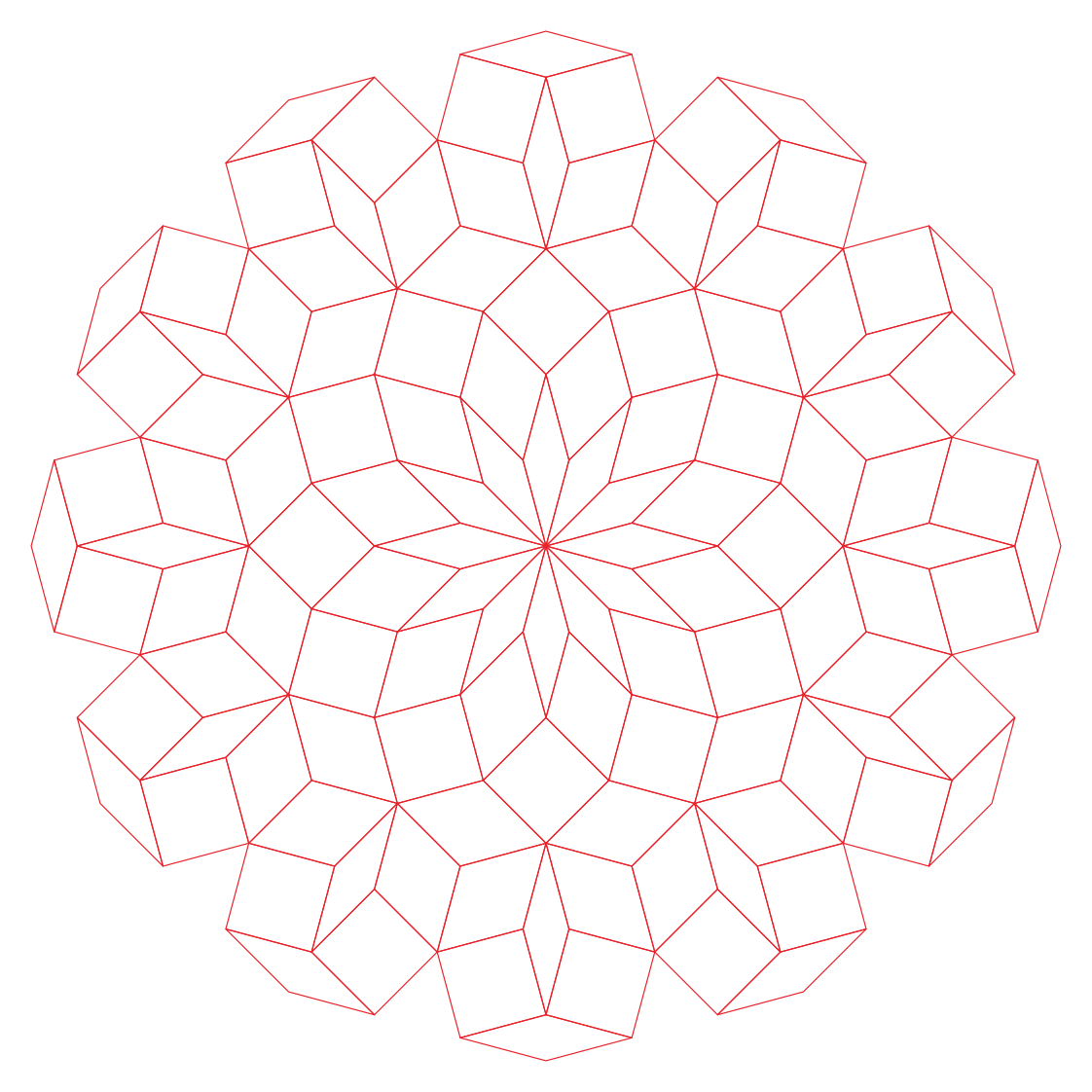} 
  \caption{{\small 
12-fold rotationally symmetric tiling by rhombuses with acute angles 
$90^ \circ $, $60^ \circ $, and $30^ \circ $ (The figure is solely a depiction 
of the area around the rotationally symmetric structure, and the tiling can 
be spread in all directions)} 
\label{Fig.5.3-1}
}
\end{figure}

In this subsection, we present pentagonal tilings based on the 12-fold rotationally 
symmetric tiling with rhombuses with of acute angles $90^ \circ$, $60^ \circ$, 
and $30^ \circ $ discovered by Socolar, as illustrated in 
Figure~\ref{Fig.5.3-1}~\cite{F_H_F_2021, Niizeki_1989}. 
Of course, a rhombus with $90^ \circ $ angles is a square. For this conversion, 
three types of pentagons corresponding to $n = 4, 6, 12$ in (\ref{eq3}) 
are utilized. As described in Subsections~\ref{subsection5.1} and \ref{subsection5.2}, 
in the figures of tilings in this study, the trapezoid corresponding to the posterior 
side is given an asterisk mark, and the pentagons with $n = 4$ in (\ref{eq3}) are not 
given an asterisk mark. 

When combining three types of pentagons corresponding to $n = 4, 6, 12$ in 
(\ref{eq3}), their shapes in a tiling change depending on the values of $\theta$ 
are as follows:

\begin{itemize}

\item $0^ \circ < \theta < 30^ \circ$: 
$n = 4$ is a convex pentagon and $n = 6, 12$ are concave pentagons 
with $E > 180^ \circ $

\item $\theta = 30^ \circ$: 
$n = 4$ is a convex pentagon, $n = 6$ is a trapezoid ($E = 180^ \circ$), 
and $n = 12$ is a concave pentagon with $E > 180^ \circ $

\item $30^ \circ < \theta < 60^ \circ$: 
$n = 4, 6$ are convex pentagons and $n = 12$ is a concave pentagon 
with $E > 180^ \circ $

\item $\theta = 60^ \circ$: 
$n = 4, 6$ are convex pentagons and $n = 12$ is a trapezoid 
($E = 180^ \circ$)

\item $60^ \circ < \theta < 90^ \circ$: 
$n = 4, 6, 12$ are all convex pentagons

\item $\theta = 90^ \circ$: 
$n = 4, 6, 12$ are all parallelograms ($B = 180^ \circ$)

\item $90^ \circ < \theta < 180^ \circ$: 
$n = 4, 6, 12$ are all concave pentagons with $B > 180^ \circ $

\end{itemize}

\renewcommand{\figurename}{{\small Figure.}}

\vspace{0.5ex}

\subsubsection{
$0^ \circ < \theta < 30^ \circ $: $n = 4$ is a convex pentagon and
 $n = 6, 12$ are concave pentagons with $E > 180^ \circ$
}
\label{subsubsection5.3.1}

In this study, as examples, the tilings in the case of $\theta = 17^ \circ$ 
are drawn. The interior angles of pentagons satisfying (\ref{eq3}) with 
$n = 4, 6, 12$ and $\theta = 17^ \circ $ are the values in 
Table~\ref{tab5.3.1}.

\begin{table}[H]
\begin{center}
{\small
\caption{
Values of interior angles of pentagons satisfying (\ref{eq3}) with 
$n = 4, 6, 12$ and $\theta = 17^ \circ $
}
\label{tab5.3.1}
}

\begin{tabular}
{wc{55pt}|wc{55pt}wc{55pt}wc{55pt}wc{55pt}wc{55pt}}
\hline
$n$& 
$A$& 
$B$& 
$C$& 
$D$& 
$E$ \\
\hline

4& 
$90^\circ${ }& 
$107^\circ${ }& 
$90^\circ${ }& 
$126.5^\circ${ }& 
$126.5^\circ${ } \\
\hline

6& 
$60^\circ${ }& 
$107^\circ${ }& 
$120^\circ${ }& 
$50.65^\circ${ }& 
$202.35^\circ${ } \\
\hline

12& 
$30^\circ${ }& 
$107^\circ${ }& 
$150^\circ${ }& 
$21.01^\circ${ }& 
$231.99^\circ${ } \\
\hline

\end{tabular}

\end{center}
\end{table}

\vspace{-1.5ex}
\noindent
Corresponding tiling figures:
\begin{itemize}

\item[$\triangleright$] Case converting a rhombus to two pentagons: 
Figure~\ref{Fig.5.3.1-1}

\item[$\triangleright$]  Case converting a rhombus divided into four by its own 
similar figures to eight pentagons: Figure~\ref{Fig.5.3.1-2}

\end{itemize}

\vspace{0.5ex}

\subsubsection{
$\theta = 30^ \circ $: $n = 4$ is a convex pentagon, $n = 6$ is a 
trapezoid ($E = 180^ \circ$), and $n = 12$ is a concave pentagon 
with $E > 180^ \circ $
}
\label{subsubsection5.3.2}

The interior angles of pentagons satisfying (\ref{eq3}) with $n = 4, 6, 12$ 
and $\theta = 30^ \circ $ are the values in Table~\ref{tab5.3.2}.

\begin{table}[H]
\begin{center}
{\small
\caption{
Values of interior angles of pentagons satisfying (\ref{eq3}) with 
$n = 4, 6, 12$ and $\theta = 30^ \circ $
}
\label{tab5.3.2}
}

\begin{tabular}
{wc{55pt}|wc{55pt}wc{55pt}wc{55pt}wc{55pt}wc{55pt}}
\hline
$n$& 
$A$& 
$B$& 
$C$& 
$D$& 
$E$ \\
\hline

4& 
$90^\circ${ }& 
$120^\circ${ }& 
$90^\circ${ }& 
$120^\circ${ }& 
$120^\circ${ } \\
\hline

6& 
$60^\circ${ }& 
$120^\circ${ }& 
$120^\circ${ }& 
$60^\circ${ }& 
$180^\circ${ } \\
\hline

12& 
$30^\circ${ }& 
$120^\circ${ }& 
$150^\circ${ }& 
$24.90^\circ${ }& 
$215.10^\circ${ } \\
\hline

\end{tabular}

\end{center}
\end{table}

\vspace{-1.5ex}
\noindent
Corresponding tiling figures (Note that the line corresponding to edge $e$ of 
pentagons degenerated into trapezoids in the figures is shown as a blue line):
\begin{itemize}

\item[$\triangleright$] Case converting a rhombus to two pentagons: 
Figure~\ref{Fig.5.3.2-1}

\item[$\triangleright$]  Case converting a rhombus divided into four by its own 
similar figures to eight pentagons: Figure~\ref{Fig.5.3.2-2}

\item[$\triangleright$] Case converting a rhombus divided into 16 by its own 
similar figures to 32 pentagons: Figure~\ref{Fig.5.3.2-3}

\end{itemize}

\vspace{0.5ex}

\subsubsection{
$30^ \circ < \theta < 60^ \circ$: $n = 4, 6$ are convex pentagons 
and $n = 12$ is a concave pentagon with $E > 180^ \circ $
}
\label{subsubsection5.3.3}

In this study, as examples, the tilings in the case of $\theta = 44^ \circ$ 
are drawn. The interior angles of pentagons satisfying (\ref{eq3}) with 
$n = 4, 6, 12$ and $\theta = 44^ \circ$ are the values in  
Table~\ref{tab5.3.3}.

\begin{table}[H]
\begin{center}
{\small
\caption{
Values of interior angles of pentagons satisfying (\ref{eq3}) with 
$n = 4, 6, 12$ and $\theta = 44^ \circ $
}
\label{tab5.3.3}
}

\begin{tabular}
{wc{55pt}|wc{55pt}wc{55pt}wc{55pt}wc{55pt}wc{55pt}}
\hline
$n$& 
$A$& 
$B$& 
$C$& 
$D$& 
$E$ \\
\hline

4& 
$90^\circ${ }& 
$134^\circ${ }& 
$90^\circ${ }& 
$113^\circ${ }& 
$113^\circ${ } \\
\hline

6& 
$60^\circ${ }& 
$134^\circ${ }& 
$120^\circ${ }& 
$64.45^\circ${ }& 
$161.55^\circ${ } \\
\hline

12& 
$30^\circ${ }& 
$134^\circ${ }& 
$150^\circ${ }& 
$27.98^\circ${ }& 
$198.02^\circ${ } \\
\hline

\end{tabular}

\end{center}
\end{table}

\vspace{-1.5ex}
\noindent
Corresponding tiling figures:
\begin{itemize}

\item[$\triangleright$] Case converting a rhombus to two pentagons: 
Figure~\ref{Fig.5.3.3-1}

\item[$\triangleright$]  Case converting a rhombus divided into four by its own 
similar figures to eight pentagons: Figure~\ref{Fig.5.3.3-2}

\end{itemize}

\vspace{0.5ex}

\subsubsection{
$\theta = 60^ \circ $: $n = 4, 6$ are convex pentagons and $n = 12$ 
is a trapezoid ($E = 180^ \circ$)
}
\label{subsubsection5.3.4}

The interior angles of pentagons satisfying (\ref{eq3}) with 
$n = 4, 6, 12$ and $\theta = 60^ \circ $ are the values in 
Table~\ref{tab5.3.4}.

\begin{table}[H]
\begin{center}
{\small
\caption{
Values of interior angles of pentagons satisfying (\ref{eq3}) with 
$n = 4, 6, 12$ and $\theta = 60^ \circ $
}
\label{tab5.3.4}
}

\begin{tabular}
{wc{55pt}|wc{55pt}wc{55pt}wc{55pt}wc{55pt}wc{55pt}}
\hline
$n$& 
$A$& 
$B$& 
$C$& 
$D$& 
$E$ \\
\hline

4& 
$90^\circ${ }& 
$150^\circ${ }& 
$90^\circ${ }& 
$105^\circ${ }& 
$105^\circ${ } \\
\hline

6& 
$60^\circ${ }& 
$150^\circ${ }& 
$120^\circ${ }& 
$65.10^\circ${ }& 
$144.90^\circ${ } \\
\hline

12& 
$30^\circ${ }& 
$150^\circ${ }& 
$150^\circ${ }& 
$30^\circ${ }& 
$180^\circ${ } \\
\hline

\end{tabular}

\end{center}
\end{table}

\vspace{-1.5ex}
\noindent
Corresponding tiling figures (Note that the line corresponding to edge $e$ of 
pentagons degenerated into trapezoids in the figures is shown as a blue line):
\begin{itemize}

\item[$\triangleright$] Case converting a rhombus to two pentagons: 
Figure~\ref{Fig.5.3.4-1}

\item[$\triangleright$]  Case converting a rhombus divided into four by its own 
similar figures to eight pentagons: Figure~\ref{Fig.5.3.4-2}

\end{itemize}

\vspace{0.5ex}

\subsubsection{
$60^ \circ < \theta < 90^ \circ $: $n = 4, 6, 12$ are all convex 
pentagons
}
\label{subsubsection5.3.5}

In this study, as examples, the tilings in the case of $\theta = 72^ \circ $ 
are drawn. The interior angles of pentagons satisfying (\ref{eq3}) with 
$n = 4, 6, 12$ and $\theta = 72^ \circ $ are the values
Table~\ref{tab5.3.5}.

\begin{table}[H]
\begin{center}
{\small
\caption{
Values of interior angles of pentagons satisfying (\ref{eq3}) with 
$n = 4, 6, 12$ and $\theta = 72^ \circ $
}
\label{tab5.3.5}
}

\begin{tabular}
{wc{55pt}|wc{55pt}wc{55pt}wc{55pt}wc{55pt}wc{55pt}}
\hline
$n$& 
$A$& 
$B$& 
$C$& 
$D$& 
$E$ \\
\hline

4& 
$90^\circ${ }& 
$162^\circ${ }& 
$90^\circ${ }& 
$99^\circ${ }& 
$99^\circ${ } \\
\hline

6& 
$60^\circ${ }& 
$162^\circ${ }& 
$120^\circ${ }& 
$63.76^\circ${ }& 
$134.24^\circ${ } \\
\hline

12& 
$30^\circ${ }& 
$162^\circ${ }& 
$150^\circ${ }& 
$30.53^\circ${ }& 
$167.47^\circ${ } \\
\hline

\end{tabular}

\end{center}
\end{table}

\vspace{-1.5ex}
\noindent
Corresponding tiling figures:
\begin{itemize}

\item[$\triangleright$] Case converting a rhombus to two pentagons: 
Figure~\ref{Fig.5.3.5-1}

\item[$\triangleright$]  Case converting a rhombus divided into four by its own 
similar figures to eight pentagons: Figure~\ref{Fig.5.3.5-2}

\item[$\triangleright$]  Case converting a rhombus divided into 16 by its own 
similar figures to 32 pentagons: Figure~\ref{Fig.5.3.5-3}

\end{itemize}

\vspace{0.5ex}

\subsubsection{
$\theta = 90^ \circ $: $n = 4, 6, 12$ are all parallelograms ($B = 180^ \circ$)
}
\label{subsubsection5.3.6}

All cases are parallelograms, and the original rhombuses are formed by 
connecting edge $e$ of the parallelograms. The interior angles of pentagons 
(parallelograms) satisfying (\ref{eq3}) with $n = 4, 6, 12$ and 
$\theta = 90^ \circ$ are the values in Table~\ref{tab5.3.6}. The converted 
tiling is equal to one of the tilings in which the rhombus of the original rhombic 
tiling is bisected into a parallelogram. The figure is omitted.

\begin{table}[H]
\begin{center}
{\small
\caption{
Values of interior angles of pentagons satisfying (\ref{eq3}) with 
$n = 4, 6, 12$ and $\theta = 90^ \circ $
}
\label{tab5.3.6}
}

\begin{tabular}
{wc{55pt}|wc{55pt}wc{55pt}wc{55pt}wc{55pt}wc{55pt}}
\hline
$n$& 
$A$& 
$B$& 
$C$& 
$D$& 
$E$ \\
\hline

4& 
$90^\circ${ }& 
$180^\circ${ }& 
$90^\circ${ }& 
$90^\circ${ }& 
$90^\circ${ } \\
\hline

6& 
$60^\circ${ }& 
$180^\circ${ }& 
$120^\circ${ }& 
$60^\circ${ }& 
$120^\circ${ } \\
\hline

12& 
$30^\circ${ }& 
$180^\circ${ }& 
$150^\circ${ }& 
$30^\circ${ }& 
$150^\circ${ } \\
\hline

\end{tabular}

\end{center}
\end{table}


\subsubsection{
$90^ \circ < \theta < 180^ \circ $: $n = 4, 6, 12$ are all concave 
pentagons with $B > 180^ \circ $
}
\label{subsubsection5.3.7}

In this study, as examples, the tilings in the case of $\theta = 125^ \circ $ 
are drawn. The interior angles of pentagons satisfying (\ref{eq3}) with 
$n = 4, 6, 12$ and $\theta = 125^ \circ $ are the values in 
Table~\ref{tab5.3.7}.

\begin{table}[H]
\begin{center}
{\small
\caption{
Values of interior angles of pentagons satisfying (\ref{eq3}) with 
$n = 4, 6, 12$ and $\theta = 125^ \circ $
}
\label{tab5.3.7}
}

\begin{tabular}
{wc{55pt}|wc{55pt}wc{55pt}wc{55pt}wc{55pt}wc{55pt}}
\hline
$n$& 
$A$& 
$B$& 
$C$& 
$D$& 
$E$ \\
\hline

4& 
$90^\circ${ }& 
$215^\circ${ }& 
$90^\circ${ }& 
$72.5^\circ${ }& 
$72.5^\circ${ } \\
\hline

6& 
$60^\circ${ }& 
$215^\circ${ }& 
$120^\circ${ }& 
$49.56^\circ${ }& 
$95.44^\circ${ } \\
\hline

12& 
$30^\circ${ }& 
$215^\circ${ }& 
$150^\circ${ }& 
$25.77^\circ${ }& 
$119.23^\circ${ } \\
\hline

\end{tabular}

\end{center}
\end{table}

\vspace{-1.5ex}
\noindent
Corresponding tiling figures:
\begin{itemize}

\item[$\triangleright$] Case converting a rhombus to two pentagons: 
Figure~\ref{Fig.5.3.7-1}

\item[$\triangleright$]  Case converting a rhombus divided into four by its own 
similar figures to eight pentagons: Figure~\ref{Fig.5.3.7-2}

\end{itemize}

\vspace{0.5ex}

\subsubsection{
Pattern generation using the properties of pentagons corresponding to 
rhombus with acute angle of $60^ \circ $
}
\label{subsubsection5.3.8}

As introduced in \cite{Sugimoto_2020_1}, the pentagons satisfying (\ref{eq3}) 
with $n = 6$, corresponding to the rhombus with an acute angle of $60^ \circ$ 
have the property of being able to form a unit that has an outline shape with 
$D_{3}$ symmetry. By using that property, the interior of the tilings can be 
replaced with different patterns (refer to Figures 14, 16, 22, and 29 in 
\cite{Sugimoto_2020_1}). For example, Figure~\ref{Fig.5.3.8-1} illustrates such 
replacements to the tilings in Figures~\ref{Fig.5.3.1-1} and \ref{Fig.5.3.5-2}. 
In particular, the pentagon satisfying (\ref{eq3}) with $n = 6$ and 
$\theta = 30^ \circ$ becomes a trapezoid, and the trapezoid has the ability 
to form equilateral triangles. Therefore, the corresponding parts in the tilings 
can be replaced by more different patterns (Figure~\ref{Fig.5.3.8-2} illustrates 
such a replacement in the tiling of Figure~\ref{Fig.5.3.2-2}).


\section{Conclusions}
\label{section6}

In \cite{Sugimoto_2020_1}, we mainly dealt with rotationally symmetric tiling 
with pentagonal tiles, and discussed the relationship between the pentagon 
and rhombus. They were mainly discussions of monohedral tiling.

In this study, we focused on the conversion of tiling with multiple types of 
rhombuses into tiling with multiple types of pentagons, and presented the 
results of the conversion in detail using rotationally symmetric tilings, 
which are also the tiling model for quasicrystals.

When a rhombus is converted into pentagons, there are two types of 
rhombuses, rhombus with acute angles corresponding to pentagons of Anterior 
side (\textit{RPA}) and rhombus with acute angles corresponding to pentagons of 
Posterior side (\textit{RPP}), and they are reflection symmetry (refer to 
Figure~\ref{Fig.4.2-1}). In the conversion that divides each rhombus of a 
rotationally symmetric tiling with multiple types of rhombuses as indicated in 
this study into two pentagons, the tiling must be done using \textit{RPA} and 
\textit{RPP}. However, if a rotationally symmetric tiling with multiple types of 
rhombuses is converted so that each rhombus is divided into eight or more 
pentagons, as indicated in this study, the tiling is formed from either \textit{RPA} 
or \textit{RPP} only. In this study, the converted pentagonal tilings were 
arranged in such a way that its center of rotational symmetry has an \textit{RPA}; 
hence, the tilings generated by a conversion that divides a rhombus into eight 
or more pentagons were formed by a group of pentagons of rhombus that 
correspond to the \textit{RPA}.

In the tiling with multiple types of rhombuses, the method of dividing a 
rhombus into its own similar figures has the rule that must be followed for 
its conversion into a pentagonal tiling. A tiling with multiple types of 
rhombuses generated by a division that does not follow the rule cannot be 
converted into a pentagonal tiling. Therefore, not all rhombic tilings can 
be converted into pentagonal tilings by the method introduced in this study. 
However, using the method described in this study, various pentagonal 
tilings will be generated from the knowledge of various patterns about 
rhombic tilings, and the patterns will be interesting.


\appendix
\def\thesection{ }


\section{}

\renewcommand{\thefigure}{\Alph{section}{-}\arabic{figure}}

\renewcommand{\figurename}{{\small Figure.}}
\begin{figure}[p]
 \centering\includegraphics[width=15cm,clip]{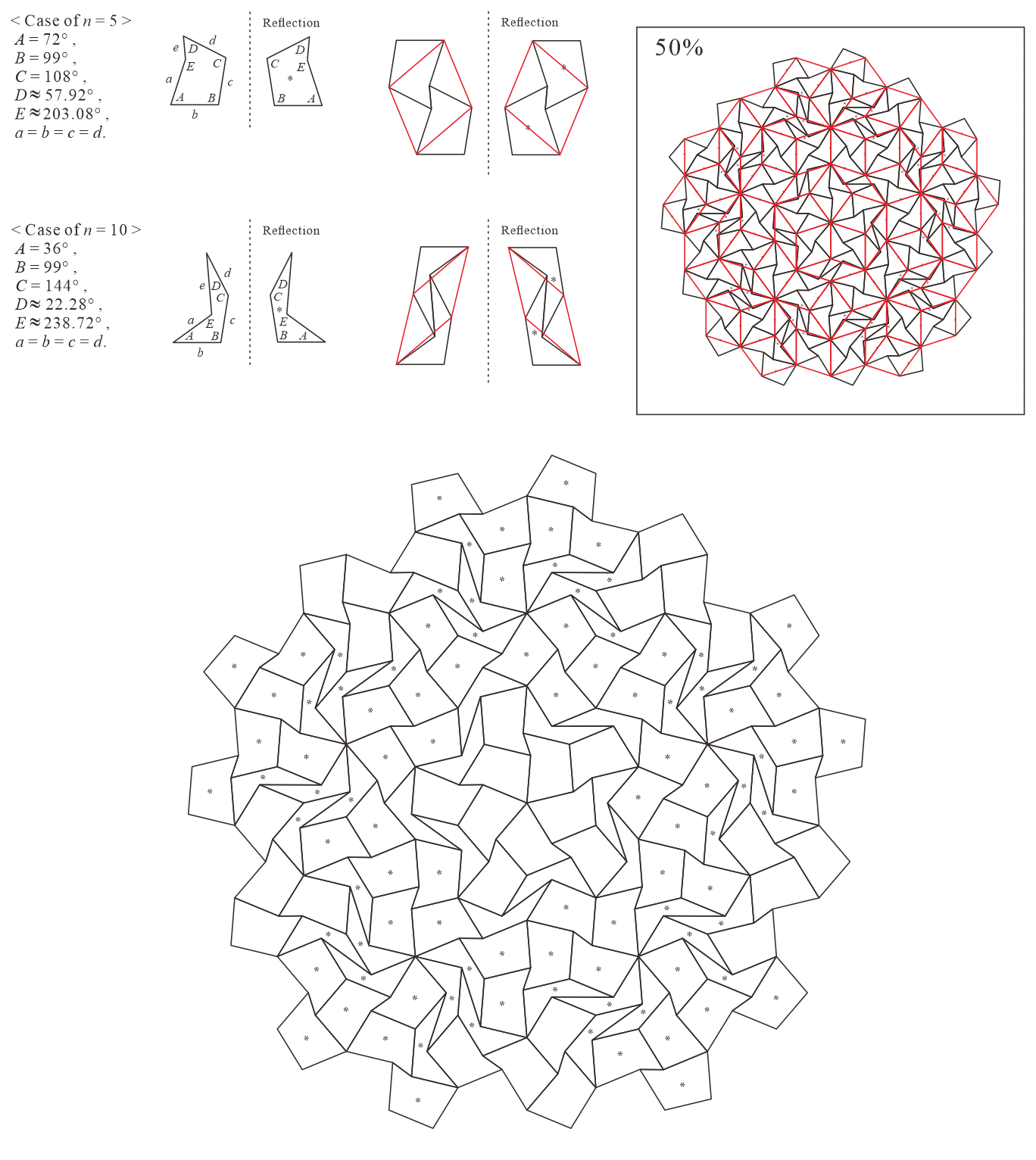} 
  \caption{{\small 
Tiling that generated based on five-fold rotationally symmetric tiling, 
with two types of rhombuses applying conversion of rhombuses into two 
pentagons satisfying (\ref{eq3}) with $n = 5, 10$ and $\theta = 9^ \circ $ 
(The figure is solely a depiction of the area around the rotationally 
symmetric structure, and the tiling can be spread in all directions)
} 
\label{Fig.5.1.1-1}
}
\end{figure}

\renewcommand{\figurename}{{\small Figure.}}
\begin{figure}[p]
 \centering\includegraphics[width=15cm,clip]{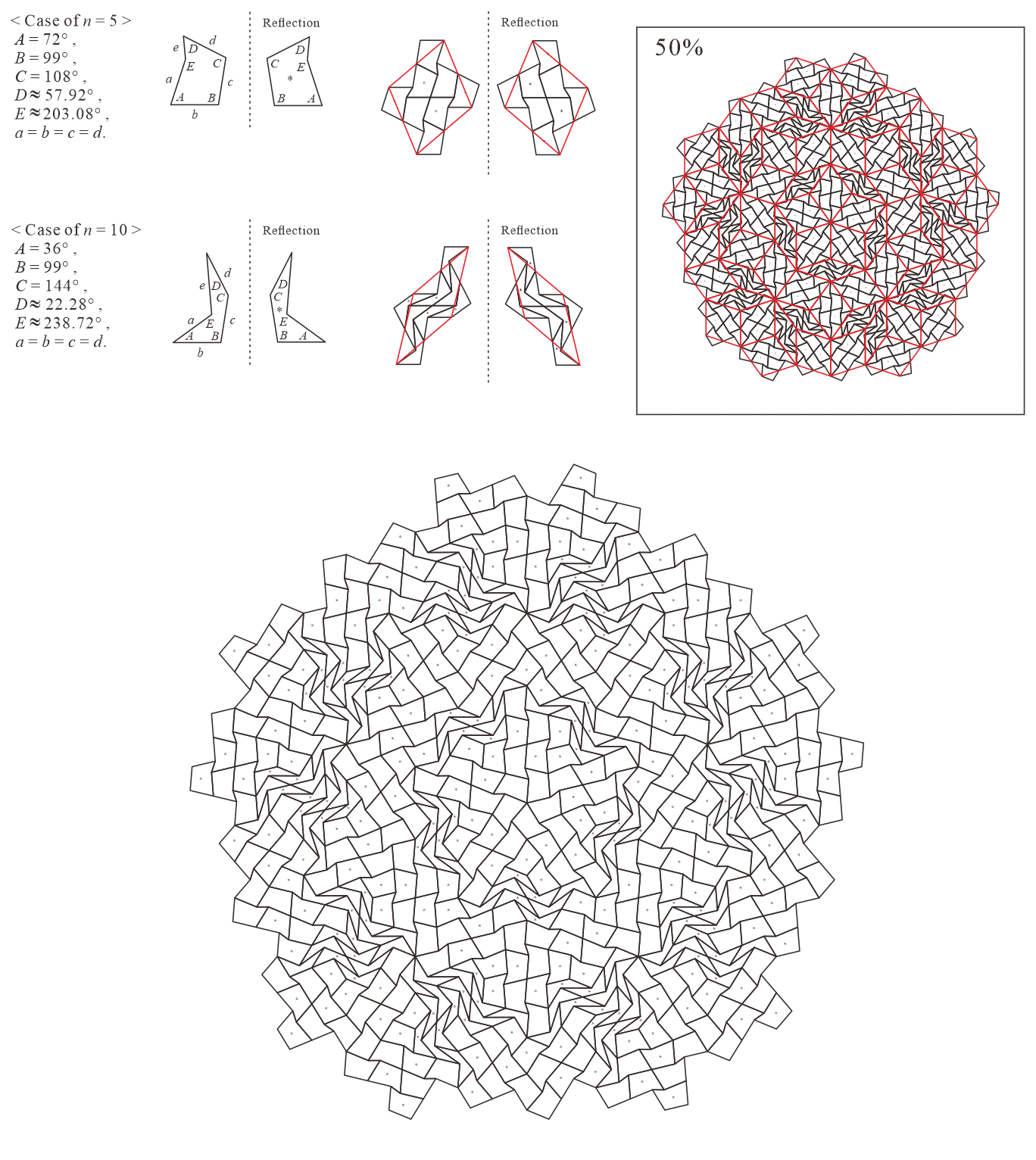} 
  \caption{{\small 
Tiling that generated based on five-fold rotationally symmetric tiling, 
with two types of rhombuses applying conversion of rhombuses into eight 
pentagons satisfying (\ref{eq3}) with $n = 5, 10$ and $\theta = 9^ \circ $ 
(The figure is solely a depiction of the area around the rotationally 
symmetric structure, and the tiling can be spread in all directions)
} 
\label{Fig.5.1.1-2}
}
\end{figure}


\renewcommand{\figurename}{{\small Figure.}}
\begin{figure}[p]
 \centering\includegraphics[width=15cm,clip]{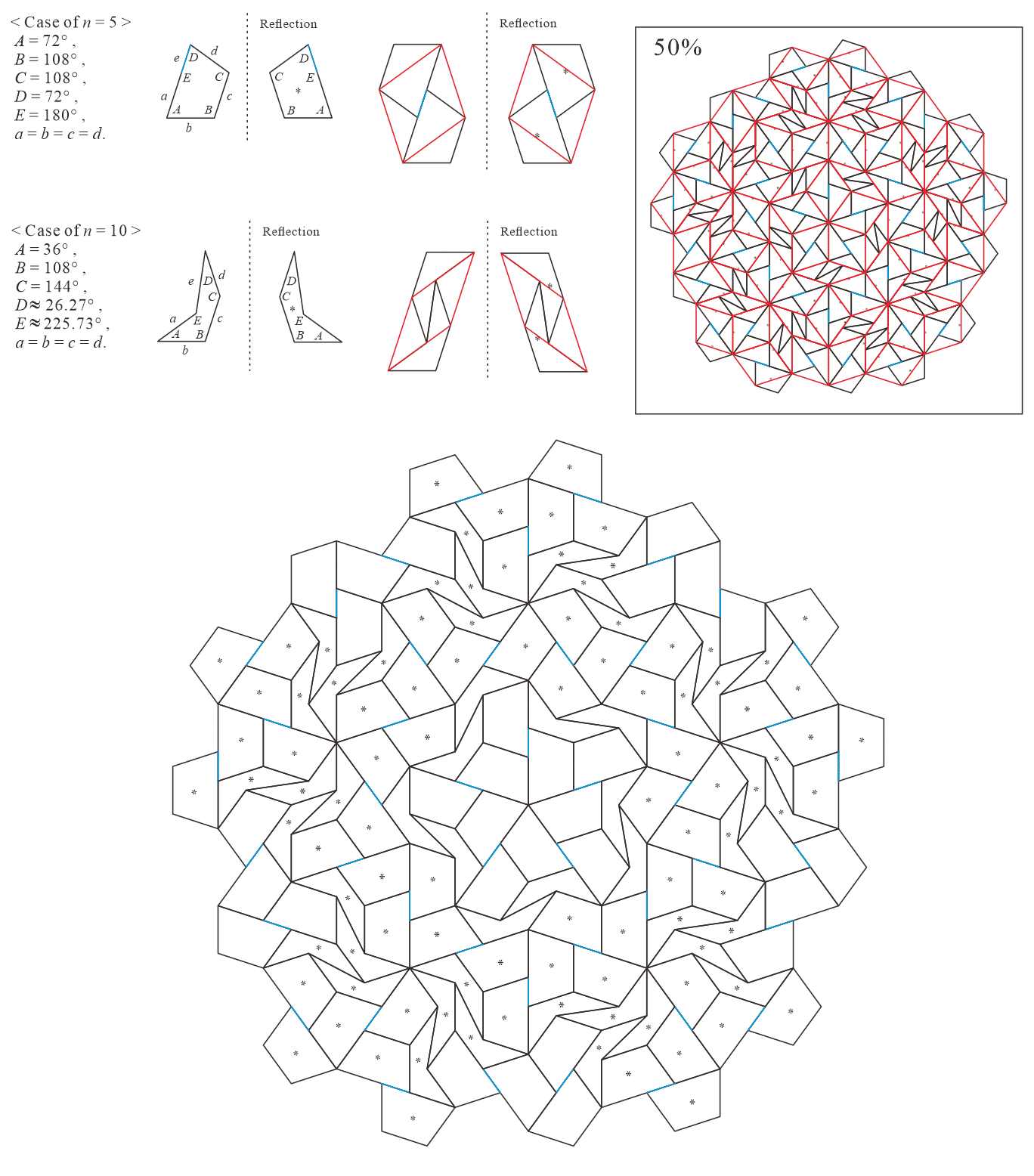} 
  \caption{{\small 
Tiling that generated based on five-fold rotationally symmetric tiling, 
with two types of rhombuses applying conversion of rhombuses into two 
pentagons satisfying (\ref{eq3}) with $n = 5, 10$ and $\theta = 18^ \circ $ 
(The figure is solely a depiction of the area around the rotationally 
symmetric structure, and the tiling can be spread in all directions)
} 
\label{Fig.5.1.2-1}
}
\end{figure}

\renewcommand{\figurename}{{\small Figure.}}
\begin{figure}[p]
 \centering\includegraphics[width=15cm,clip]{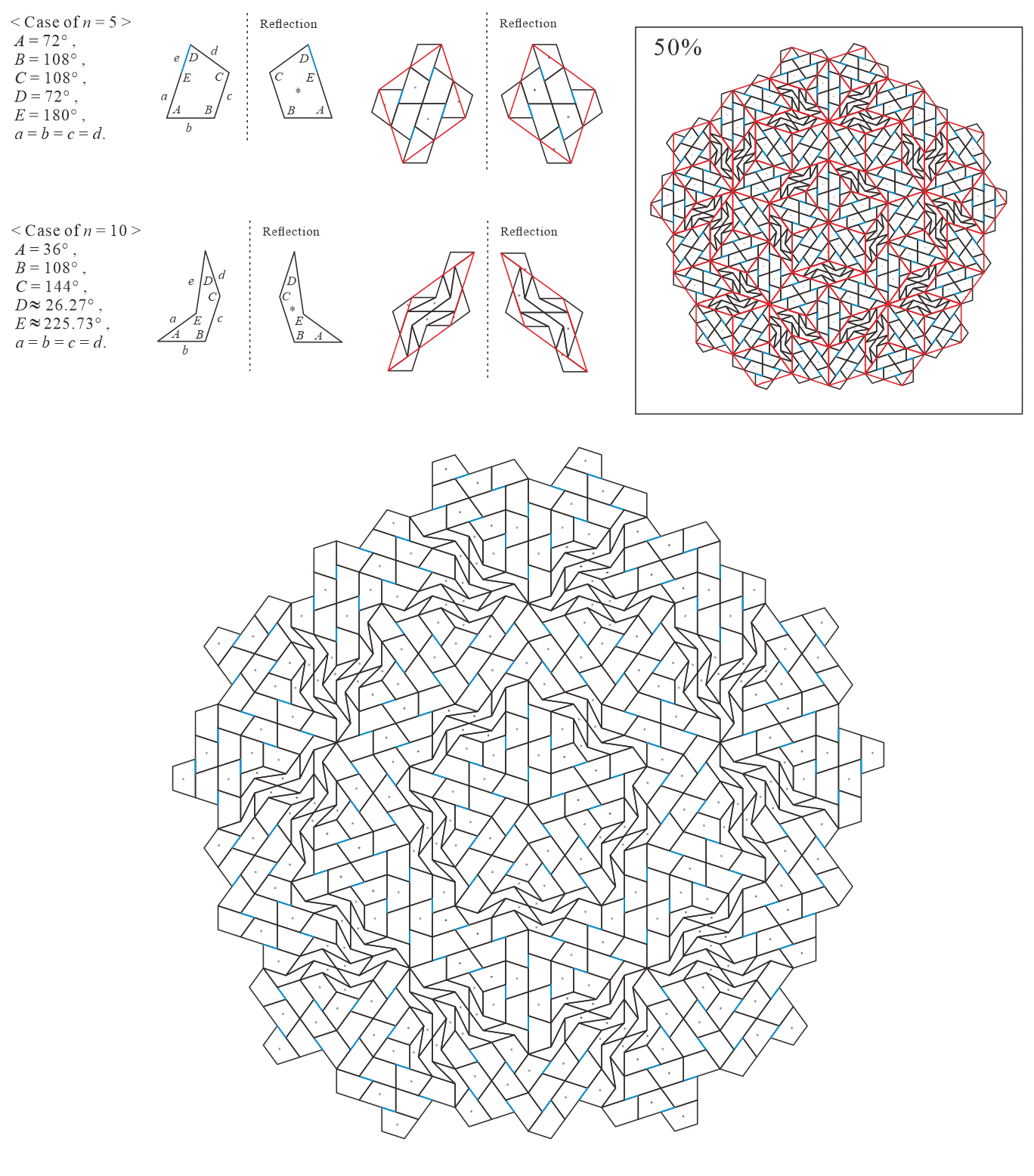} 
  \caption{{\small 
Tiling that generated based on five-fold rotationally symmetric tiling, 
with two types of rhombuses applying conversion of rhombuses into eight 
pentagons satisfying (\ref{eq3}) with $n = 5, 10$ and $\theta = 18^ \circ $ 
(The figure is solely a depiction of the area around the rotationally 
symmetric structure, and the tiling can be spread in all directions)
} 
\label{Fig.5.1.2-2}
}
\end{figure}


\renewcommand{\figurename}{{\small Figure.}}
\begin{figure}[p]
 \centering\includegraphics[width=15cm,clip]{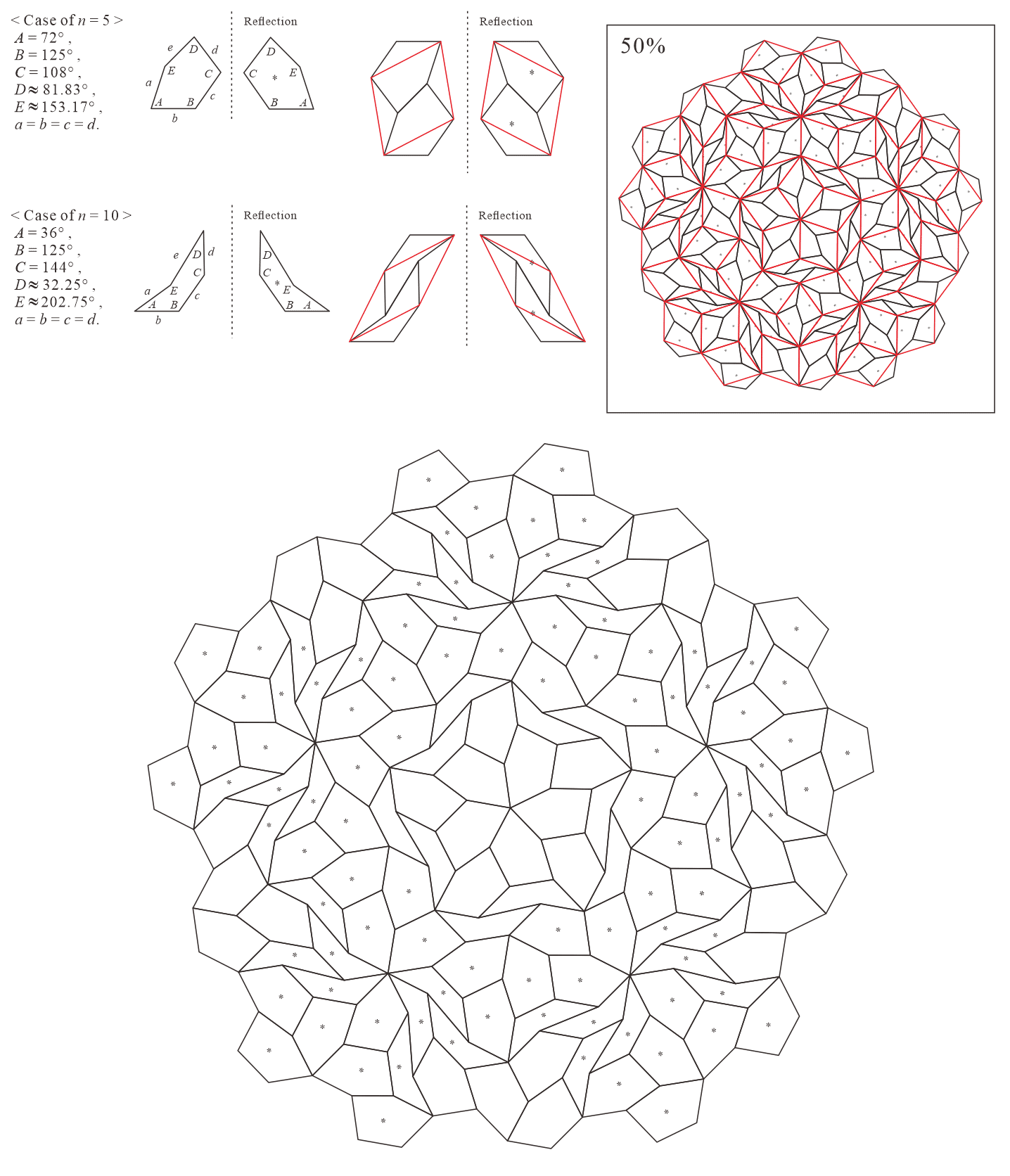} 
  \caption{{\small 
Tiling that generated based on five-fold rotationally symmetric tiling, 
with two types of rhombuses applying conversion of rhombuses into two 
pentagons satisfying (\ref{eq3}) with $n = 5, 10$ and $\theta = 35^ \circ $ 
(The figure is solely a depiction of the area around the rotationally 
symmetric structure, and the tiling can be spread in all directions)
} 
\label{Fig.5.1.3-1}
}
\end{figure}

\renewcommand{\figurename}{{\small Figure.}}
\begin{figure}[p]
 \centering\includegraphics[width=15cm,clip]{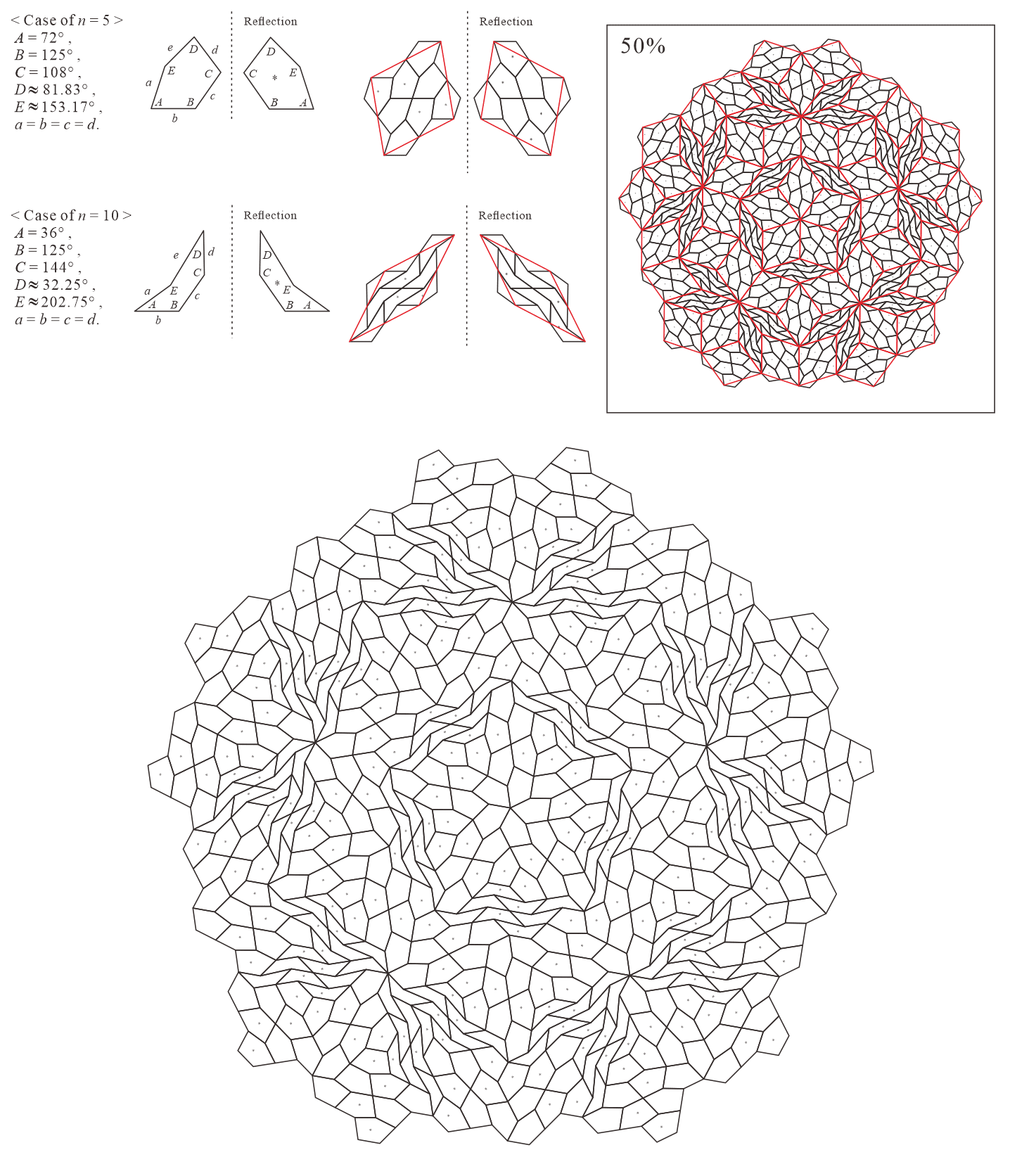} 
  \caption{{\small 
Tiling that generated based on five-fold rotationally symmetric tiling, 
with two types of rhombuses applying conversion of rhombuses into eight 
pentagons satisfying (\ref{eq3}) with $n = 5, 10$ and $\theta = 35^ \circ $ 
(The figure is solely a depiction of the area around the rotationally 
symmetric structure, and the tiling can be spread in all directions)
} 
\label{Fig.5.1.3-2}
}
\end{figure}


\renewcommand{\figurename}{{\small Figure.}}
\begin{figure}[p]
 \centering\includegraphics[width=15cm,clip]{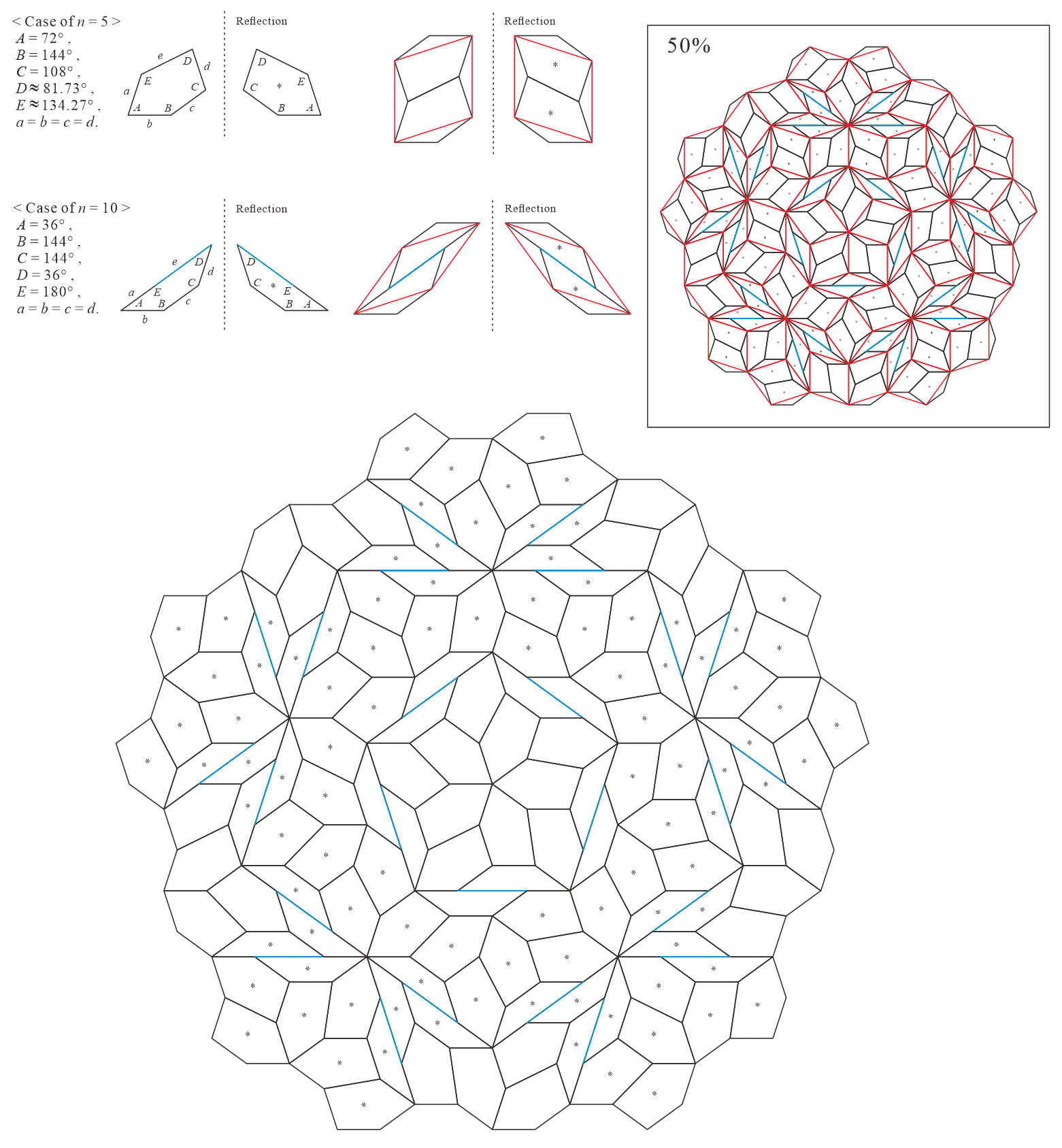} 
  \caption{{\small 
Tiling that generated based on five-fold rotationally symmetric tiling, 
with two types of rhombuses applying conversion of rhombuses into two 
pentagons satisfying (\ref{eq3}) with $n = 5, 10$ and $\theta = 54^ \circ $ 
(The figure is solely a depiction of the area around the rotationally 
symmetric structure, and the tiling can be spread in all directions)
} 
\label{Fig.5.1.4-1}
}
\end{figure}

\renewcommand{\figurename}{{\small Figure.}}
\begin{figure}[p]
 \centering\includegraphics[width=15cm,clip]{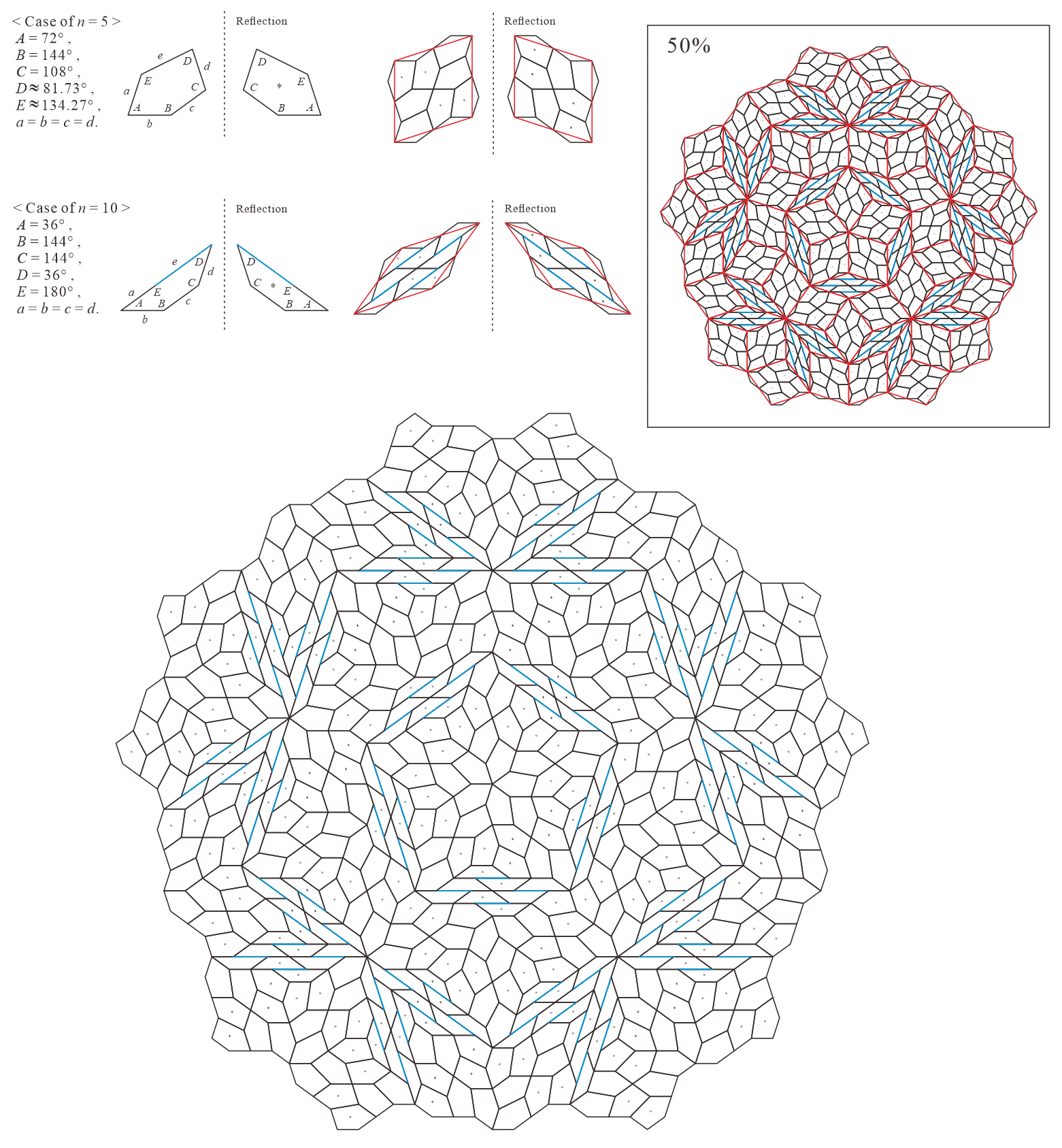} 
  \caption{{\small 
Tiling that generated based on five-fold rotationally symmetric tiling, 
with two types of rhombuses applying conversion of rhombuses into eight 
pentagons satisfying (\ref{eq3}) with $n = 5, 10$ and $\theta = 54^ \circ $ 
(The figure is solely a depiction of the area around the rotationally 
symmetric structure, and the tiling can be spread in all directions)
} 
\label{Fig.5.1.4-2}
}
\end{figure}


\renewcommand{\figurename}{{\small Figure.}}
\begin{figure}[p]
 \centering\includegraphics[width=15cm,clip]{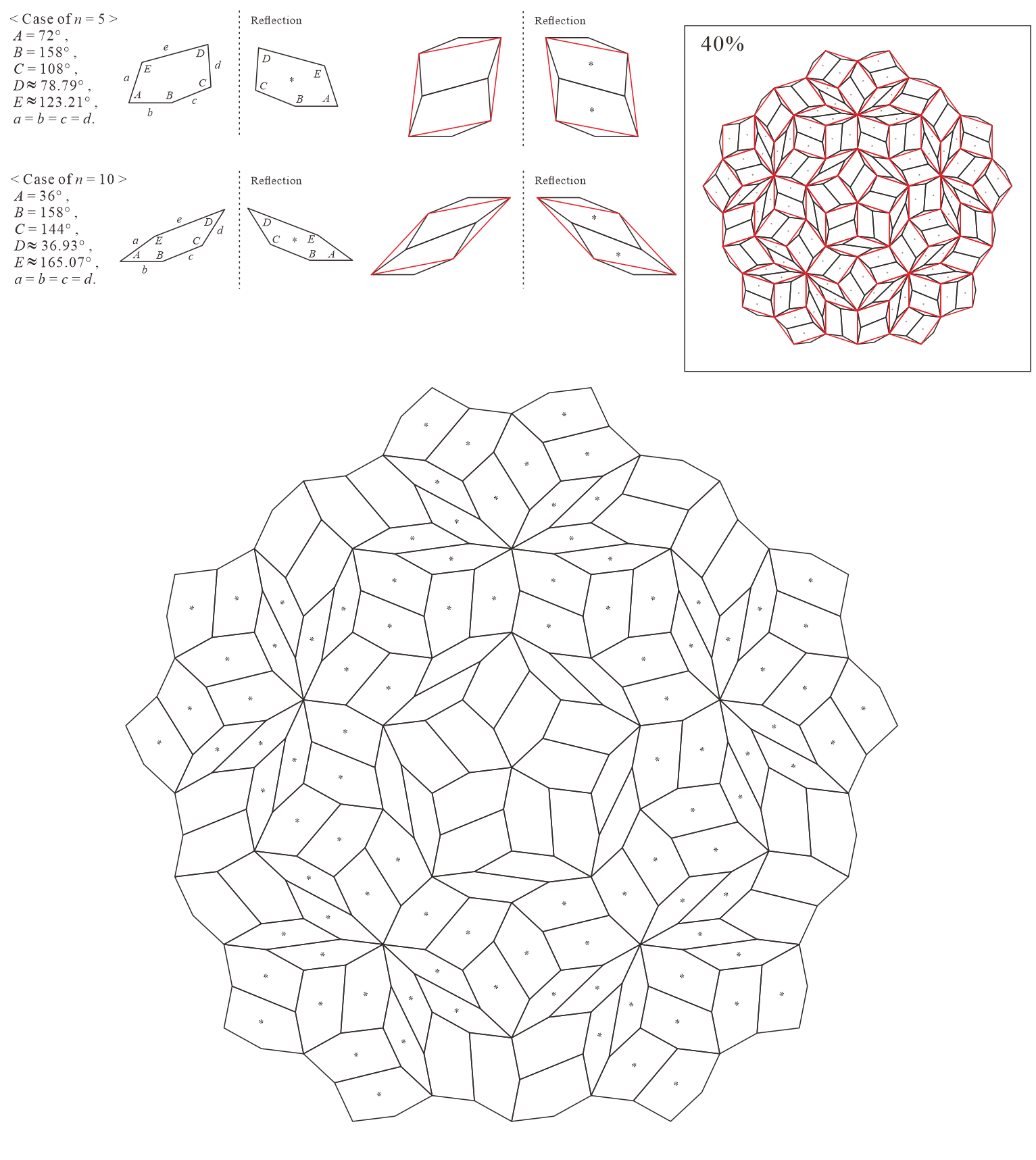} 
  \caption{{\small 
Tiling that generated based on five-fold rotationally symmetric tiling, 
with two types of rhombuses applying conversion of rhombuses into two 
pentagons satisfying (\ref{eq3}) with $n = 5, 10$ and $\theta = 68^ \circ $ 
(The figure is solely a depiction of the area around the rotationally 
symmetric structure, and the tiling can be spread in all directions)
} 
\label{Fig.5.1.5-1}
}
\end{figure}

\renewcommand{\figurename}{{\small Figure.}}
\begin{figure}[p]
 \centering\includegraphics[width=15cm,clip]{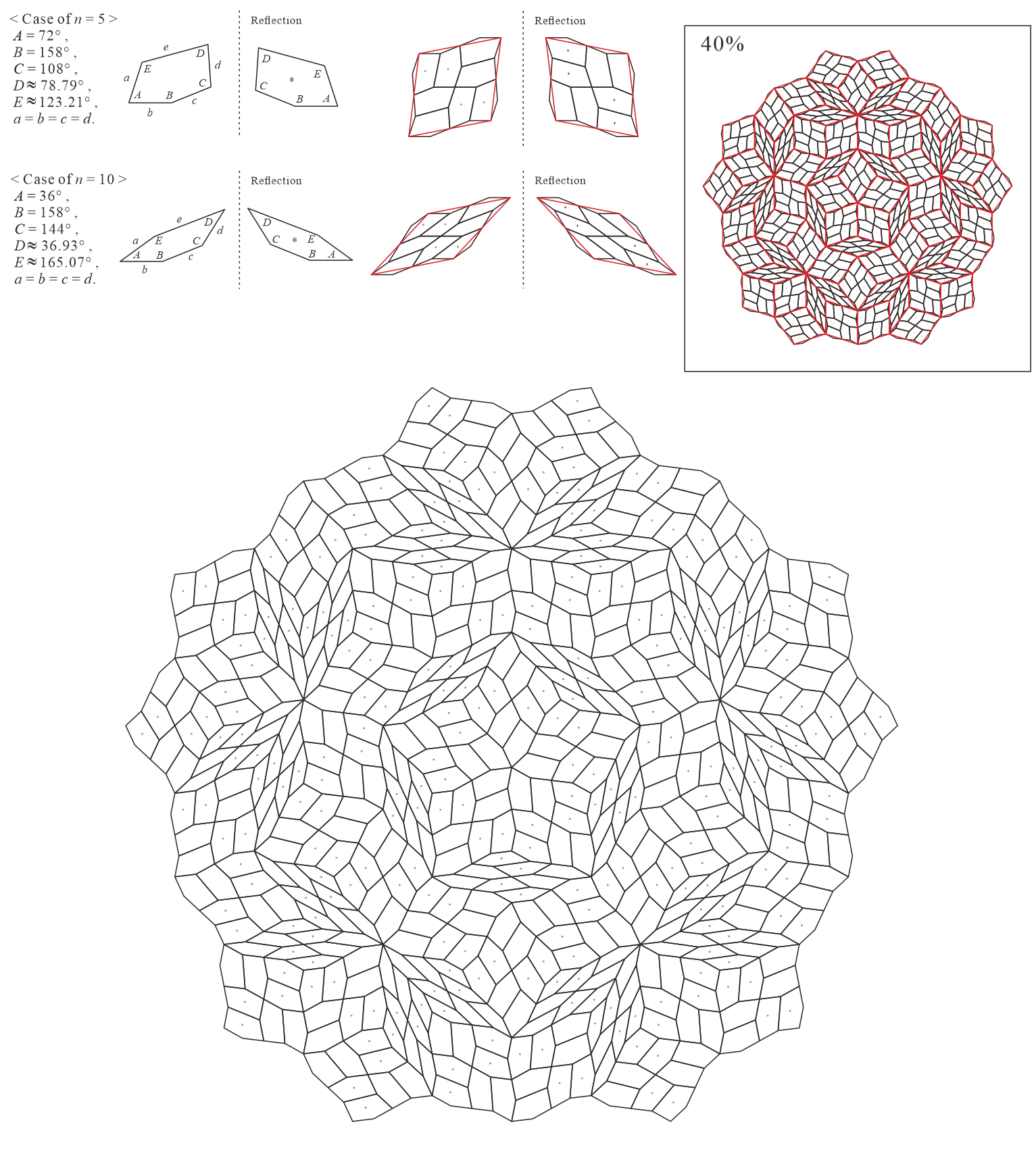} 
  \caption{{\small 
Tiling that generated based on five-fold rotationally symmetric tiling, 
with two types of rhombuses applying conversion of rhombuses into eight 
pentagons satisfying (\ref{eq3}) with $n = 5, 10$ and $\theta = 68^ \circ $ 
(The figure is solely a depiction of the area around the rotationally 
symmetric structure, and the tiling can be spread in all directions)
} 
\label{Fig.5.1.5-2}
}
\end{figure}

\renewcommand{\figurename}{{\small Figure.}}
\begin{figure}[p]
 \centering\includegraphics[width=15cm,clip]{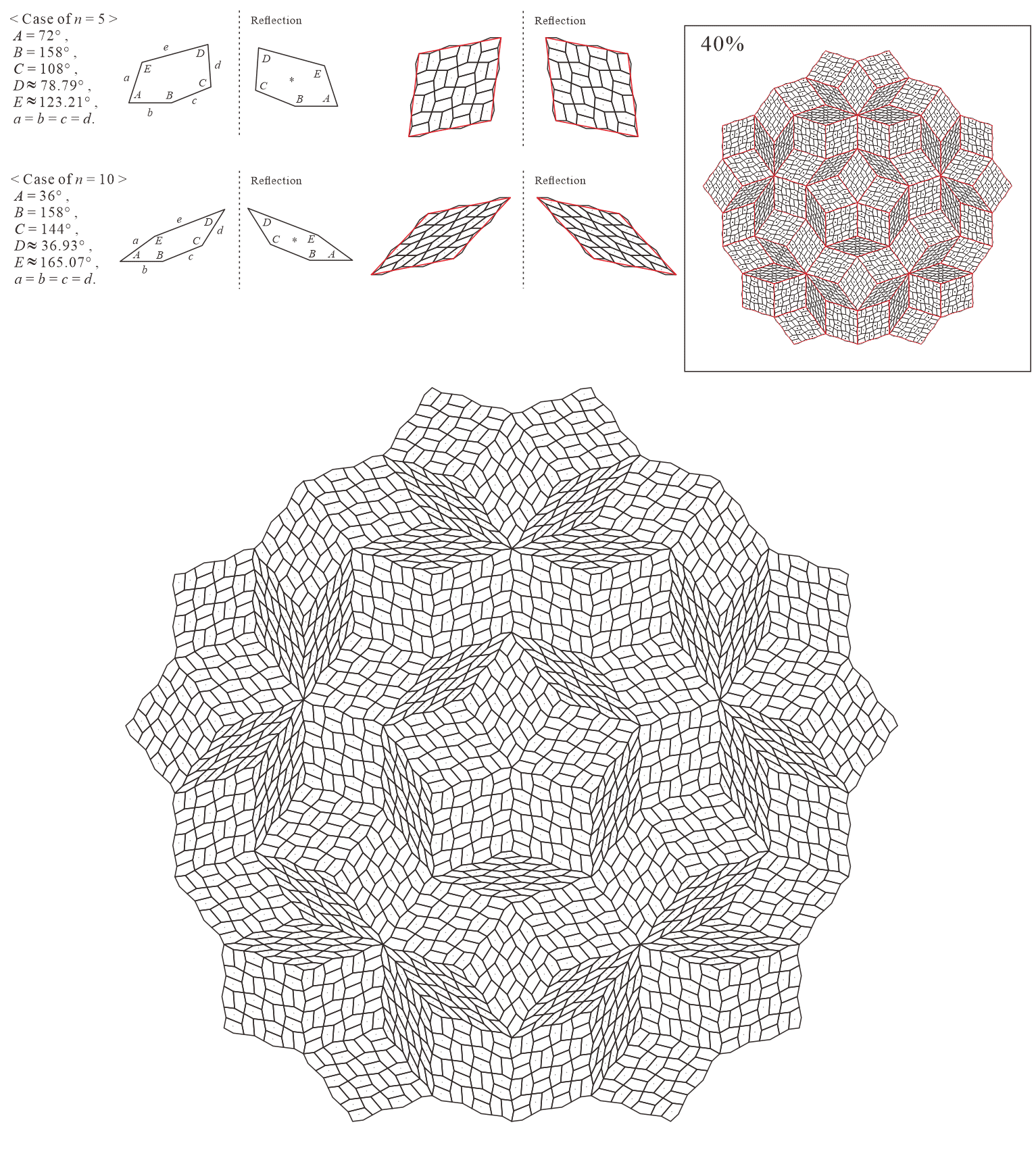} 
  \caption{{\small 
Tiling that generated based on five-fold rotationally symmetric tiling, 
with two types of rhombuses applying conversion of rhombuses into 32 
pentagons satisfying (\ref{eq3}) with $n = 5, 10$ and $\theta = 68^ \circ $ 
(The figure is solely a depiction of the area around the rotationally 
symmetric structure, and the tiling can be spread in all directions)
} 
\label{Fig.5.1.5-3}
}
\end{figure}

\renewcommand{\figurename}{{\small Figure.}}
\begin{figure}[p]
 \centering\includegraphics[width=15cm,clip]{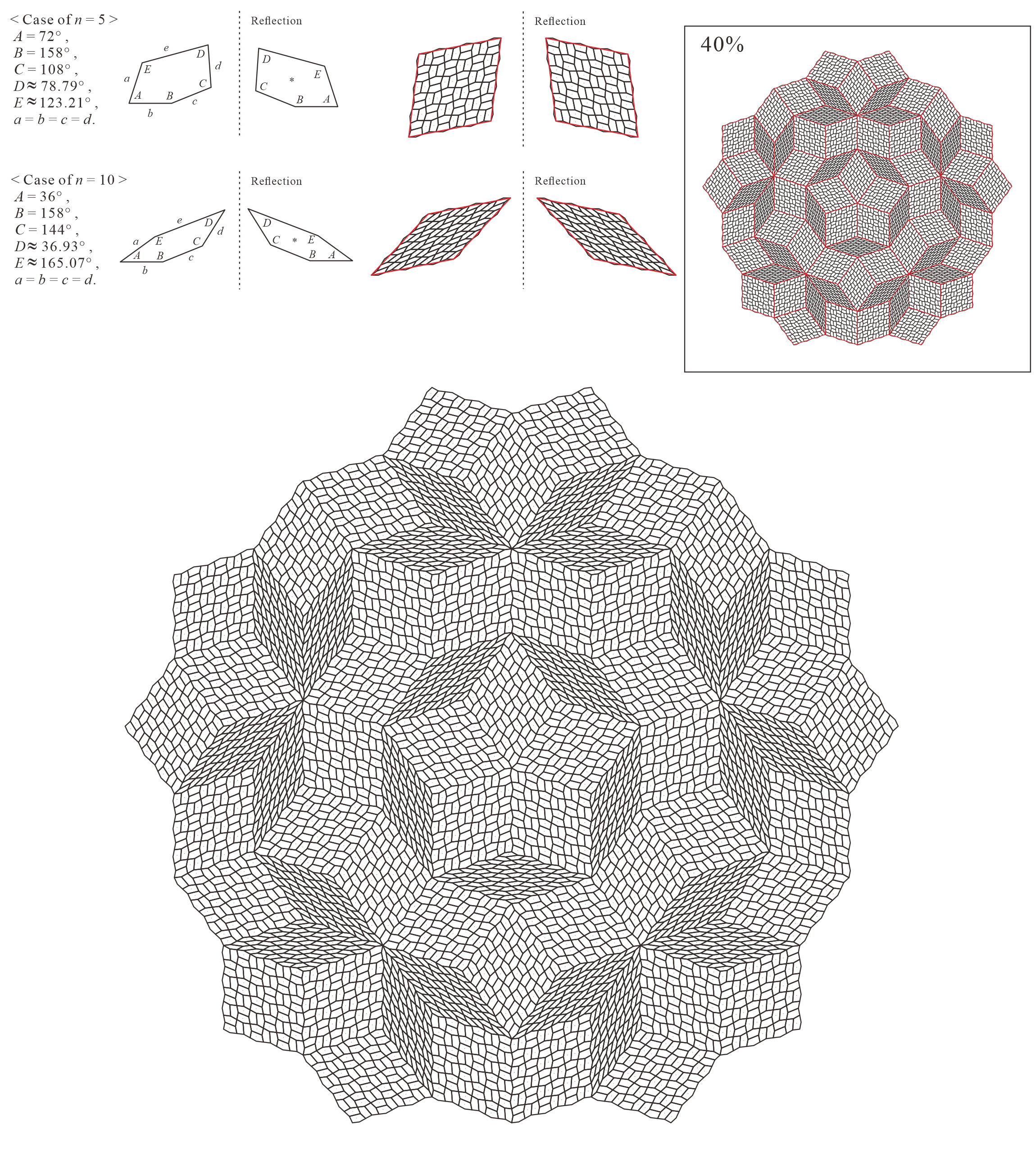} 
  \caption{{\small 
Tiling that generated based on five-fold rotationally symmetric tiling, 
with two types of rhombuses applying conversion of rhombuses into 72 
pentagons satisfying (\ref{eq3}) with $n = 5, 10$ and $\theta = 68^ \circ $ 
(The figure is solely a depiction of the area around the rotationally 
symmetric structure, and the tiling can be spread in all directions)
} 
\label{Fig.5.1.5-4}
}
\end{figure}


\renewcommand{\figurename}{{\small Figure.}}
\begin{figure}[p]
 \centering\includegraphics[width=15cm,clip]{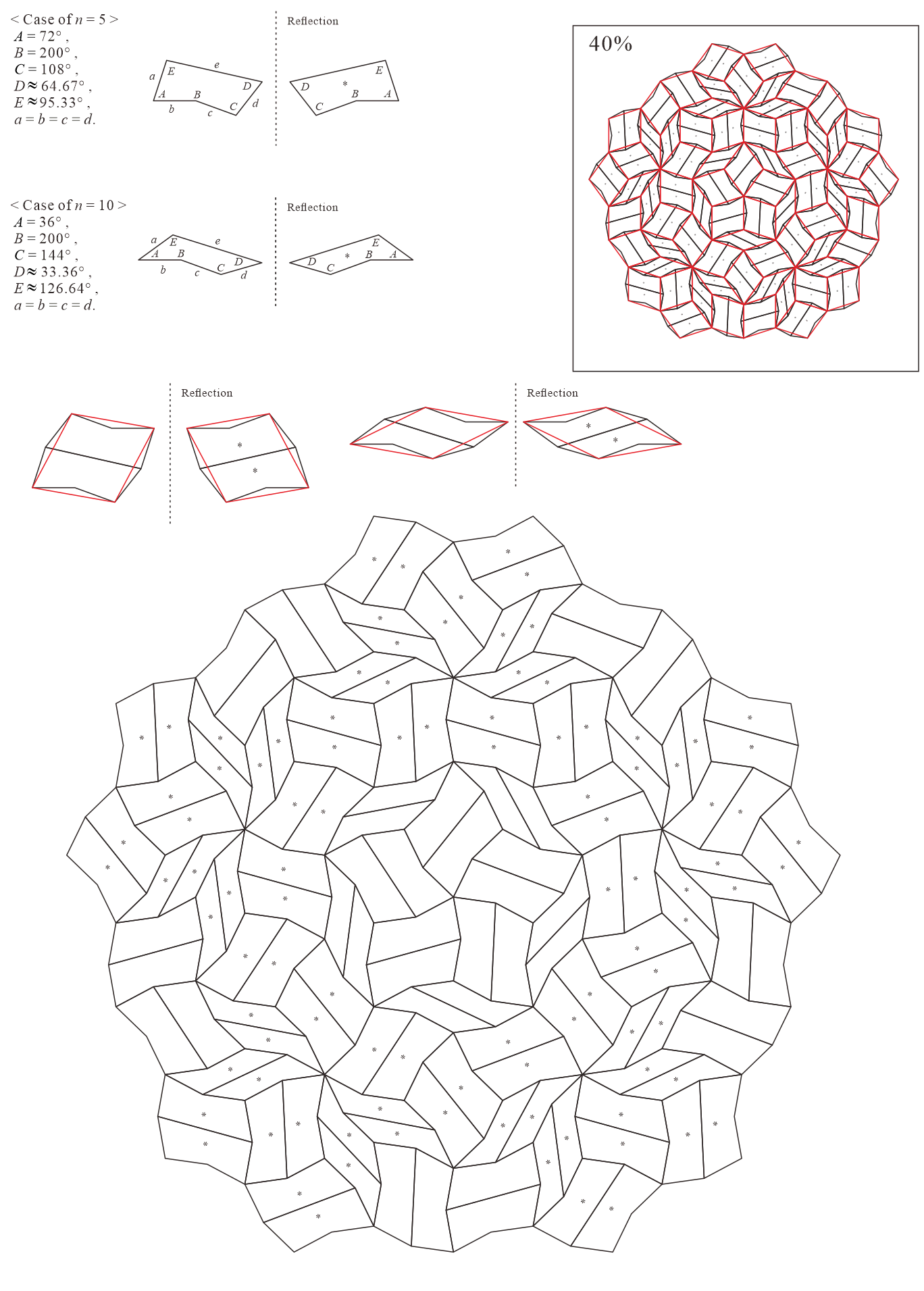} 
  \caption{{\small 
Tiling that generated based on five-fold rotationally symmetric tiling, 
with two types of rhombuses applying conversion of rhombuses into two 
pentagons satisfying (\ref{eq3}) with $n = 5, 10$ and $\theta = 110^ \circ $ 
(The figure is solely a depiction of the area around the rotationally 
symmetric structure, and the tiling can be spread in all directions)
} 
\label{Fig.5.1.7-1}
}
\end{figure}

\renewcommand{\figurename}{{\small Figure.}}
\begin{figure}[p]
 \centering\includegraphics[width=15cm,clip]{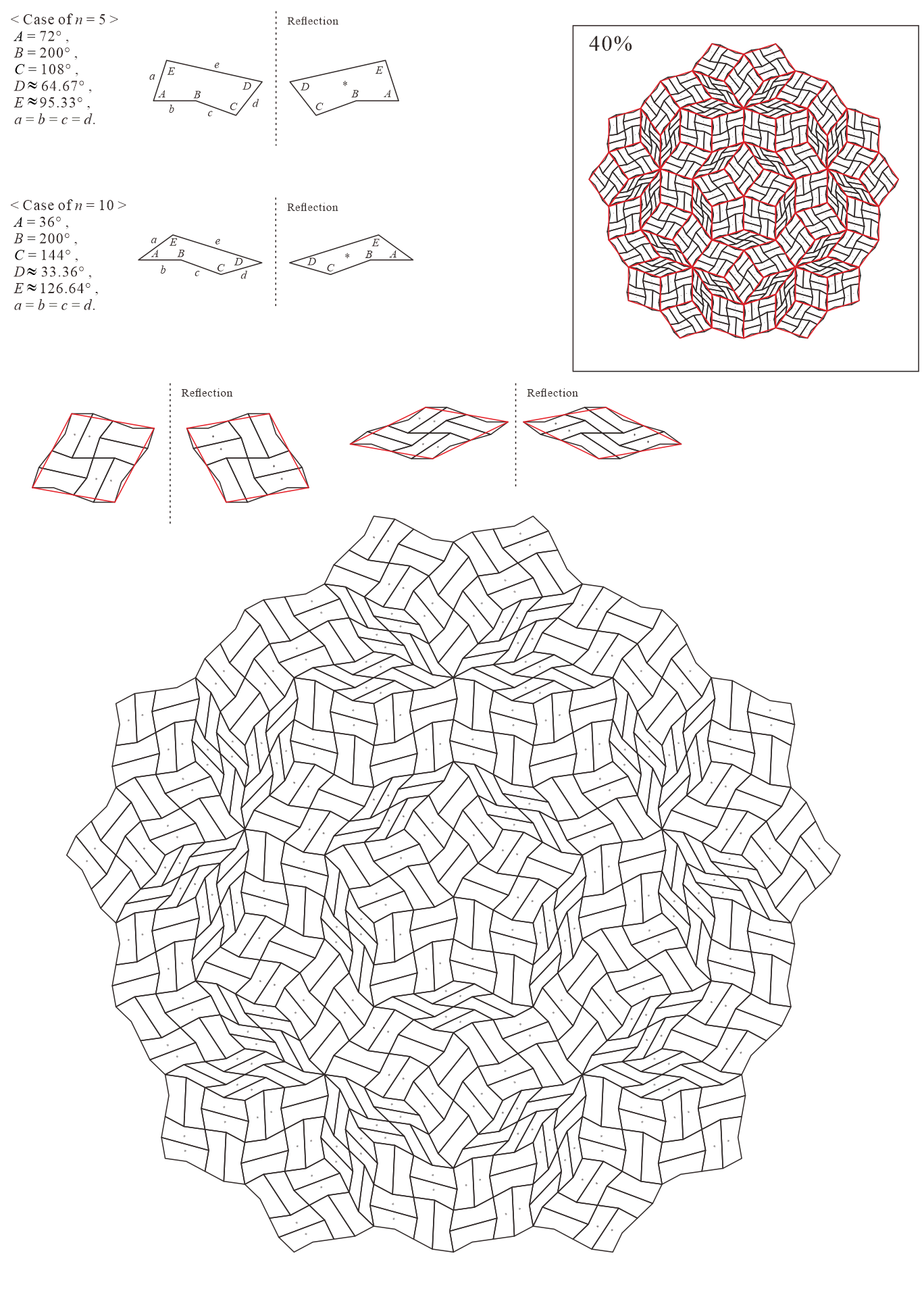} 
  \caption{{\small 
Tiling that generated based on five-fold rotationally symmetric tiling, 
with two types of rhombuses applying conversion of rhombuses into eight 
pentagons satisfying (\ref{eq3}) with $n = 5, 10$ and $\theta = 110^ \circ $ 
(The figure is solely a depiction of the area around the rotationally 
symmetric structure, and the tiling can be spread in all directions)
} 
\label{Fig.5.1.7-2}
}
\end{figure}



\renewcommand{\figurename}{{\small Figure.}}
\begin{figure}[p]
 \centering\includegraphics[width=15cm,clip]{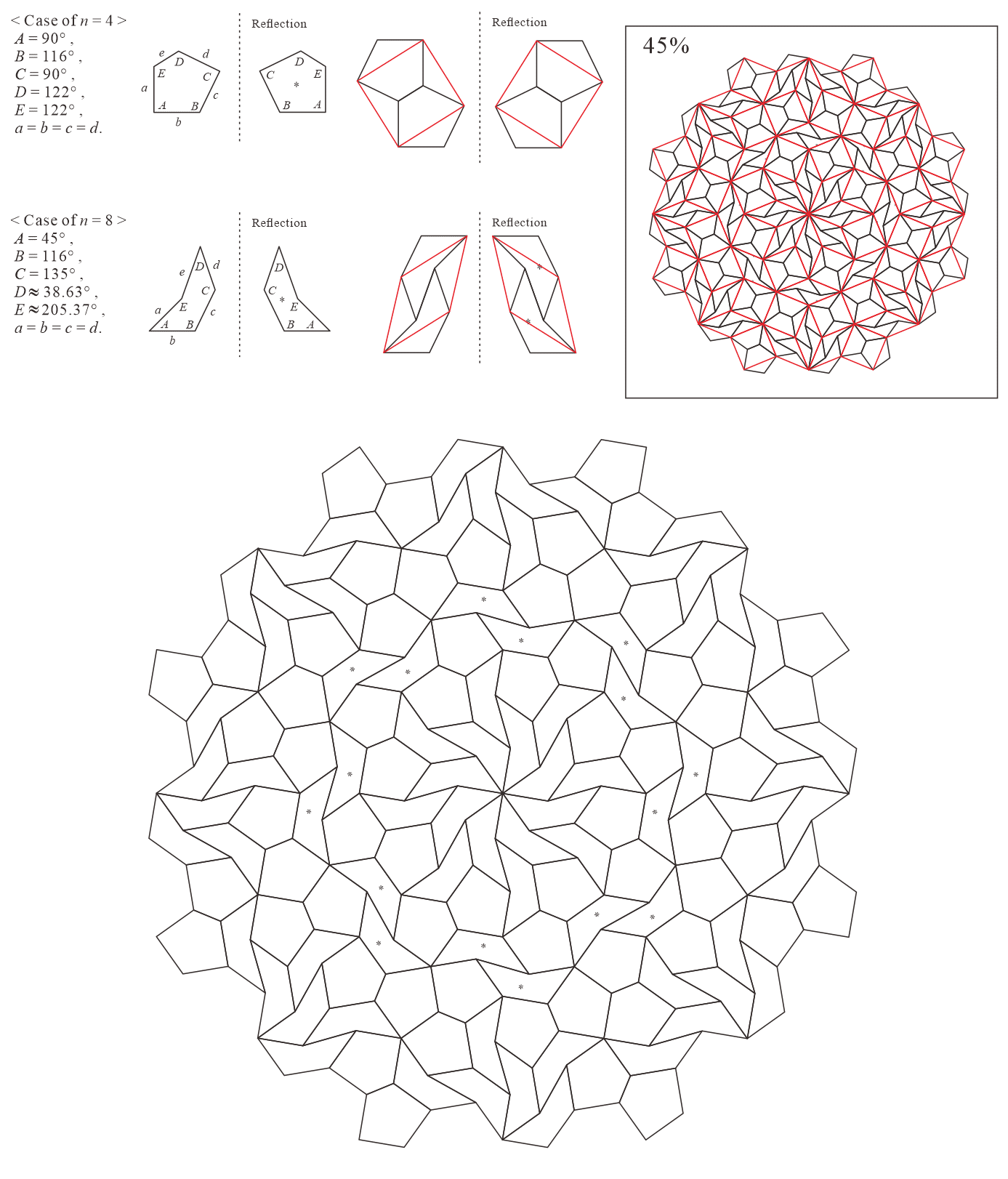} 
  \caption{{\small 
Tiling that generated based on eight-fold rotationally symmetric tiling, 
with two types of rhombuses applying conversion of rhombuses into two 
pentagons satisfying (\ref{eq3}) with $n = 4, 8$ and $\theta = 26^ \circ $ 
(The figure is solely a depiction of the area around the rotationally 
symmetric structure, and the tiling can be spread in all directions)
} 
\label{Fig.5.2.1-1}
}
\end{figure}

\renewcommand{\figurename}{{\small Figure.}}
\begin{figure}[p]
 \centering\includegraphics[width=15cm,clip]{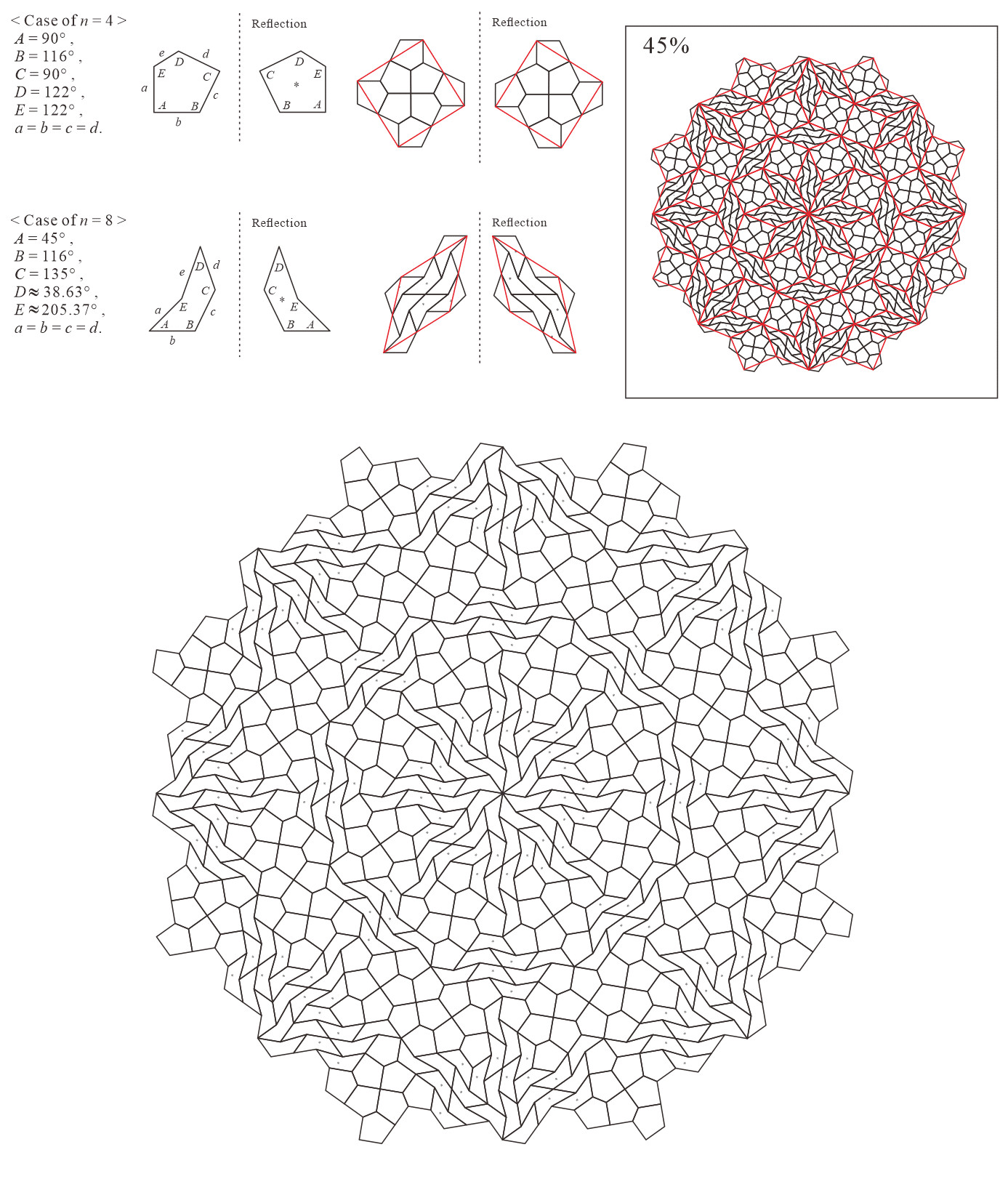} 
  \caption{{\small 
Tiling that generated based on eight-fold rotationally symmetric tiling, 
with two types of rhombuses applying conversion of rhombuses into eight 
pentagons satisfying (\ref{eq3}) with $n = 4, 8$ and $\theta = 26^ \circ $ 
(The figure is solely a depiction of the area around the rotationally 
symmetric structure, and the tiling can be spread in all directions)
} 
\label{Fig.5.2.1-2}
}
\end{figure}


\renewcommand{\figurename}{{\small Figure.}}
\begin{figure}[p]
 \centering\includegraphics[width=15cm,clip]{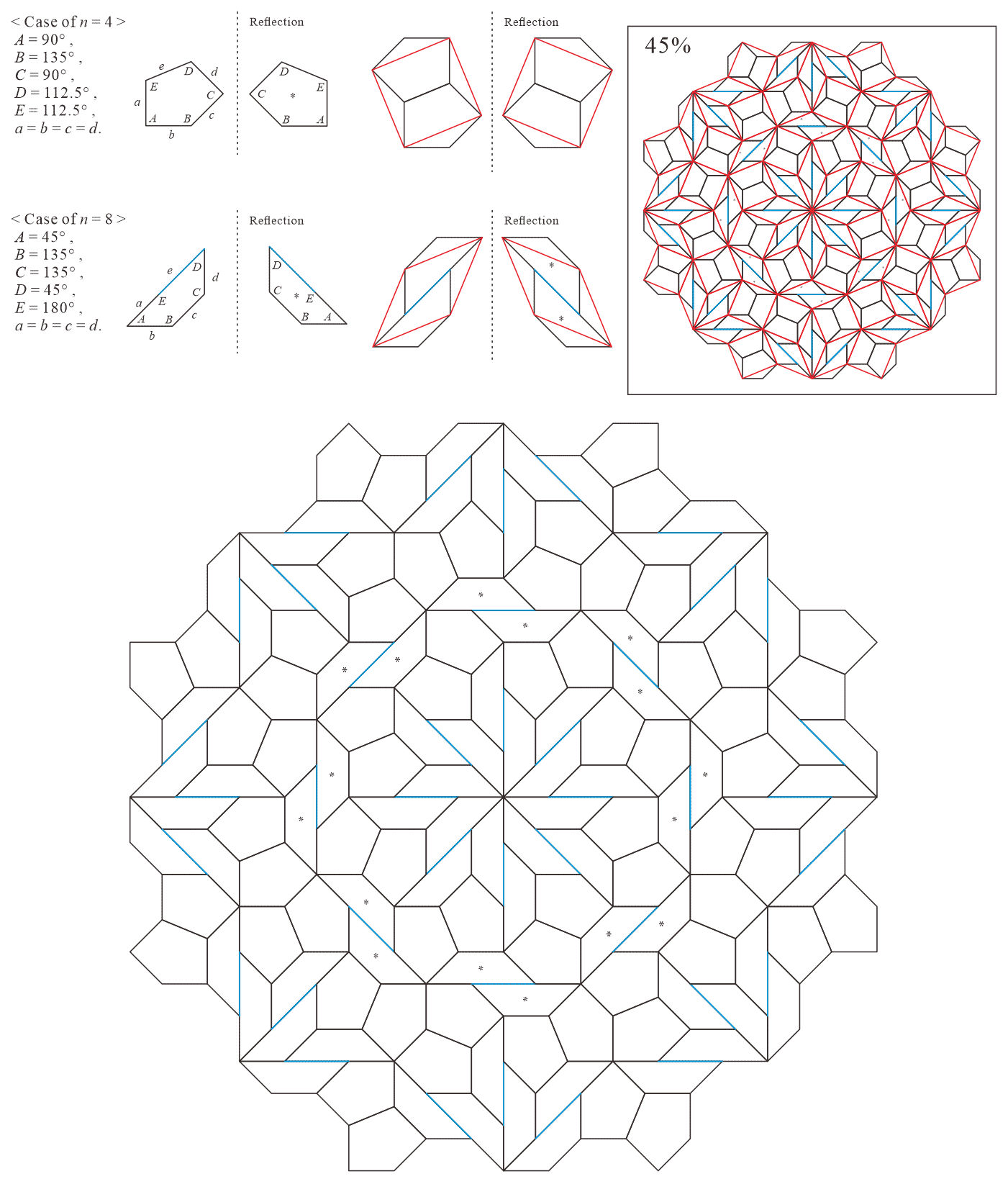} 
  \caption{{\small 
Tiling that generated based on  eight-fold rotationally symmetric tiling, 
with two types of rhombuses applying conversion of rhombuses into two 
pentagons satisfying (\ref{eq3}) with $n = 4, 8$ and $\theta = 45^ \circ $ 
(The figure is solely a depiction of the area around the rotationally 
symmetric structure, and the tiling can be spread in all directions)
} 
\label{Fig.5.2.2-1}
}
\end{figure}

\renewcommand{\figurename}{{\small Figure.}}
\begin{figure}[p]
 \centering\includegraphics[width=15cm,clip]{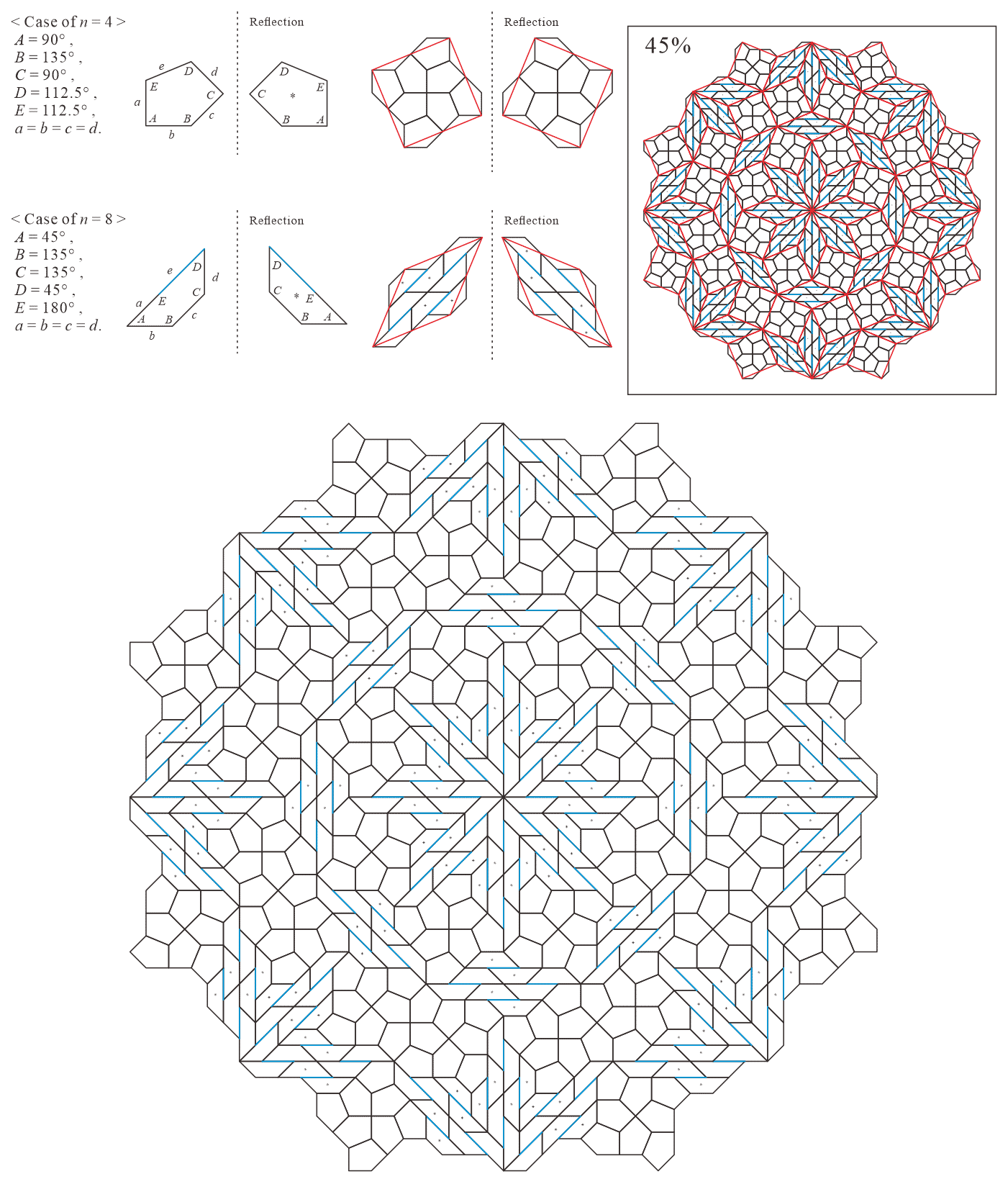} 
  \caption{{\small 
Tiling that generated based on eight-fold rotationally symmetric tiling, 
with two types of rhombuses applying conversion of rhombuses into eight 
pentagons satisfying (\ref{eq3}) with $n = 4, 8$ and $\theta = 45^ \circ $ 
(The figure is solely a depiction of the area around the rotationally 
symmetric structure, and the tiling can be spread in all directions)
} 
\label{Fig.5.2.2-2}
}
\end{figure}

\renewcommand{\figurename}{{\small Figure.}}
\begin{figure}[p]
 \centering\includegraphics[width=15cm,clip]{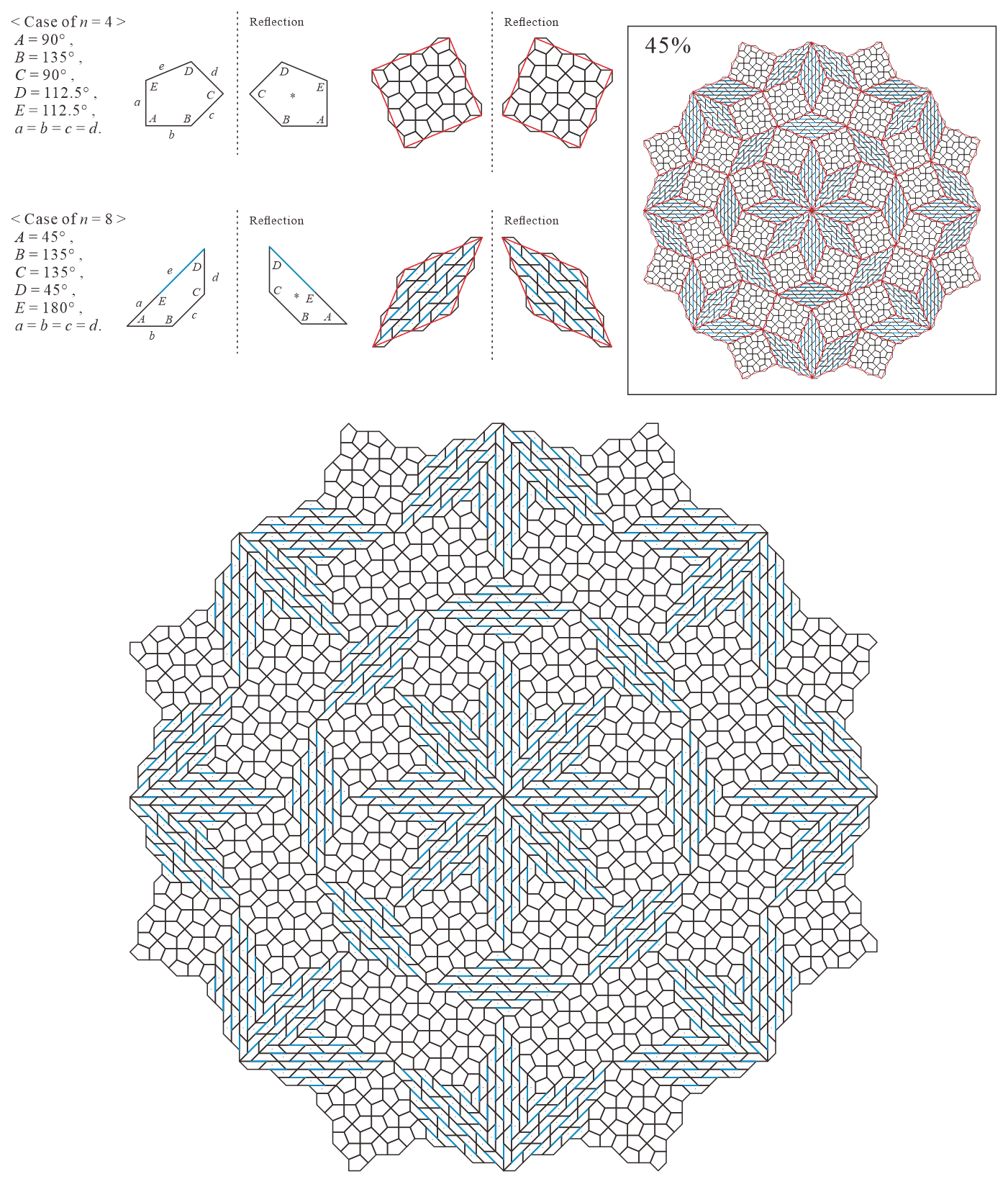} 
  \caption{{\small 
Tiling that generated based on eight-fold rotationally symmetric tiling, 
with two types of rhombuses applying conversion of rhombuses into 32 
pentagons satisfying (\ref{eq3}) with $n = 4, 8$ and $\theta = 45^ \circ $ 
(The figure is solely a depiction of the area around the rotationally 
symmetric structure, and the tiling can be spread in all directions)
} 
\label{Fig.5.2.2-3}
}
\end{figure}


\renewcommand{\figurename}{{\small Figure.}}
\begin{figure}[p]
 \centering\includegraphics[width=15cm,clip]{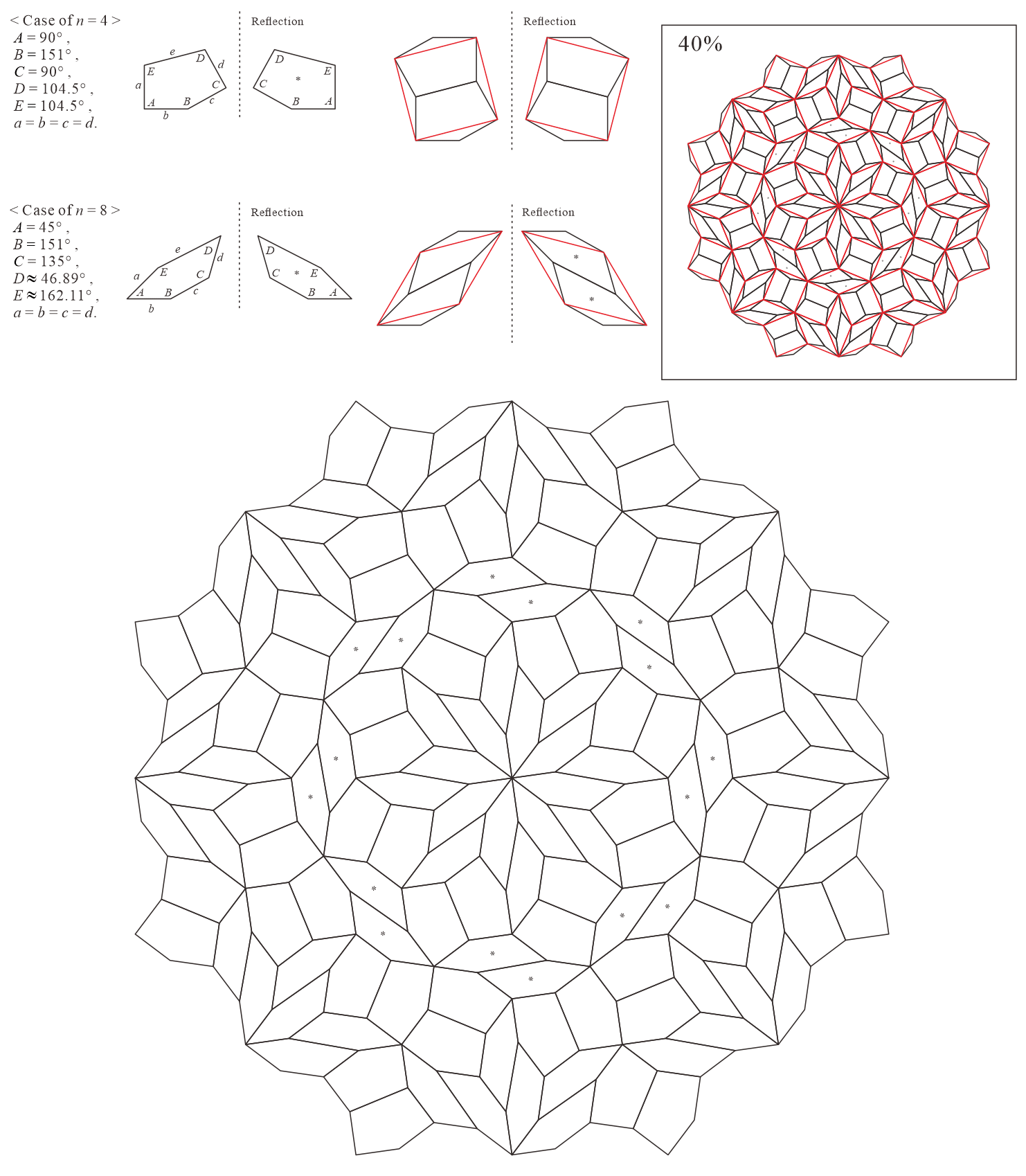} 
  \caption{{\small 
Tiling that generated based on eight-fold rotationally symmetric tiling, 
with two types of rhombuses applying conversion of rhombuses into two 
pentagons satisfying (\ref{eq3}) with $n = 4, 8$ and $\theta = 61^ \circ $ 
(The figure is solely a depiction of the area around the rotationally 
symmetric structure, and the tiling can be spread in all directions)
} 
\label{Fig.5.2.3-1}
}
\end{figure}

\renewcommand{\figurename}{{\small Figure.}}
\begin{figure}[p]
 \centering\includegraphics[width=15cm,clip]{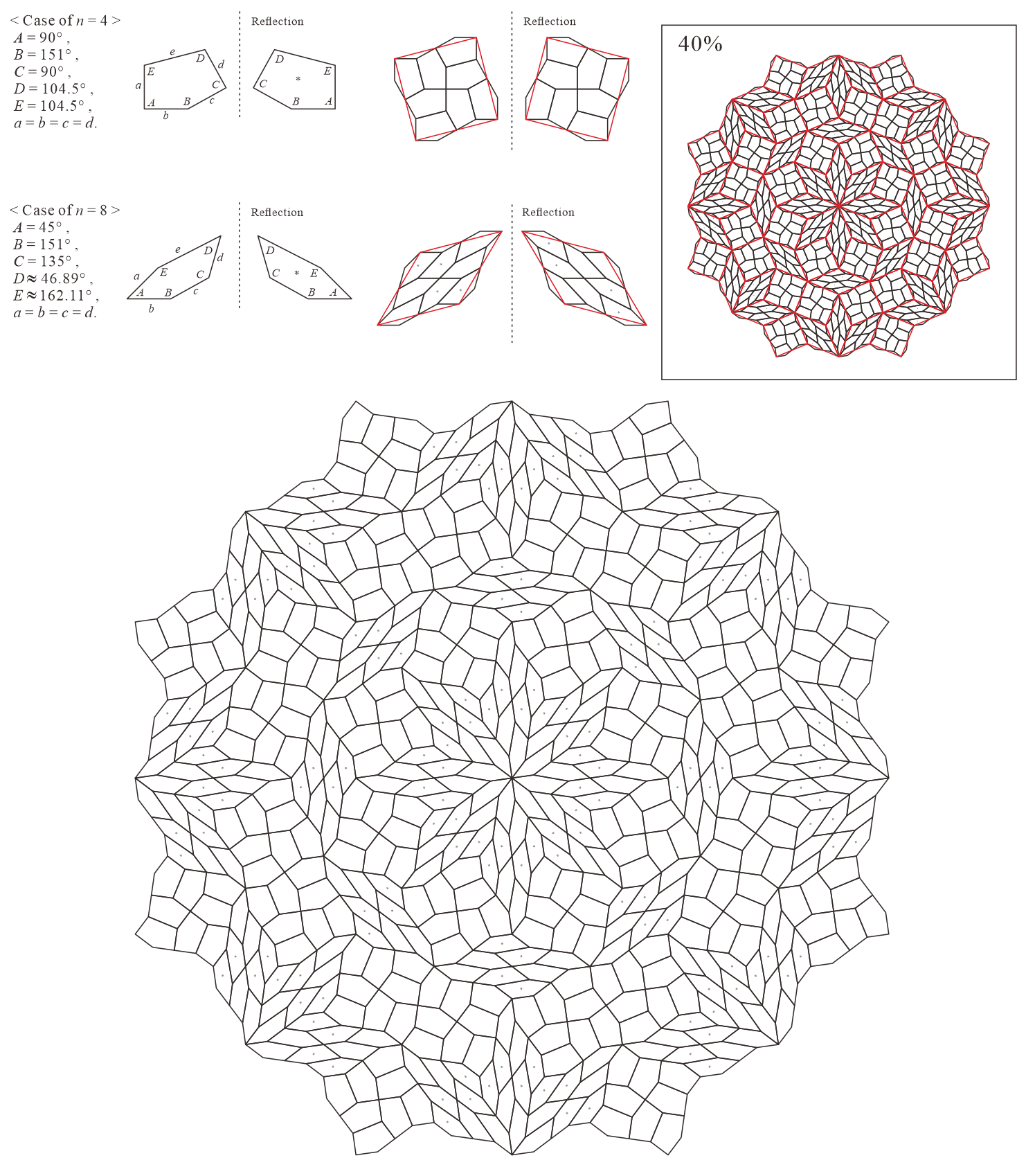} 
  \caption{{\small 
Tiling that generated based on eight-fold rotationally symmetric tiling, 
with two types of rhombuses applying conversion of rhombuses into eight 
pentagons satisfying (\ref{eq3}) with $n = 4, 8$ and $\theta = 61^ \circ $ 
(The figure is solely a depiction of the area around the rotationally 
symmetric structure, and the tiling can be spread in all directions)
} 
\label{Fig.5.2.3-2}
}
\end{figure}

\renewcommand{\figurename}{{\small Figure.}}
\begin{figure}[p]
 \centering\includegraphics[width=15cm,clip]{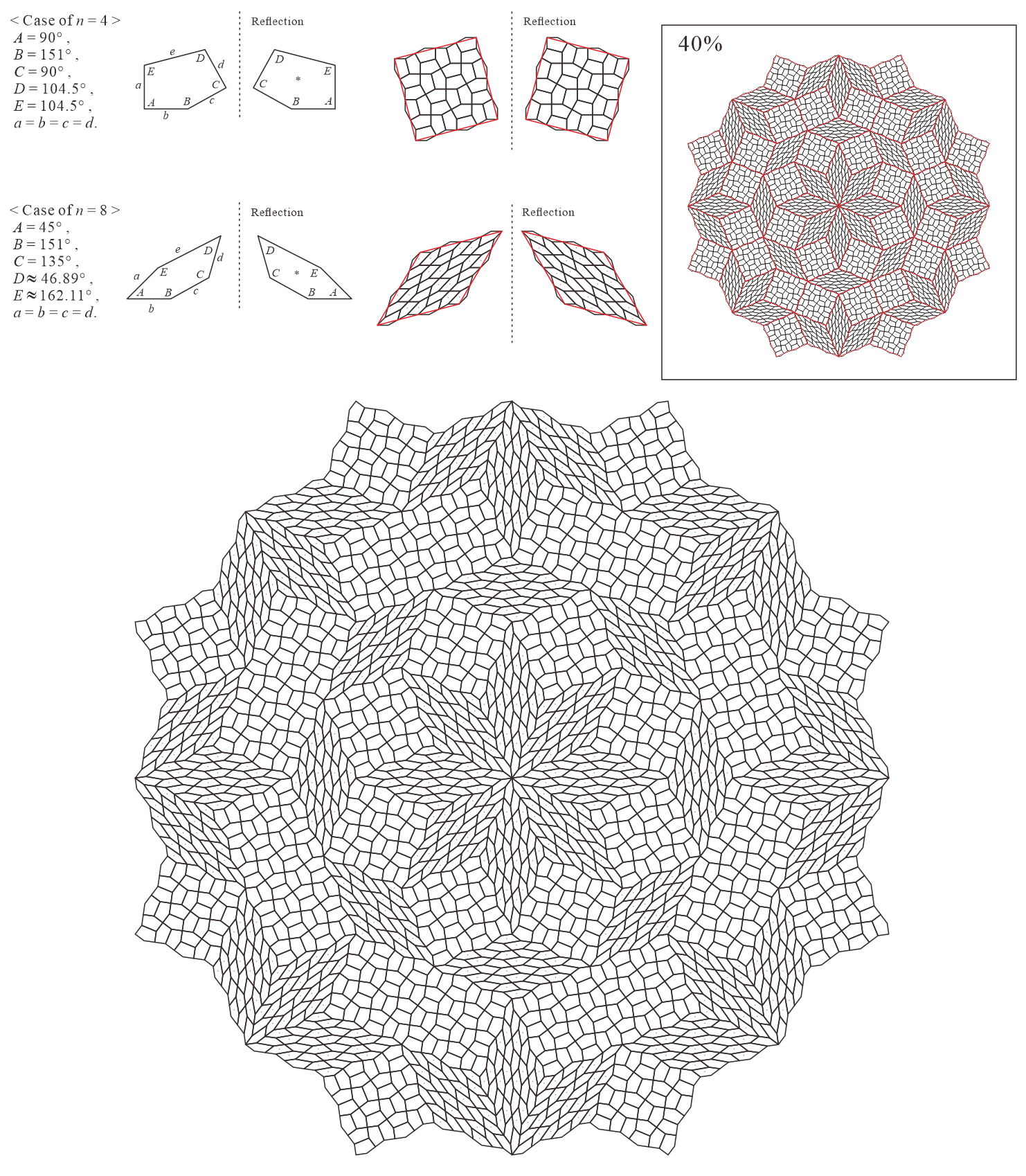} 
  \caption{{\small 
Tiling that generated based on eight-fold rotationally symmetric tiling, 
with two types of rhombuses applying conversion of rhombuses into 32 
pentagons satisfying (\ref{eq3}) with $n = 4, 8$ and $\theta = 61^ \circ $ 
(The figure is solely a depiction of the area around the rotationally 
symmetric structure, and the tiling can be spread in all directions)
} 
\label{Fig.5.2.3-3}
}
\end{figure}


\renewcommand{\figurename}{{\small Figure.}}
\begin{figure}[p]
 \centering\includegraphics[width=15cm,clip]{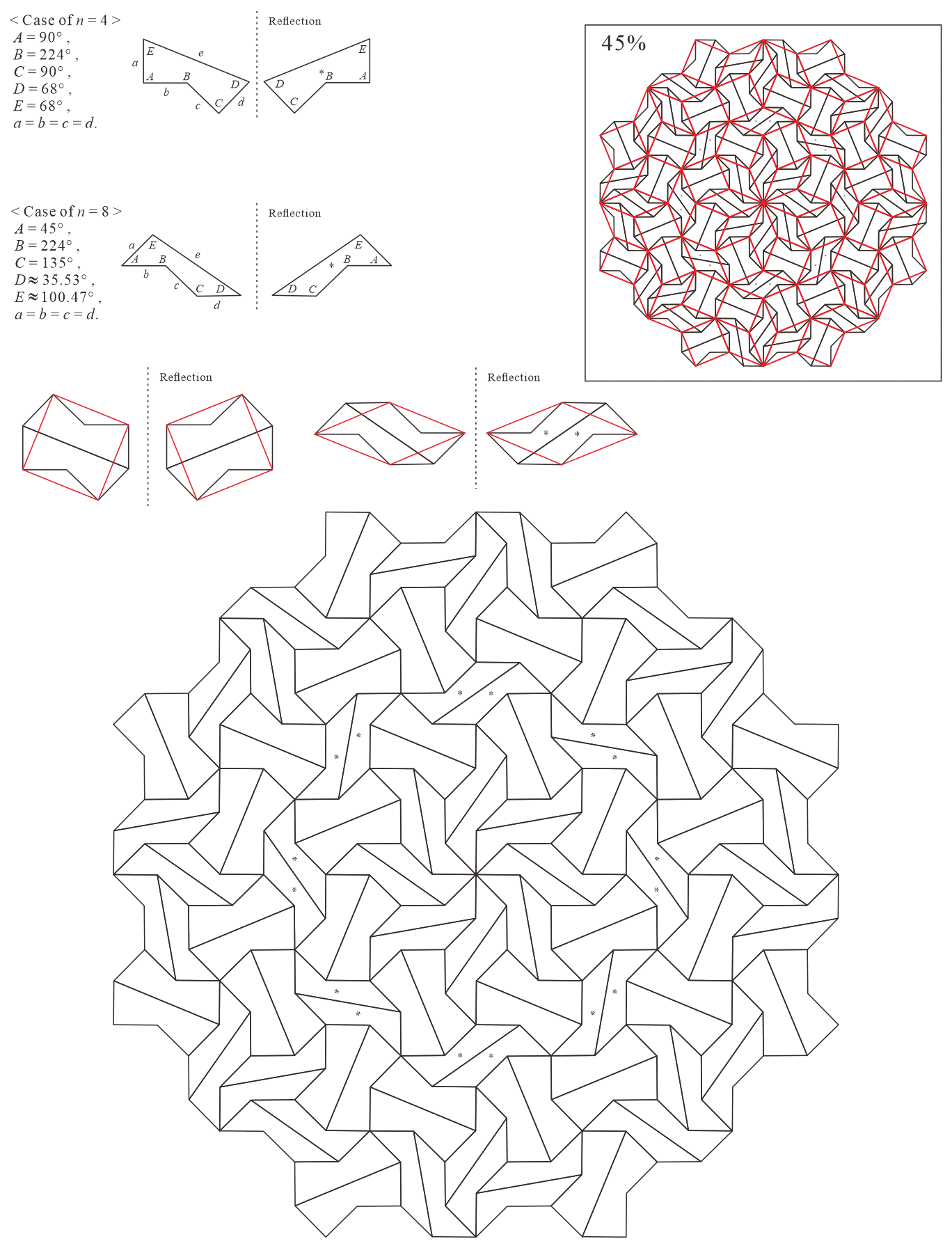} 
  \caption{{\small 
Tiling that generated based on eight-fold rotationally symmetric tiling, 
with two types of rhombuses applying conversion of rhombuses into two 
pentagons satisfying (\ref{eq3}) with $n = 4, 8$ and $\theta = 134^ \circ $ 
(The figure is solely a depiction of the area around the rotationally 
symmetric structure, and the tiling can be spread in all directions)
} 
\label{Fig.5.2.5-1}
}
\end{figure}

\renewcommand{\figurename}{{\small Figure.}}
\begin{figure}[p]
 \centering\includegraphics[width=15cm,clip]{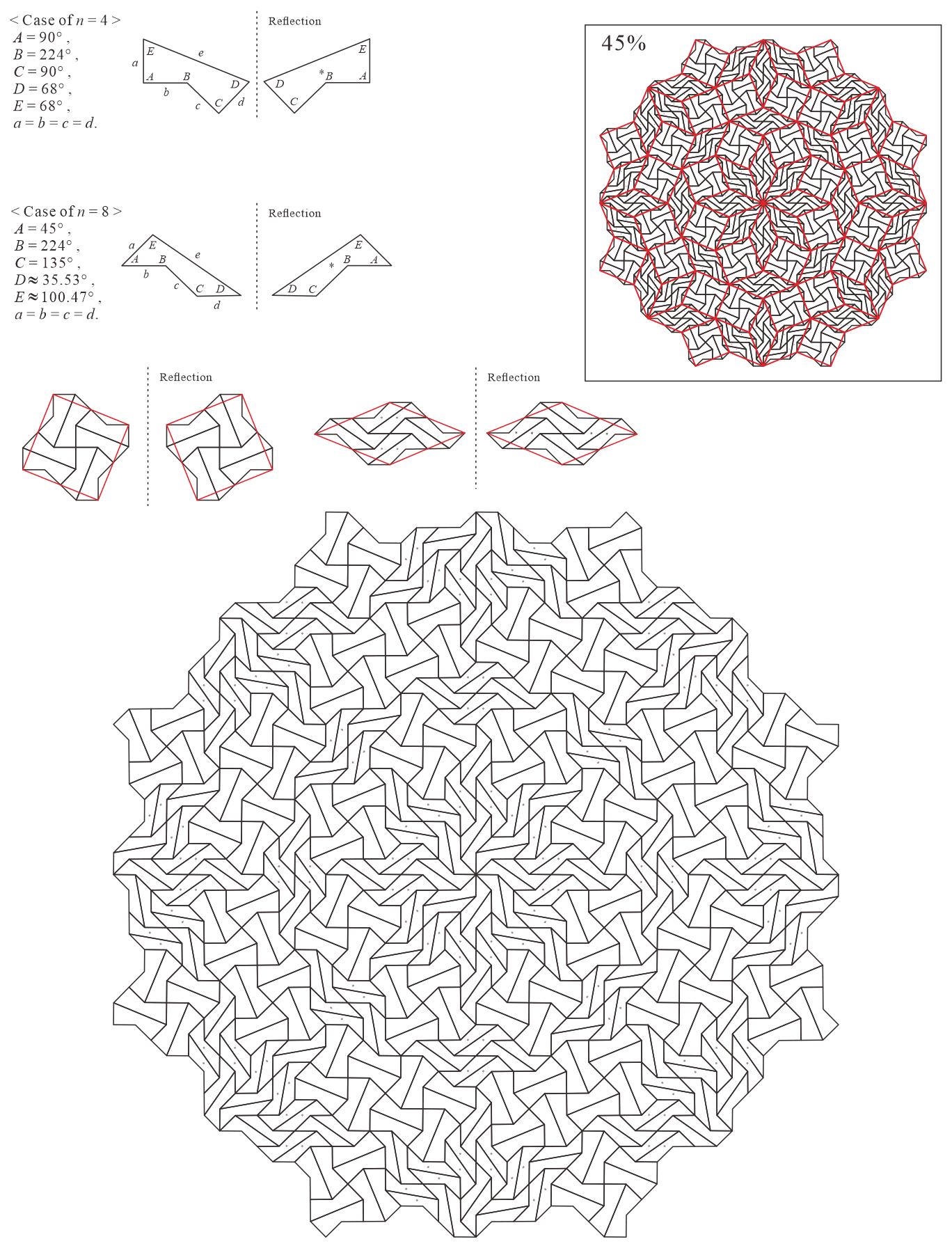} 
  \caption{{\small 
Tiling that generated based on eight-fold rotationally symmetric tiling, 
with two types of rhombuses applying conversion of rhombuses into eight 
pentagons satisfying (\ref{eq3}) with $n = 4, 8$ and $\theta = 134^ \circ $ 
(The figure is solely a depiction of the area around the rotationally 
symmetric structure, and the tiling can be spread in all directions)
} 
\label{Fig.5.2.5-2}
}
\end{figure}



\renewcommand{\figurename}{{\small Figure.}}
\begin{figure}[p]
 \centering\includegraphics[width=15cm,clip]{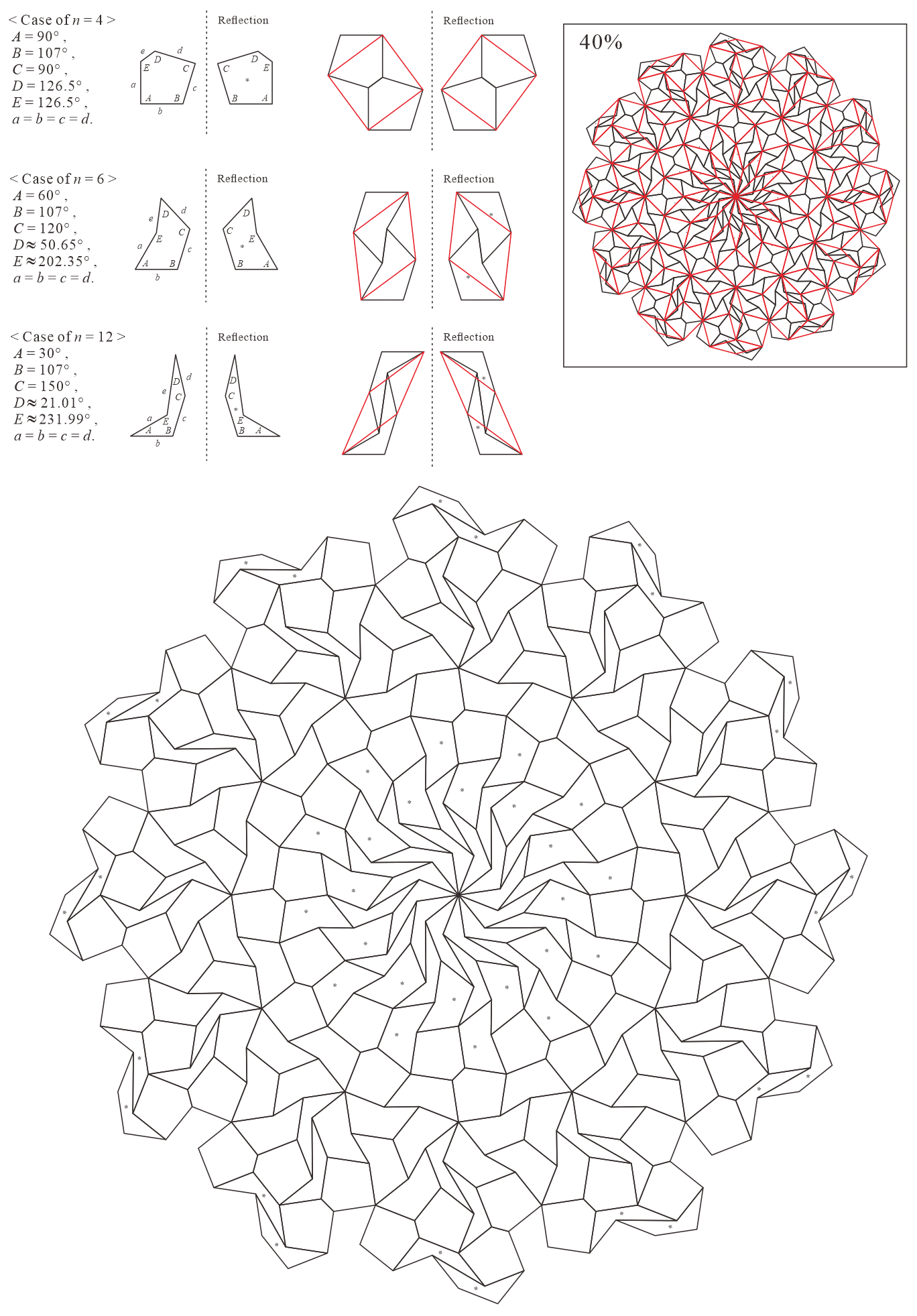} 
  \caption{{\small 
Tiling that generated based on 12-fold rotationally symmetric tiling, 
with three types of rhombuses applying conversion of rhombuses into two 
pentagons satisfying (\ref{eq3}) with $n = 4, 6, 12$ and $\theta = 17^ \circ $ 
(The figure is solely a depiction of the area around the rotationally 
symmetric structure, and the tiling can be spread in all directions)
} 
\label{Fig.5.3.1-1}
}
\end{figure}

\renewcommand{\figurename}{{\small Figure.}}
\begin{figure}[p]
 \centering\includegraphics[width=15cm,clip]{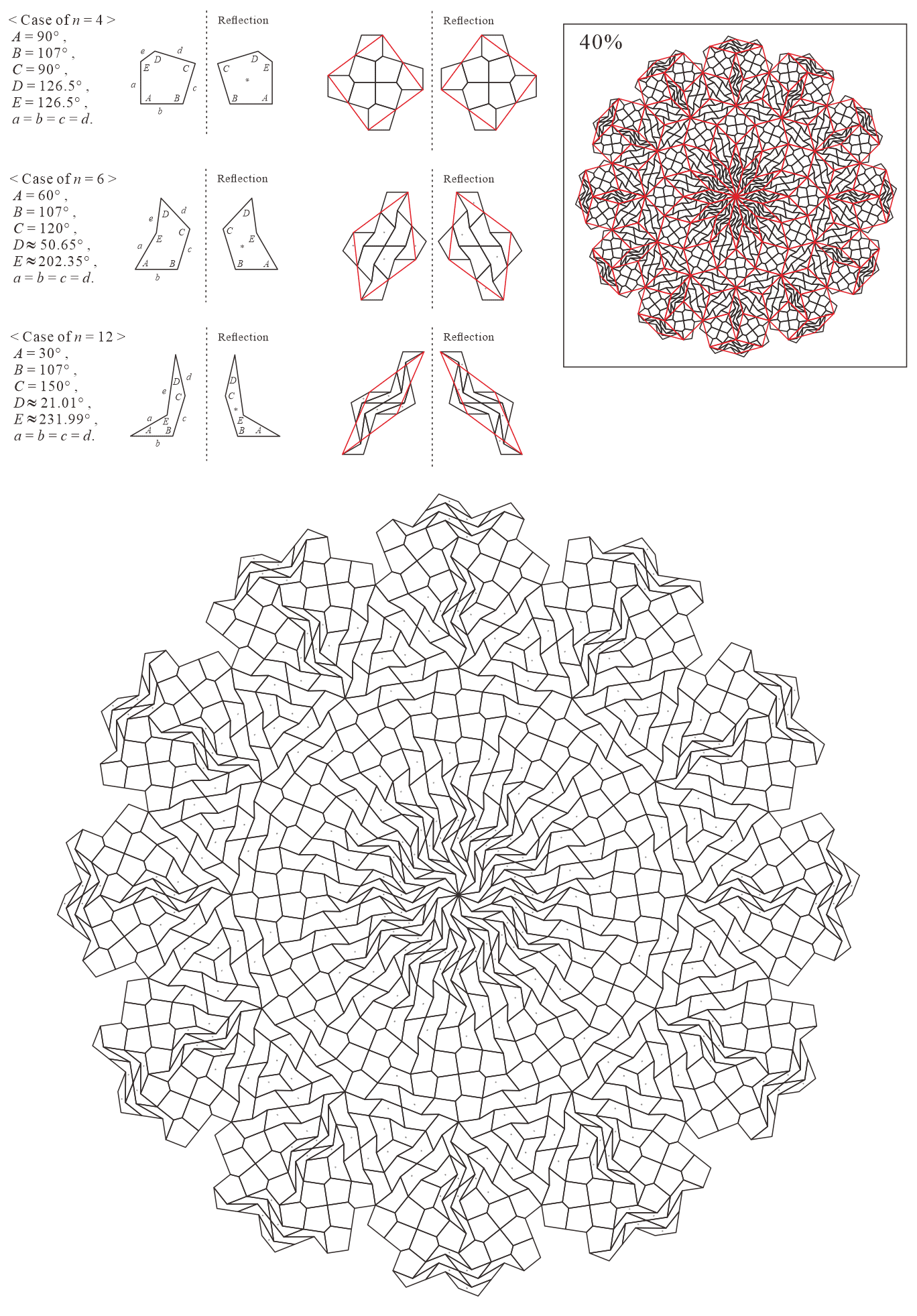} 
  \caption{{\small 
Tiling that generated based on 12-fold rotationally symmetric tiling, 
with three types of rhombuses applying conversion of rhombuses into eight 
pentagons satisfying (\ref{eq3}) with $n = 4, 6, 12$ and $\theta = 17^ \circ $ 
(The figure is solely a depiction of the area around the rotationally 
symmetric structure, and the tiling can be spread in all directions)
} 
\label{Fig.5.3.1-2}
}
\end{figure}


\renewcommand{\figurename}{{\small Figure.}}
\begin{figure}[p]
 \centering\includegraphics[width=15cm,clip]{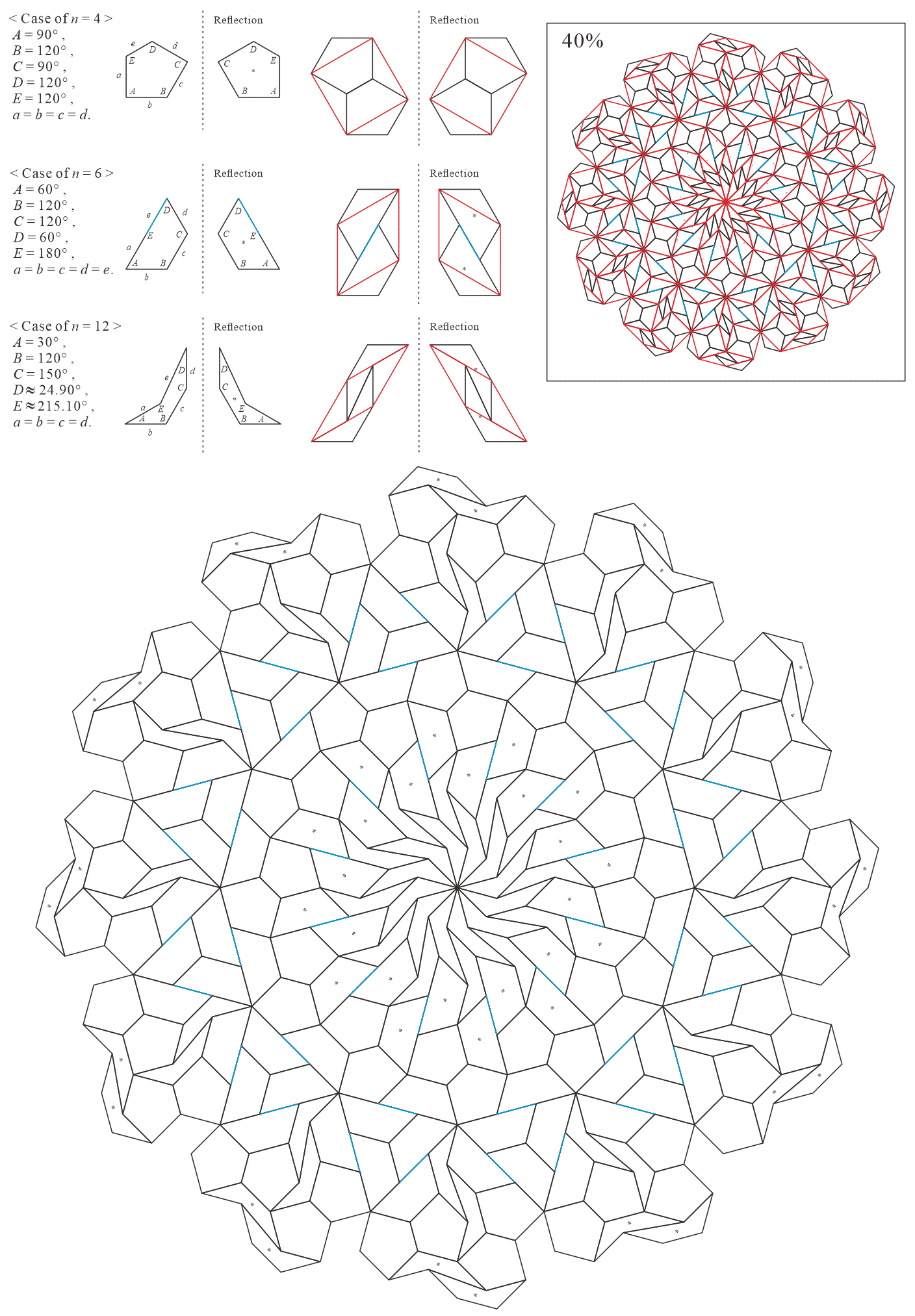} 
  \caption{{\small 
Tiling that generated based on 12-fold rotationally symmetric tiling, 
with three types of rhombuses applying conversion of rhombuses into two 
pentagons satisfying (\ref{eq3}) with $n = 4, 6, 12$ and $\theta = 30^ \circ $ 
(The figure is solely a depiction of the area around the rotationally 
symmetric structure, and the tiling can be spread in all directions)
} 
\label{Fig.5.3.2-1}
}
\end{figure}

\renewcommand{\figurename}{{\small Figure.}}
\begin{figure}[p]
 \centering\includegraphics[width=15cm,clip]{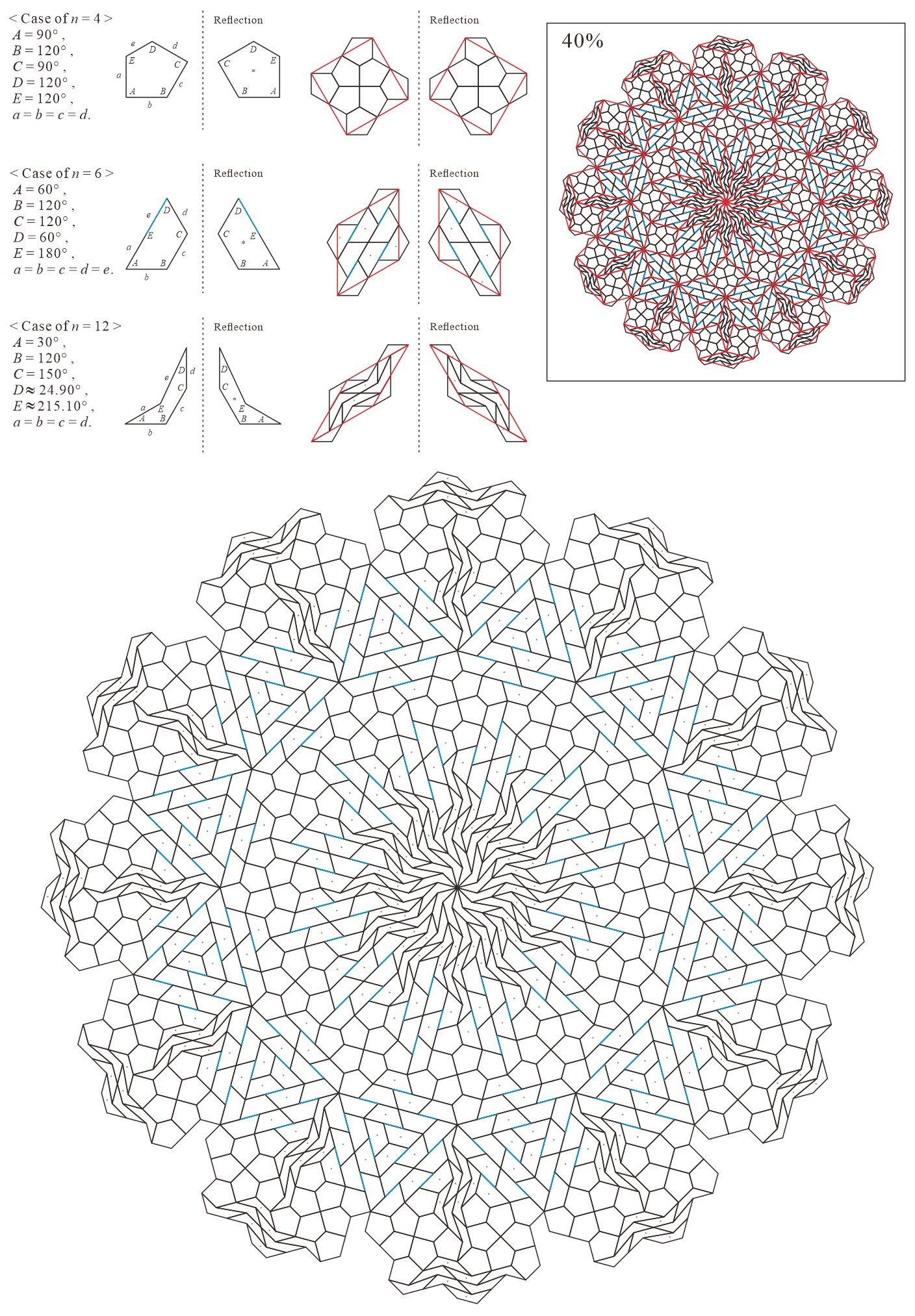} 
  \caption{{\small 
Tiling that generated based on 12-fold rotationally symmetric tiling, 
with three types of rhombuses applying conversion of rhombuses into eight 
pentagons satisfying (\ref{eq3}) with $n = 4, 6, 12$ and $\theta = 30^ \circ $ 
(The figure is solely a depiction of the area around the rotationally 
symmetric structure, and the tiling can be spread in all directions)
} 
\label{Fig.5.3.2-2}
}
\end{figure}

\renewcommand{\figurename}{{\small Figure.}}
\begin{figure}[p]
 \centering\includegraphics[width=15cm,clip]{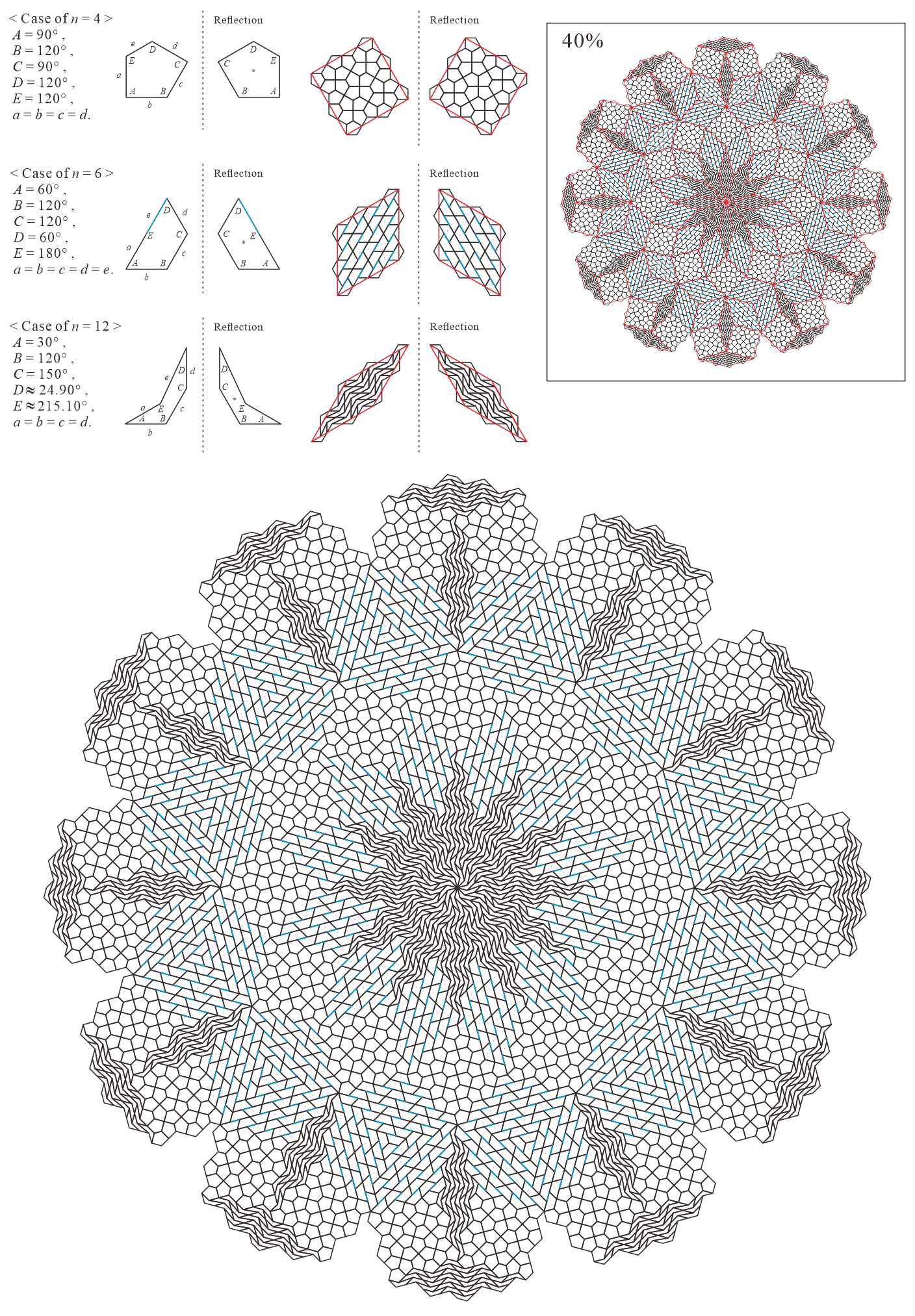} 
  \caption{{\small 
Tiling that generated based on 12-fold rotationally symmetric tiling, 
with three types of rhombuses applying conversion of rhombuses into 32 
pentagons satisfying (\ref{eq3}) with $n = 4, 6, 12$ and $\theta = 30^ \circ $ 
(The figure is solely a depiction of the area around the rotationally 
symmetric structure, and the tiling can be spread in all directions)
} 
\label{Fig.5.3.2-3}
}
\end{figure}


\renewcommand{\figurename}{{\small Figure.}}
\begin{figure}[p]
 \centering\includegraphics[width=15cm,clip]{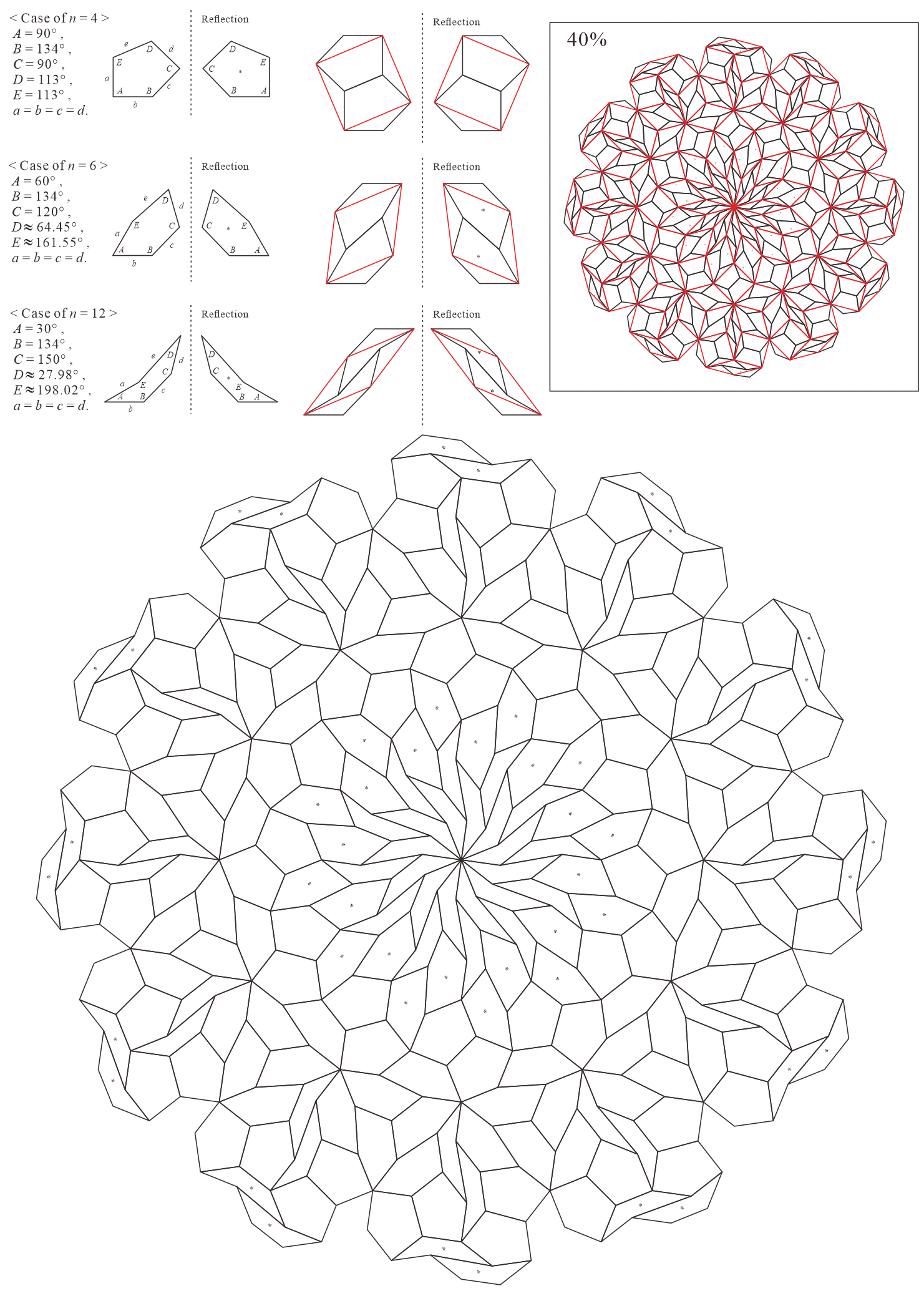} 
  \caption{{\small 
Tiling that generated based on 12-fold rotationally symmetric tiling, 
with three types of rhombuses applying conversion of rhombuses into two 
pentagons satisfying (\ref{eq3}) with $n = 4, 6, 12$ and $\theta = 44^ \circ $ 
(The figure is solely a depiction of the area around the rotationally 
symmetric structure, and the tiling can be spread in all directions)
} 
\label{Fig.5.3.3-1}
}
\end{figure}

\renewcommand{\figurename}{{\small Figure.}}
\begin{figure}[p]
 \centering\includegraphics[width=15cm,clip]{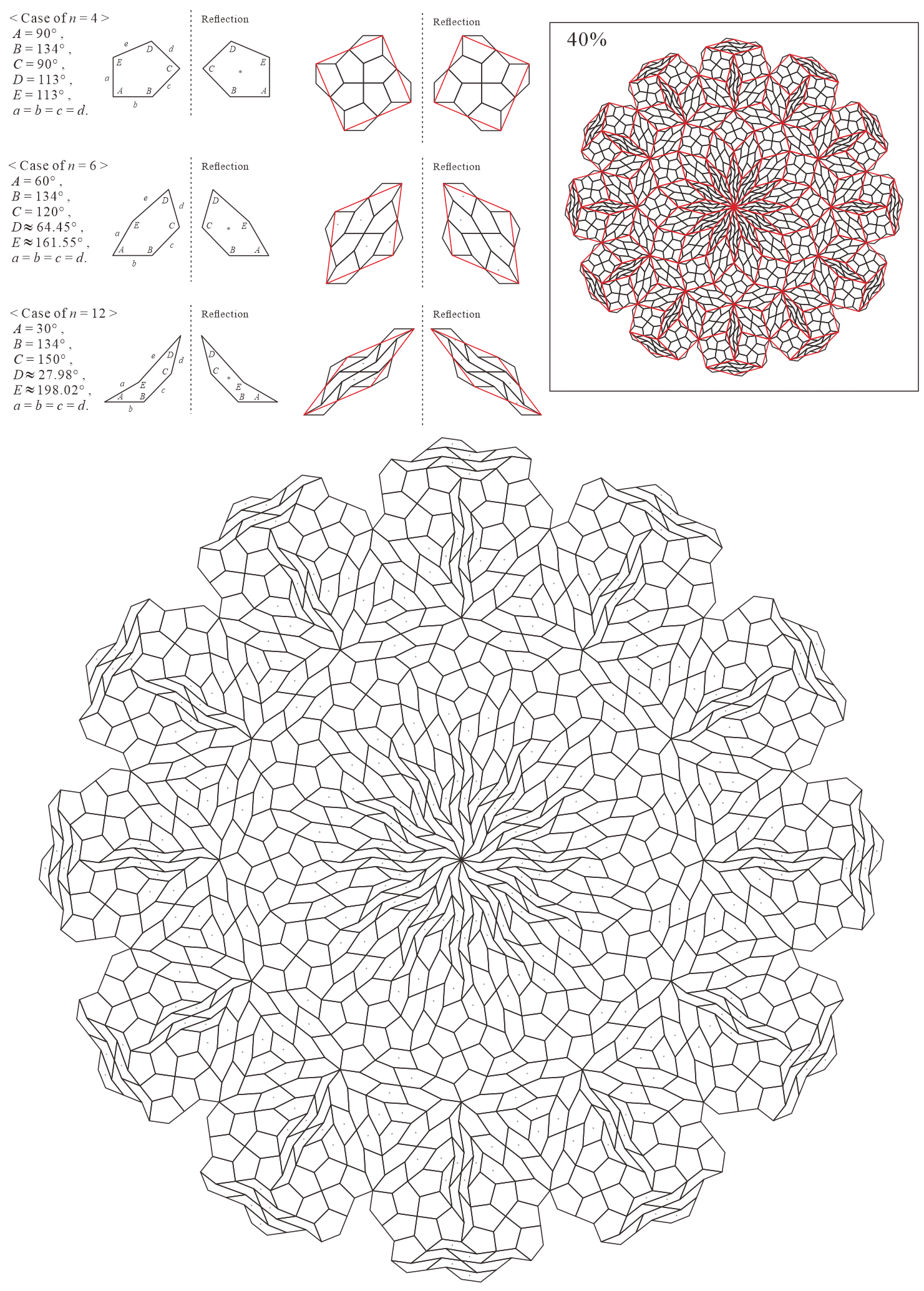} 
  \caption{{\small 
Tiling that generated based on 12-fold rotationally symmetric tiling, 
with three types of rhombuses applying conversion of rhombuses into eight 
pentagons satisfying (\ref{eq3}) with $n = 4, 6, 12$ and $\theta = 44^ \circ $ 
(The figure is solely a depiction of the area around the rotationally 
symmetric structure, and the tiling can be spread in all directions)
} 
\label{Fig.5.3.3-2}
}
\end{figure}


\renewcommand{\figurename}{{\small Figure.}}
\begin{figure}[p]
 \centering\includegraphics[width=15cm,clip]{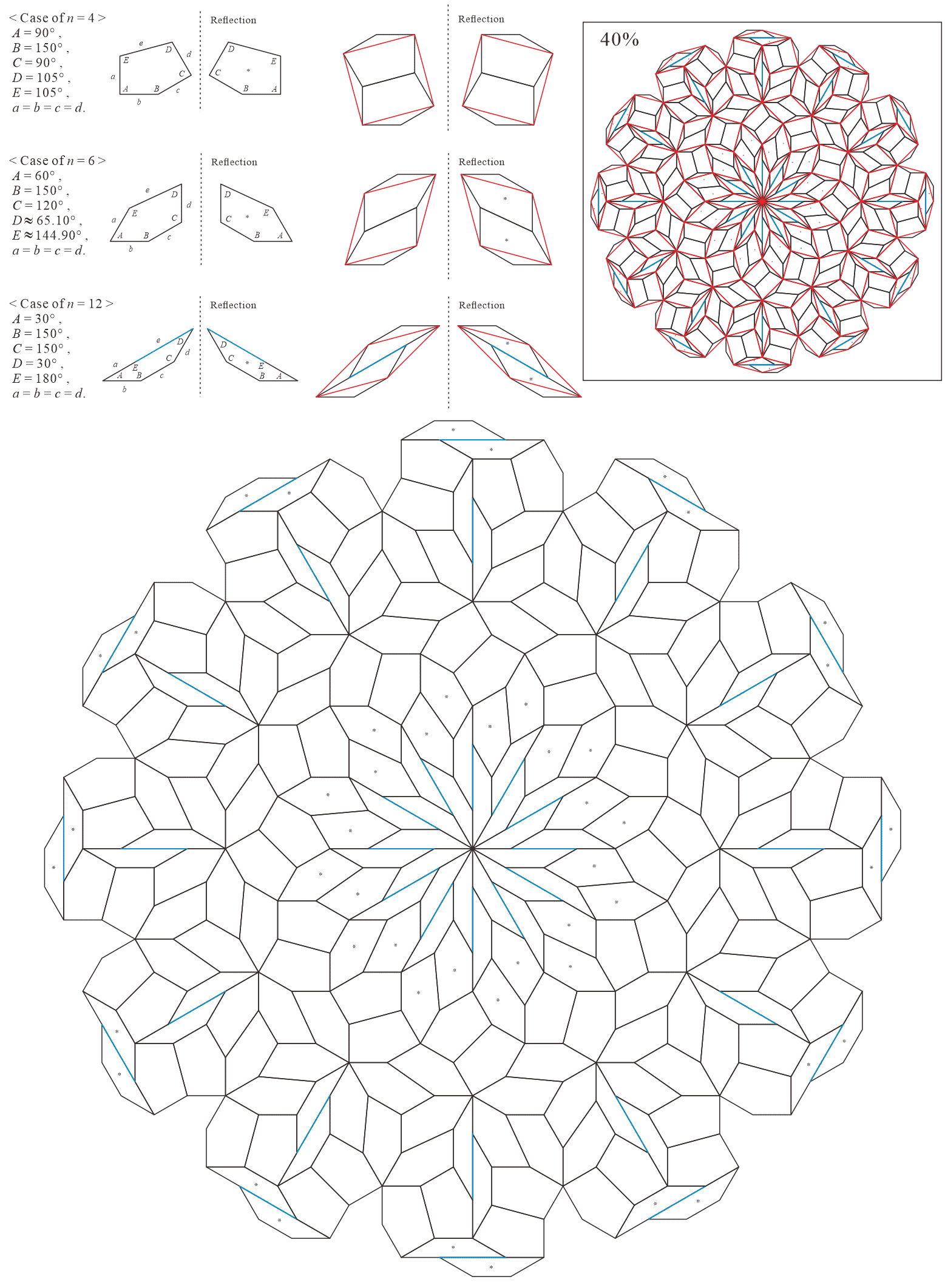} 
  \caption{{\small 
Tiling that generated based on 12-fold rotationally symmetric tiling, 
with three types of rhombuses applying conversion of rhombuses into two 
pentagons satisfying (\ref{eq3}) with $n = 4, 6, 12$ and $\theta = 60^ \circ $ 
(The figure is solely a depiction of the area around the rotationally 
symmetric structure, and the tiling can be spread in all directions)
} 
\label{Fig.5.3.4-1}
}
\end{figure}

\renewcommand{\figurename}{{\small Figure.}}
\begin{figure}[p]
 \centering\includegraphics[width=15cm,clip]{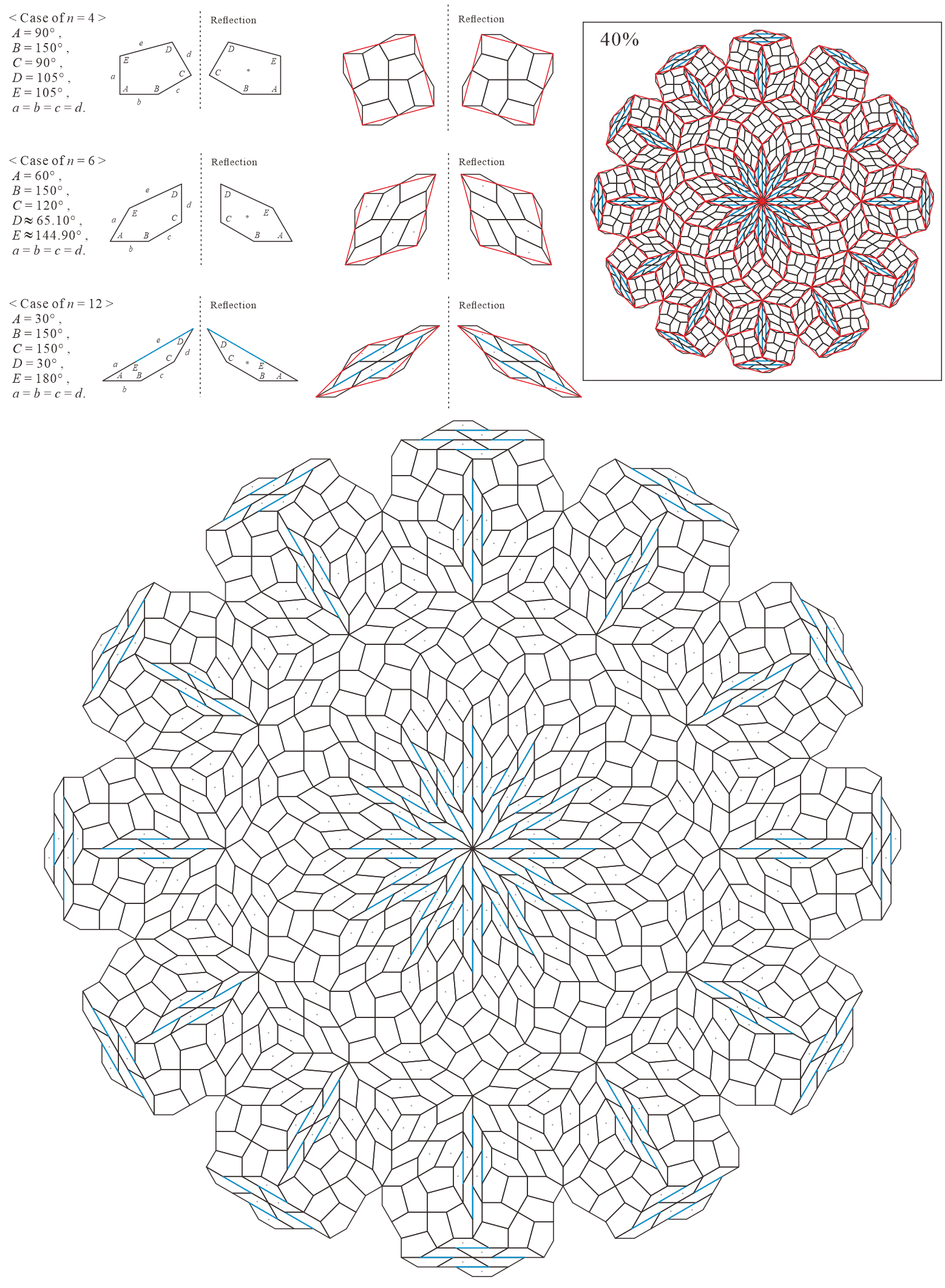} 
  \caption{{\small 
Tiling that generated based on 12-fold rotationally symmetric tiling, 
with three types of rhombuses applying conversion of rhombuses into eight 
pentagons satisfying (\ref{eq3}) with $n = 4, 6, 12$ and $\theta = 60^ \circ $ 
(The figure is solely a depiction of the area around the rotationally 
symmetric structure, and the tiling can be spread in all directions)
} 
\label{Fig.5.3.4-2}
}
\end{figure}


\renewcommand{\figurename}{{\small Figure.}}
\begin{figure}[p]
 \centering\includegraphics[width=15cm,clip]{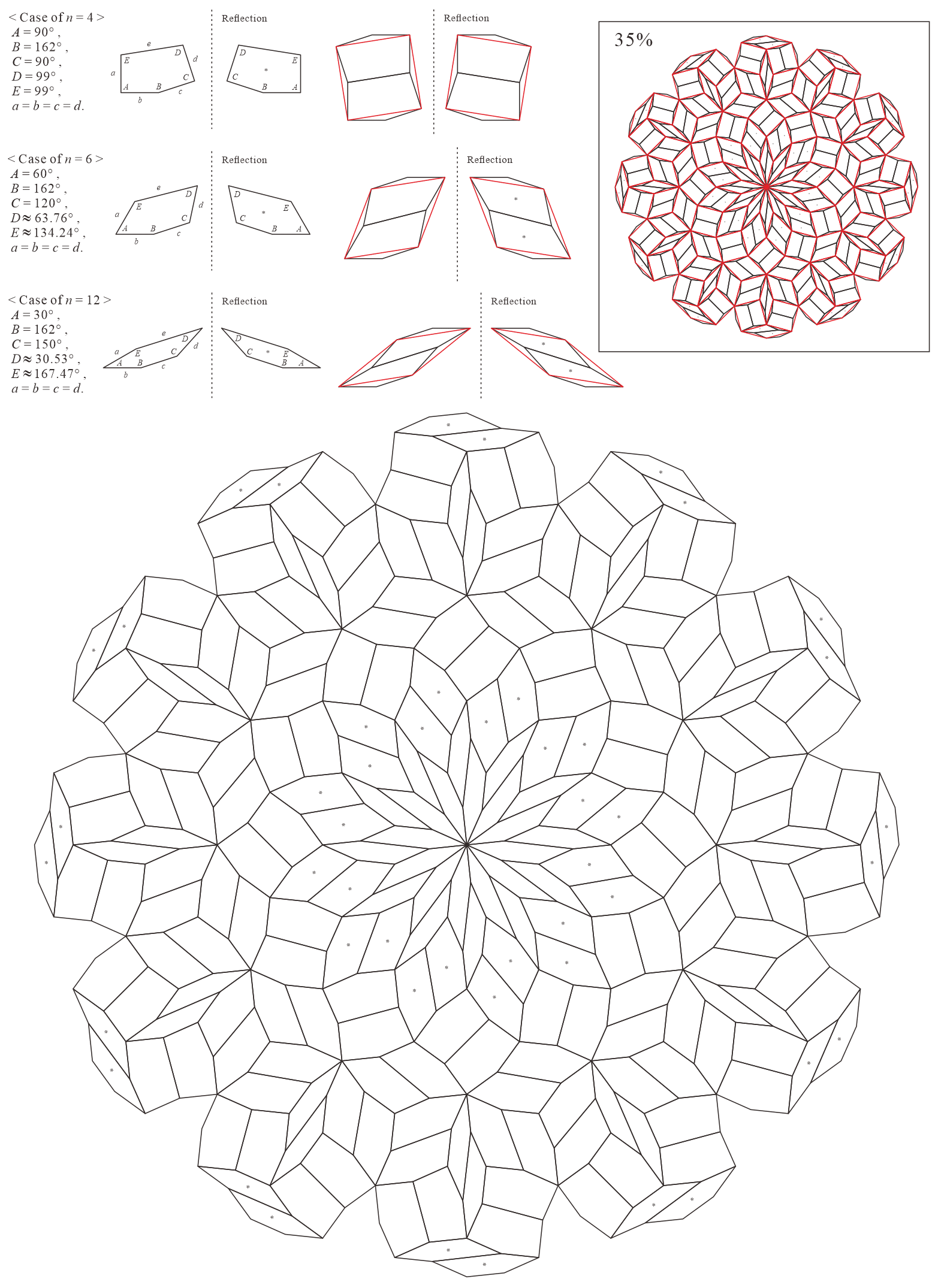} 
  \caption{{\small 
Tiling that generated based on 12-fold rotationally symmetric tiling, 
with three types of rhombuses applying conversion of rhombuses into two 
pentagons satisfying (\ref{eq3}) with $n = 4, 6, 12$ and $\theta = 72^ \circ $ 
(The figure is solely a depiction of the area around the rotationally 
symmetric structure, and the tiling can be spread in all directions)
} 
\label{Fig.5.3.5-1}
}
\end{figure}

\renewcommand{\figurename}{{\small Figure.}}
\begin{figure}[p]
 \centering\includegraphics[width=15cm,clip]{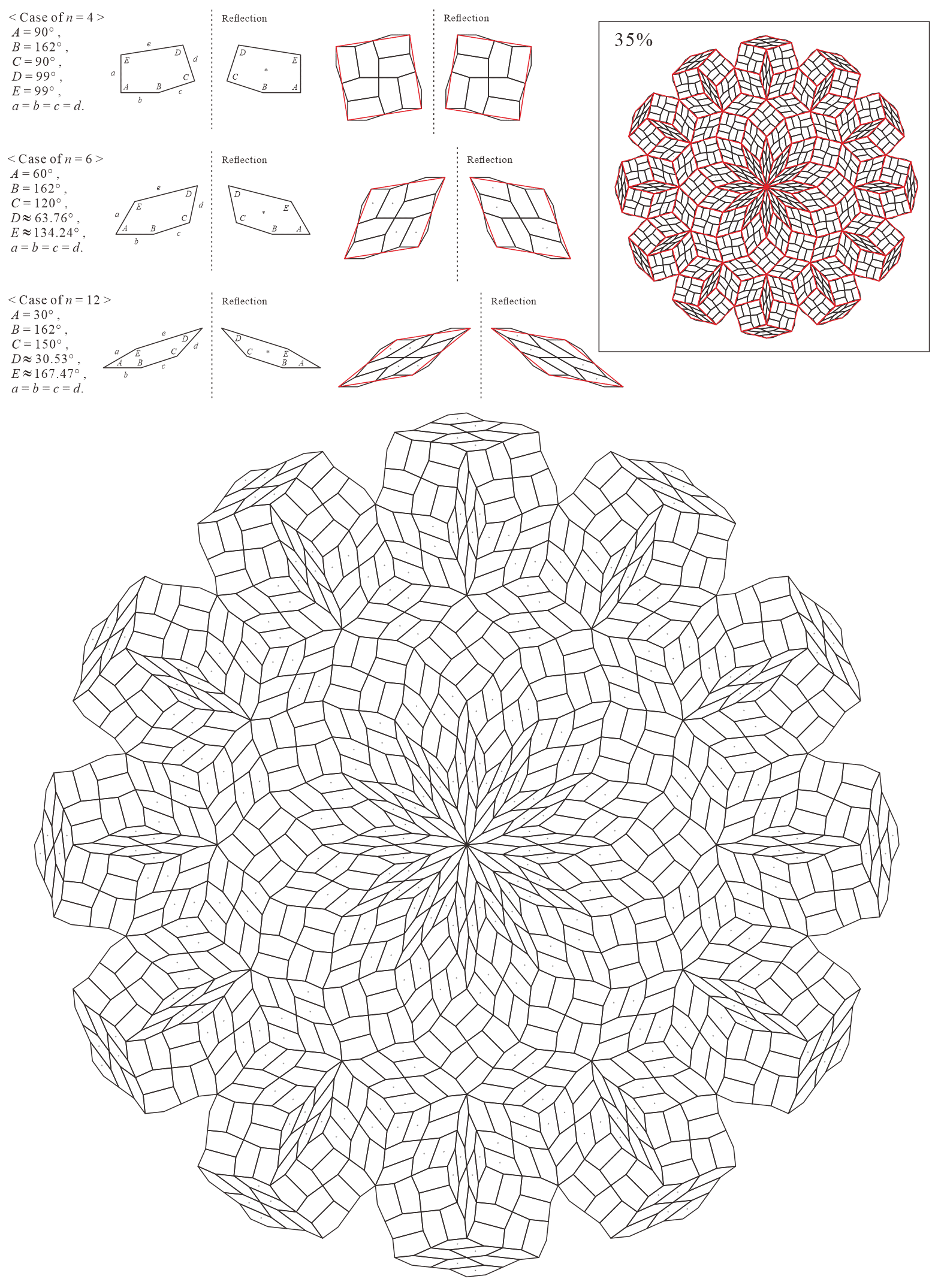} 
  \caption{{\small 
Tiling that generated based on 12-fold rotationally symmetric tiling, 
with three types of rhombuses applying conversion of rhombuses into eight 
pentagons satisfying (\ref{eq3}) with $n = 4, 6, 12$ and $\theta = 72^ \circ $ 
(The figure is solely a depiction of the area around the rotationally 
symmetric structure, and the tiling can be spread in all directions)
} 
\label{Fig.5.3.5-2}
}
\end{figure}

\renewcommand{\figurename}{{\small Figure.}}
\begin{figure}[p]
 \centering\includegraphics[width=15cm,clip]{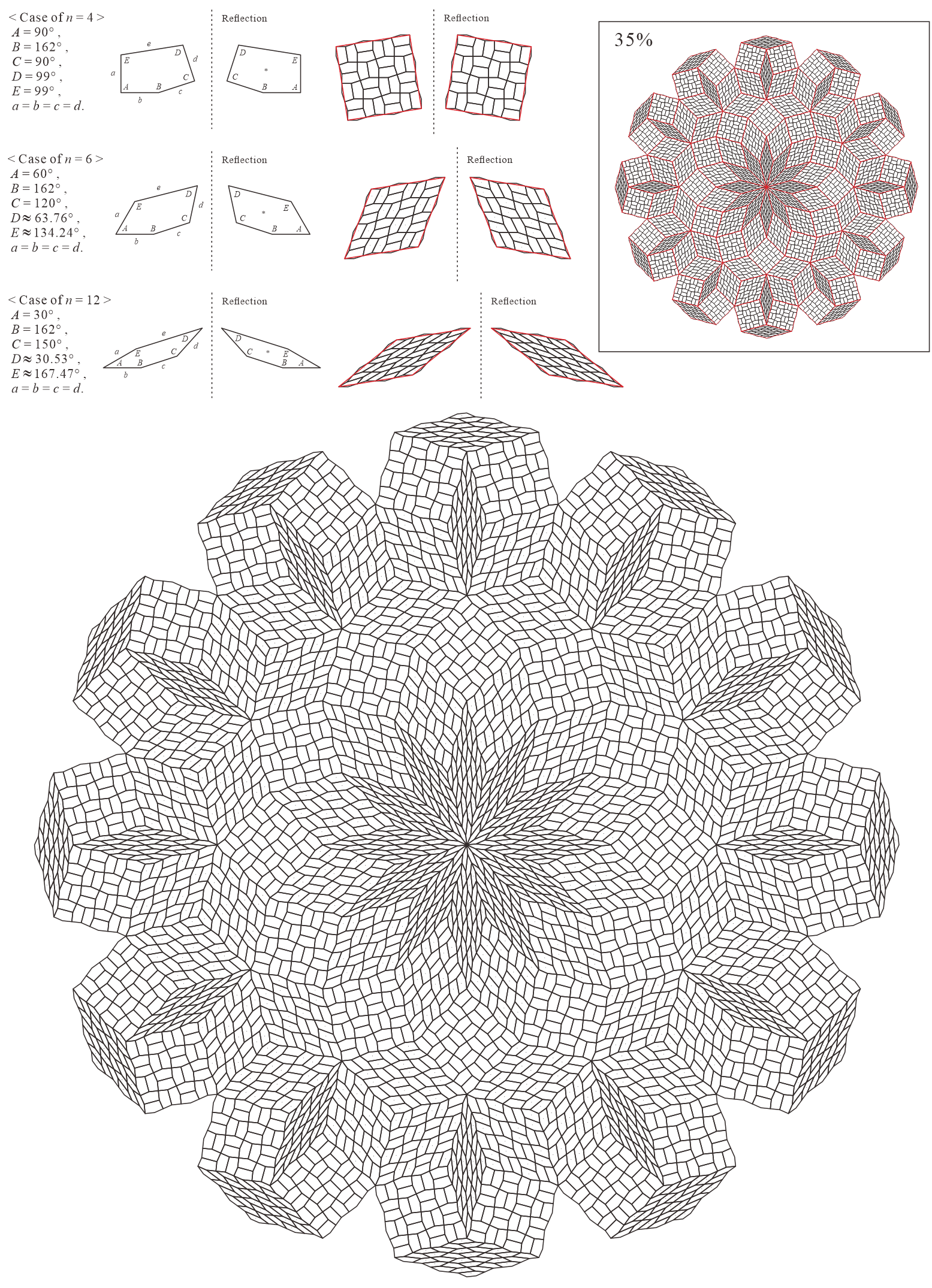} 
  \caption{{\small 
Tiling that generated based on 12-fold rotationally symmetric tiling, 
with three types of rhombuses applying conversion of rhombuses into 32 
pentagons satisfying (\ref{eq3}) with $n = 4, 6, 12$ and $\theta = 72^ \circ $ 
(The figure is solely a depiction of the area around the rotationally 
symmetric structure, and the tiling can be spread in all directions)
} 
\label{Fig.5.3.5-3}
}
\end{figure}


\renewcommand{\figurename}{{\small Figure.}}
\begin{figure}[p]
 \centering\includegraphics[width=15cm,clip]{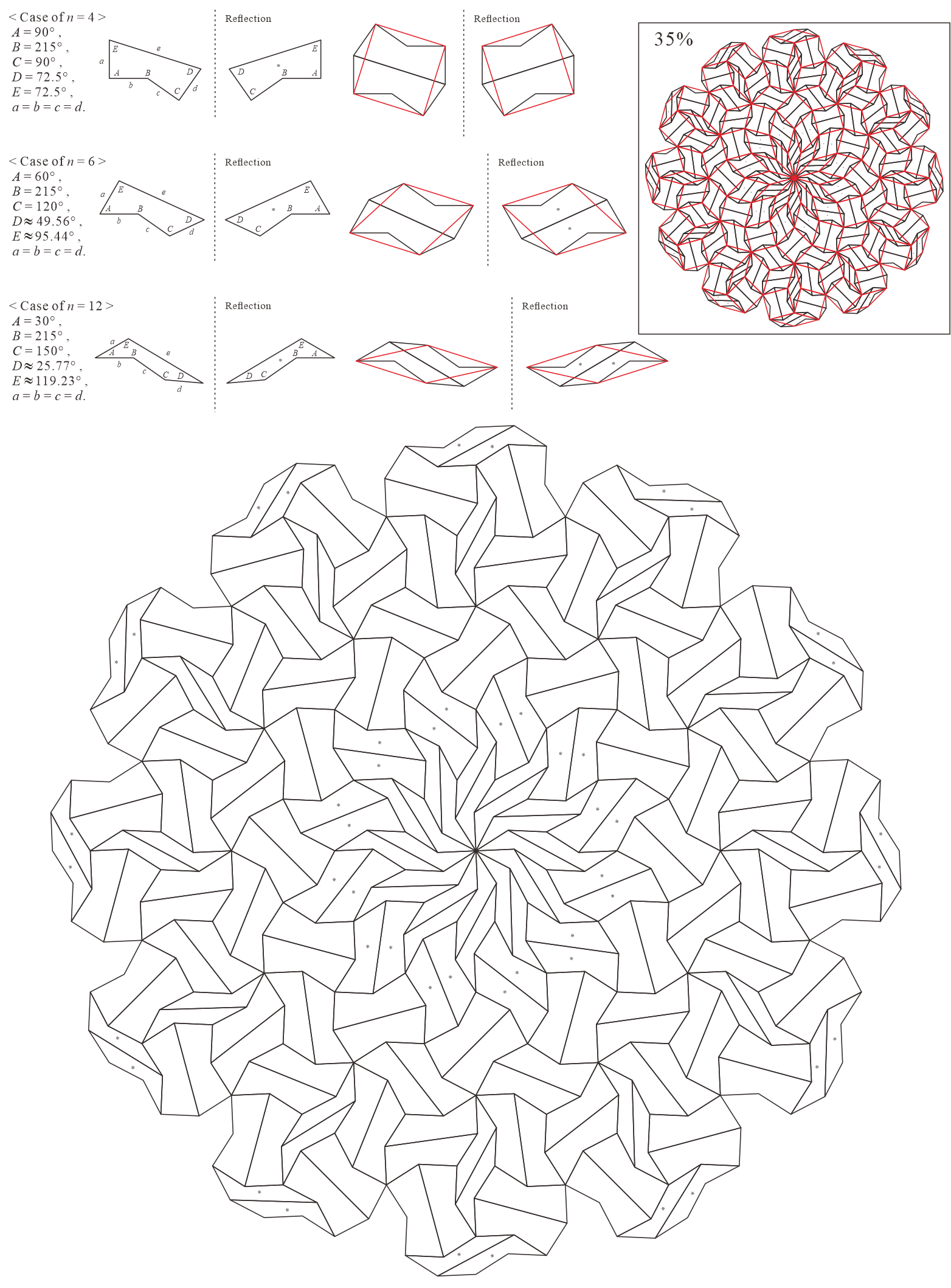} 
  \caption{{\small 
Tiling that generated based on 12-fold rotationally symmetric tiling, 
with three types of rhombuses applying conversion of rhombuses into two 
pentagons satisfying (\ref{eq3}) with $n = 4, 6, 12$ and $\theta = 125^ \circ $ 
(The figure is solely a depiction of the area around the rotationally 
symmetric structure, and the tiling can be spread in all directions)
} 
\label{Fig.5.3.7-1}
}
\end{figure}

\renewcommand{\figurename}{{\small Figure.}}
\begin{figure}[p]
 \centering\includegraphics[width=15cm,clip]{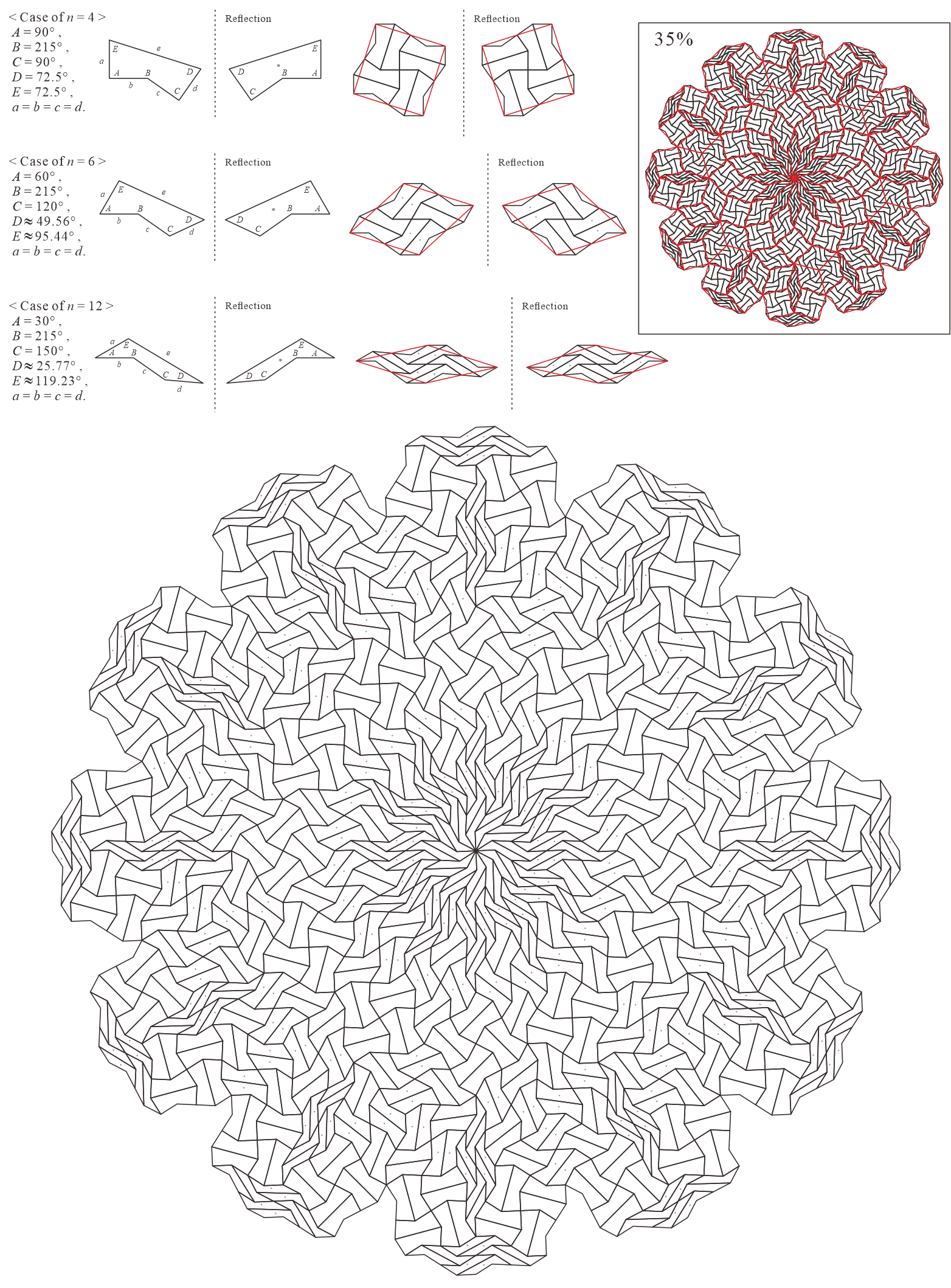} 
  \caption{{\small 
Tiling that generated based on 12-fold rotationally symmetric tiling, 
with three types of rhombuses applying conversion of rhombuses into eight 
pentagons satisfying (\ref{eq3}) with $n = 4, 6, 12$ and $\theta = 125^ \circ $ 
(The figure is solely a depiction of the area around the rotationally 
symmetric structure, and the tiling can be spread in all directions)
} 
\label{Fig.5.3.7-2}
}
\end{figure}


\renewcommand{\figurename}{{\small Figure.}}
\begin{figure}[p]
 \centering\includegraphics[width=15cm,clip]{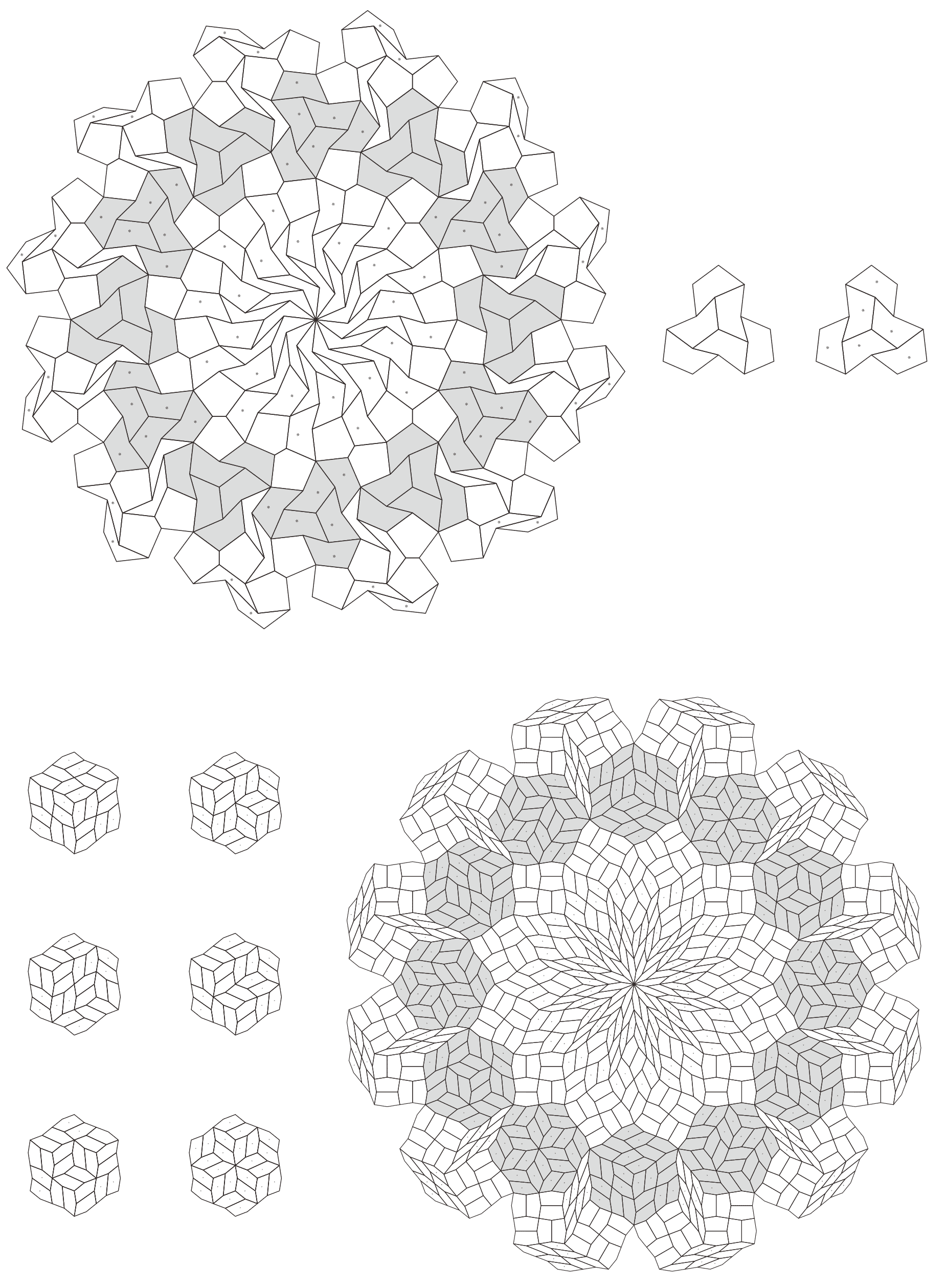} 
  \caption{{\small 
Examples of replacing tilings with different pattern using 
units with $D_{3}$ symmetry that can form pentagons satisfying (\ref{eq3}) with 
$n = 6$
} 
\label{Fig.5.3.8-1}
}
\end{figure}

\renewcommand{\figurename}{{\small Figure.}}
\begin{figure}[p]
 \centering\includegraphics[width=15cm,clip]{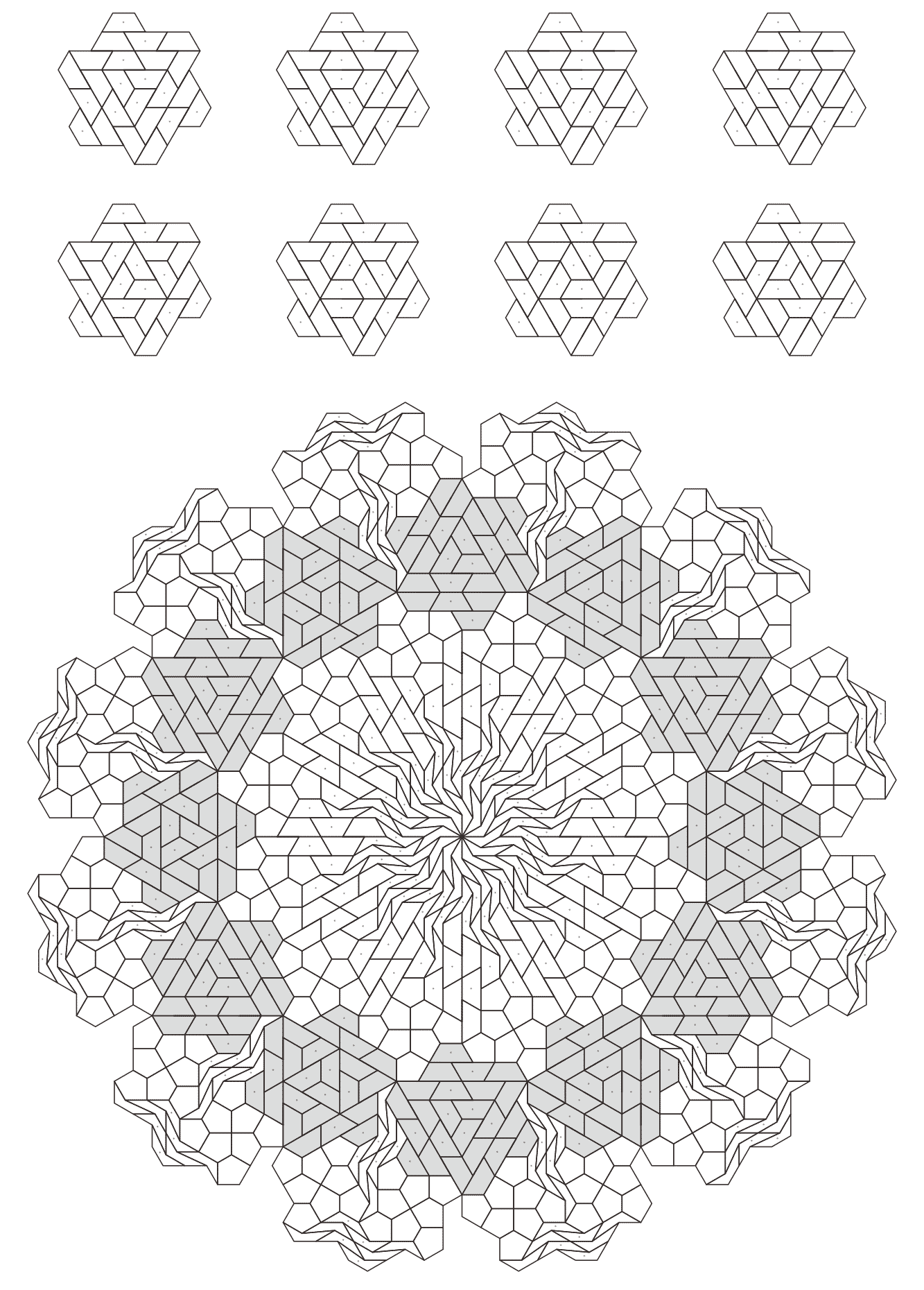} 
  \caption{{\small 
Example of replacing tiling with different pattern using 
equilateral triangles that can form pentagons satisfying (\ref{eq3}) with 
$n = 6$ and $\theta = 30^ \circ $
} 
\label{Fig.5.3.8-2}
}
\end{figure}

\end{document}